\RequirePackage{fix-cm}
\documentclass[smallextended]{svjour3}       
\smartqed  
\usepackage{graphicx}
\usepackage[top=2cm, bottom=2.5cm, left=3cm, right=3cm]{geometry}

\usepackage[ngerman, english]{babel}   

\usepackage[latin1]{inputenc}
\usepackage{graphicx}                      
\usepackage{lmodern}                        
\usepackage{caption}
\usepackage{tikz}
 \usetikzlibrary{calc}                                                                                                                        
  \usetikzlibrary{decorations.markings}
 \usetikzlibrary{arrows}
\usetikzlibrary{arrows.meta}
\usepackage{xcolor}
\usetikzlibrary{fit}
\usepackage{subcaption}
\usepackage{mathtools}
\usepackage{graphicx}
\newlength{\imagewidth}
\usepackage{float} 

\usepackage{booktabs}
\usepackage{textcomp}
\usepackage[all,cmtip]{xy} 
\usepackage{amsmath}
\usepackage{lscape}
\usepackage{amssymb,amsmath} 
\usepackage{lastpage} 
\usepackage{braket}

\usepackage{setspace}
\usepackage{nicefrac}
\usepackage{siunitx}
\usepackage{footnote}

\allowdisplaybreaks

\DeclareMathOperator{\id}{id}

\DeclareMathOperator{\Iso}{Iso}
\DeclareMathOperator{\p}{\mathfrak{p}}

\DeclareMathOperator{\Aut}{Aut}

\DeclareMathOperator{\Hom}{Hom}
\DeclareMathOperator{\Inn}{Inn}
\DeclareMathOperator{\Out}{Out}
\DeclareMathOperator{\Mod}{Mod}
\DeclareMathOperator{\PMod}{PMod}

\newcommand{\orb}[1]{\normalsize{\ensuremath{\mathsf{#1}}}\normalsize}

\usepackage{url}
\usepackage{hyperref}
\usepackage{amssymb,amsmath}
\usepackage{caption}
\usepackage{subcaption}
\begin{document}

\title{Enumerating Isotopy Classes of Tilings guided by the symmetry of Triply-Periodic Minimal Surfaces\thanks{This research was funded by the Emmy Noether Programme of the Deutsche Forschungsgemeinschaft. B.K was supported by the Deutscher Akademischer Austauschdienst for a research stay at the Australian National University.}
}


\author{Benedikt Kolbe         \and
        Myfanwy E. Evans 
}


\institute{Benedikt Kolbe \at
Université de Lorraine, CNRS, Inria, LORIA, F-54000 Nancy, France \\
				\email{benedikt.kolbe@inria.fr}
           \and
           Myfanwy E. Evans \at
              Universit\"{a}t Potsdam\\
                Karl-Liebknecht-Str. 24-25\\
		14476 Potsdam-Golm\\
	\email{evans@uni-potsdam.de}
}

\maketitle

\begin{abstract}
We present a technique for the enumeration of all isotopically distinct ways of tiling a hyperbolic surface of finite genus, possibly nonorientable and with punctures and boundary. This generalizes the enumeration using Delaney-Dress combinatorial tiling theory of combinatorial classes of tilings to isotopy classes of tilings. To accomplish this, we derive an action of the mapping class group of the orbifold associated to the symmetry group of a tiling on the set of tilings. We explicitly give descriptions and presentations of semi-pure mapping class groups and of tilings as decorations on orbifolds. We apply this enumerative result to generate an array of isotopically distinct tilings of the hyperbolic plane with symmetries generated by rotations that are commensurate with the three-dimensional symmetries of the primitive, diamond and gyroid triply-periodic minimal surfaces, which have relevance to a variety of physical systems.
\end{abstract}

\begin{keywords}
Isotopic tiling theory, mapping class group, orbifolds, triply-periodic minimal surface, Delaney-Dress tiling theory, Hyperbolic tilings
\end{keywords}

\subclass{
  05B45, 05C30, 52C20, 57M07, 68U05, 82D25
}

\section{Introduction}
Hyperbolic tilings have found a spectacular niche in describing complicated three-dimensional structure in both structural chemistry and polymer self-assembly~\cite{Hyde1991,Hyde1993,Chen2001,Kirkensgaard2014}. These disparate systems are united by a common geometric thread of symmetric tilings on some particular hyperbolic triply-periodic minimal surfaces (TPMS), minimal surfaces embedded in $\mathbb{R}^3$ invariant under three linearly independent translations. This provides motivation for a geometric exploration of related structures, in particular, the enumeration of possible tilings of the gyroid, the diamond, the primitive, and other genus $3$ TPMS. 

\begin{figure}[!htbp]
\imagewidth=0.31\textwidth
\captionsetup[subfigure]{width=0.9\imagewidth}
  \begin{subfigure}[t]{0.3\textwidth}
  \centering
    \includegraphics[width=\textwidth]{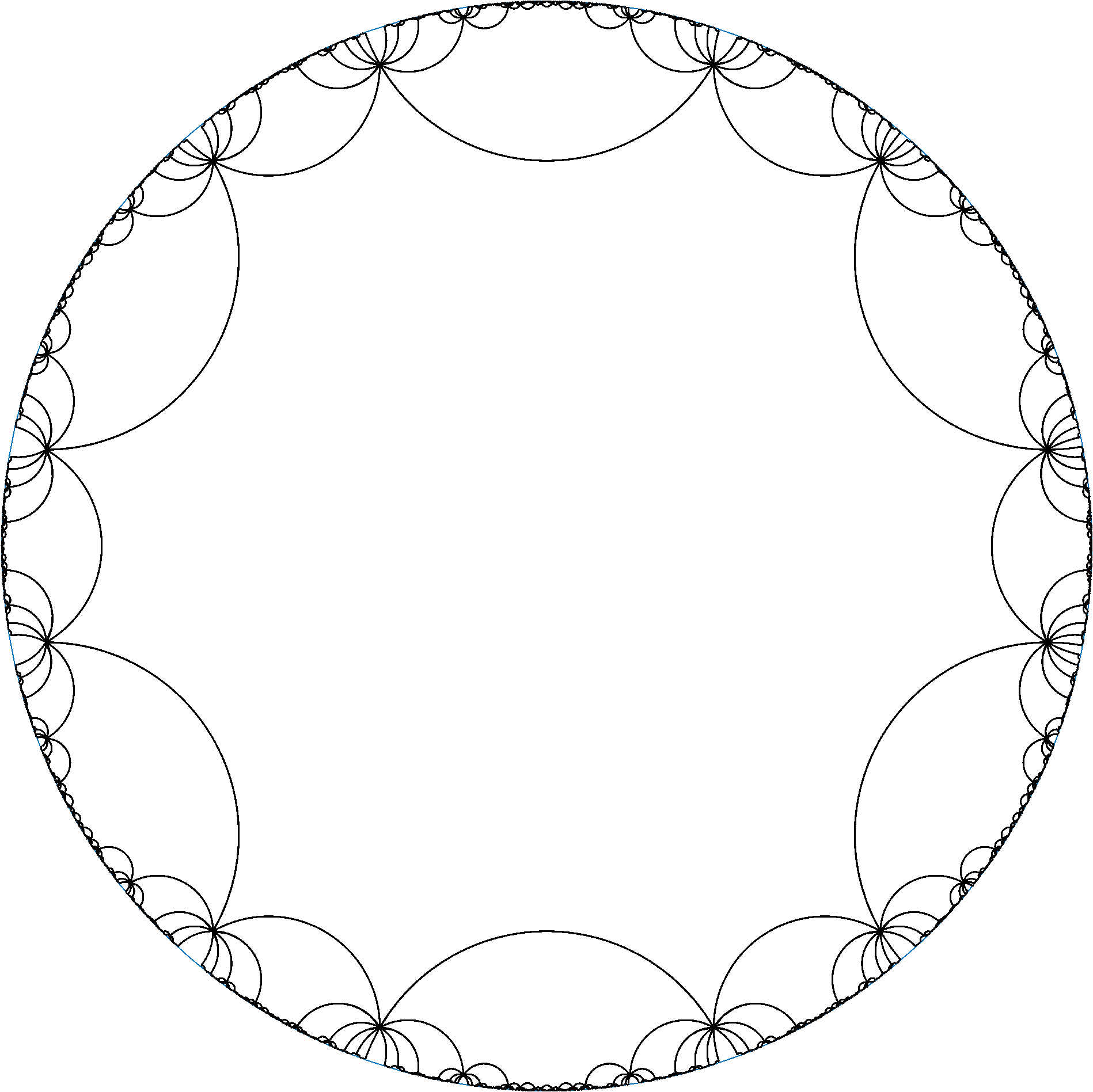}
  \caption{} \label{subfig:dodectess}
  \end{subfigure}
\hfill
  \begin{subfigure}[t]{0.3\textwidth}
  \centering
    \includegraphics[width=\textwidth]{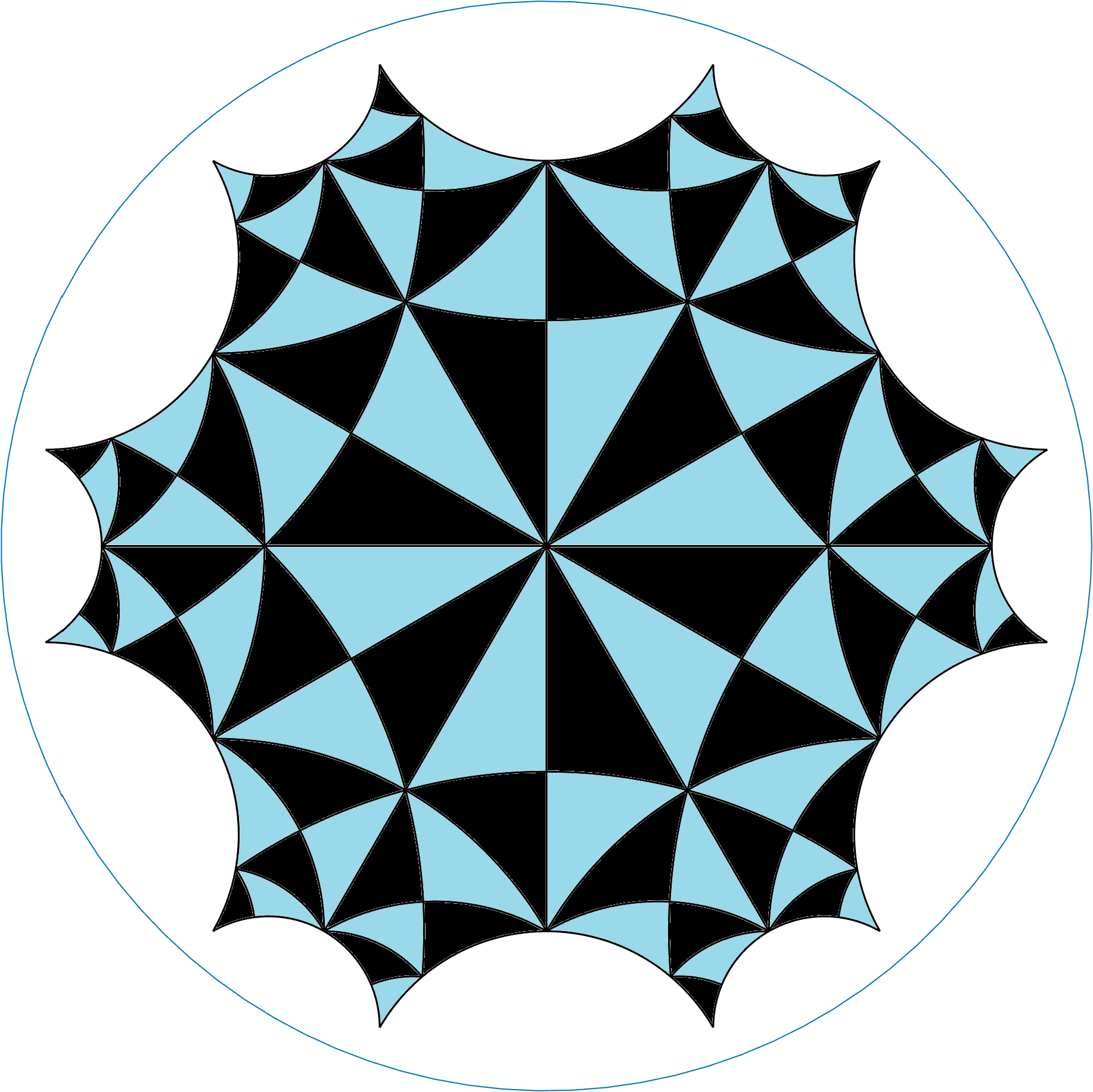}
 \caption{} \label{subfig:dodec246}
  \end{subfigure}
    \hfill
  \begin{subfigure}[t]{0.3\textwidth}
  \centering
    \includegraphics[width=\textwidth]{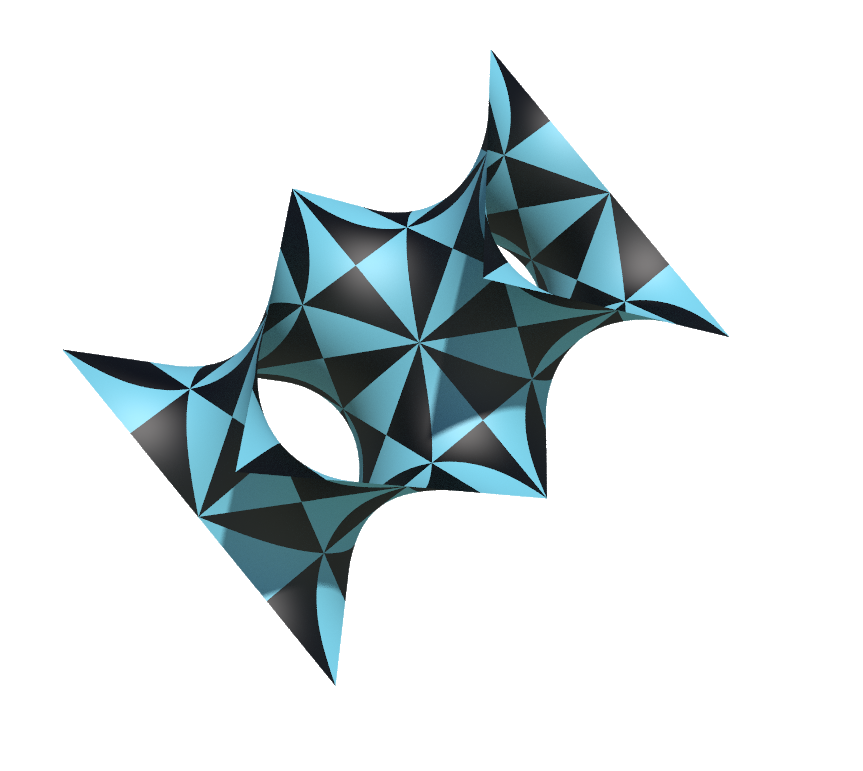}
 \caption{} \label{subfig:dsurf}
  \end{subfigure}
  \caption{The symmetries of the diamond triply-periodic minimal surface (D-surface) in $\mathbb{R}^3$ and its uniformization in $\mathbb{H}^2$. (a) Tesselation of $\mathbb{H}^2$ by dodecagons corresponding to the genus $3$ hyperbolic surface (after identifying opposite edges) that gives rise to the D-surface. (b) The tiling of $\mathbb{H}^2$ by triangles with symmetry group $\star 246.$, which is the symmetry of the D-surface. Each line represents a mirror symmetry. (c) A periodic unit cell of the D-surface in $\mathbb{R}^3$, together with its smallest asymmetric triangle patches. These correspond to the $\star 246$ tiles in $\mathbb{H}^2$ in the neighbouring image.}\label{fig:symmetries}
\end{figure}

\begin{figure}[!htbp]
\imagewidth=0.31\textwidth
\captionsetup[subfigure]{width=0.9\imagewidth}
  \begin{subfigure}[t]{0.33\textwidth}
  \centering
    \includegraphics[width=0.9\textwidth, height=0.9\textwidth]{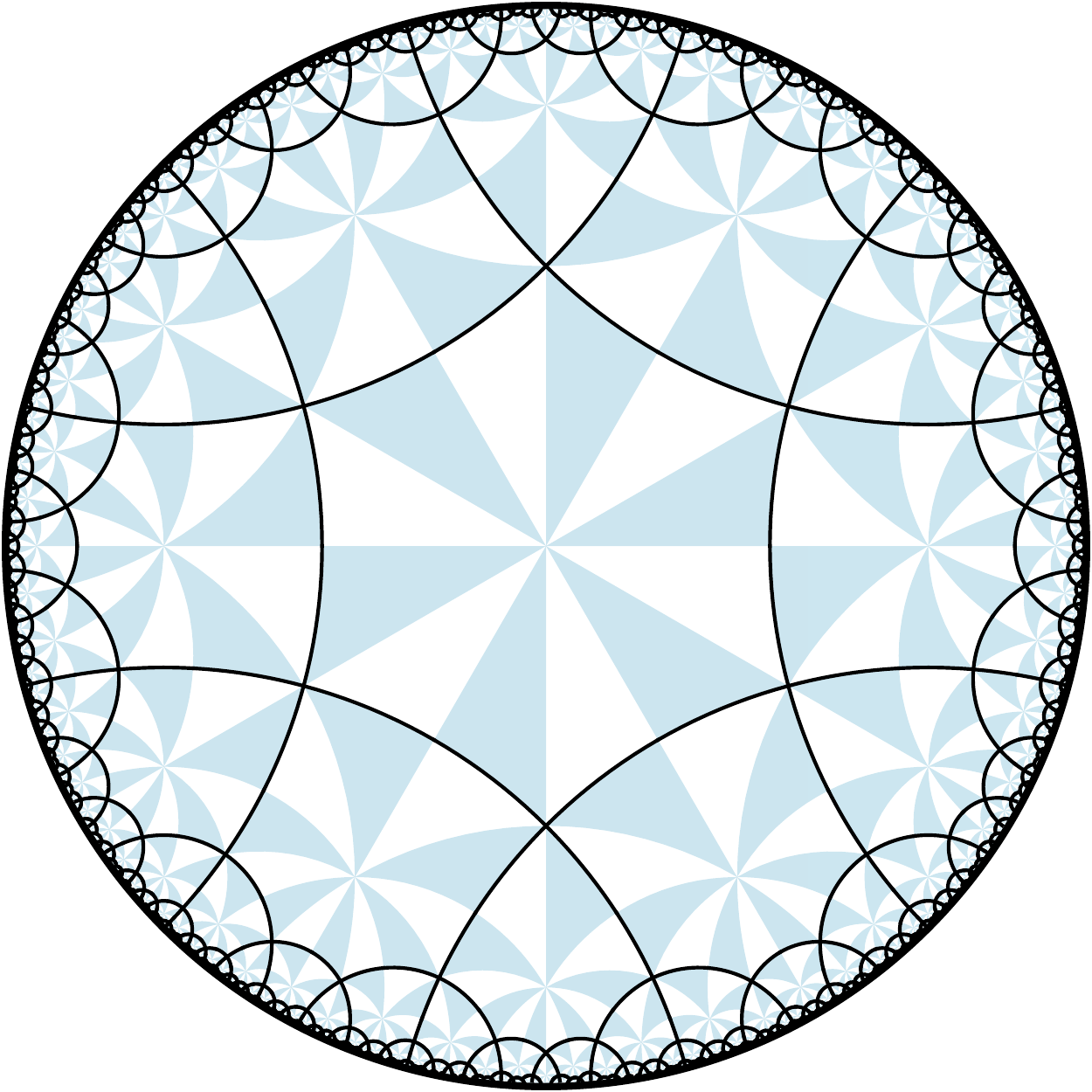}
    \caption{} 
  \end{subfigure}
  \begin{subfigure}[t]{0.32\textwidth}
  \centering
    \includegraphics[width=\textwidth]{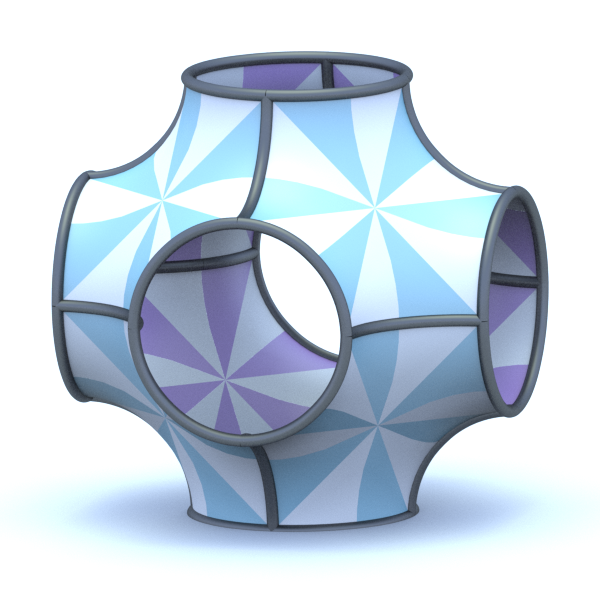}
    \caption{} 
  \end{subfigure}
  \begin{subfigure}[t]{0.32\textwidth}
  \centering
    \includegraphics[width=\textwidth,height=0.9\textwidth]{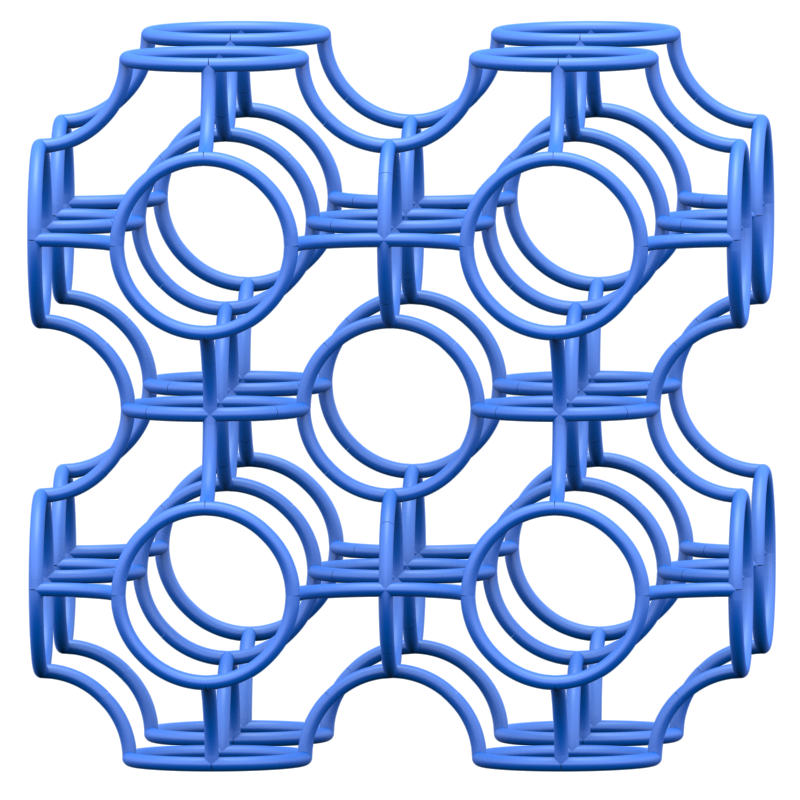}
    \caption{} 
  \end{subfigure}
  \caption{The cubic structure of the mineral Sodalite (c) can be described by a hyperbolic tiling projected to the primitive cubic triply-periodic minimal surface (P-surface). (a) A tiling of the hyperbolic plane with symmetry group $*246$, represented by the solid black lines, with a tiling by triangles with $\star 246$ symmetry shown in blue in the background. (b) The tiling from (a) shown as a decoration of the P-surface. The triangles illustrate the symmetries of the surface. (c) The resulting net in $\mathbb{R}^3$ when the tile boundaries are considered as curves in 3-dimensional Euclidean space rather than curves on the surface, which is the structure of Sodalite.}\label{fig:epinetidea}
\end{figure}

An extensive enumerative technique, called the EPINET project~\cite{epinet}, uses TPMS decorated with hyperbolic tilings as a blueprint for complicated structures, including crystallographic nets, in $\mathbb{R}^3$~\cite{Sadoc1989,Nesper2001,Robins2005,Robins2006,Ramsden2009,Castle2012}. Upon discarding the surface and only retaining the boundary of the tile edges, one obtains a net in $\mathbb{R}^3$. The hyperbolic in-surface symmetries of TPMS manifest as ambient Euclidean symmetries of $\mathbb{R}^3$~\cite{Perez2002} (see Figure~\ref{fig:symmetries}), so that symmetric tilings of TPMS give rise to symmetric graph embeddings in $\mathbb{R}^3.$  These ideas have been used to enumerate crystalline structures resulting from bounded hyperbolic tilings with symmetry group generated entirely by reflections~\cite{Ramsden2003,Ramsden2009}, also known as a \emph{Coxeter} group. An example of the construction process is shown in Figure~\ref{fig:epinetidea}.

\begin{figure}[!htbp]
\centering
\includegraphics[width=0.32\linewidth]{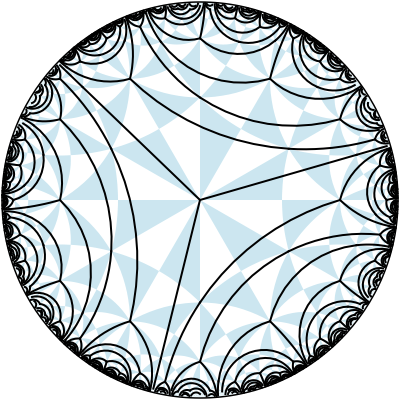}
\includegraphics[width=0.32\linewidth]{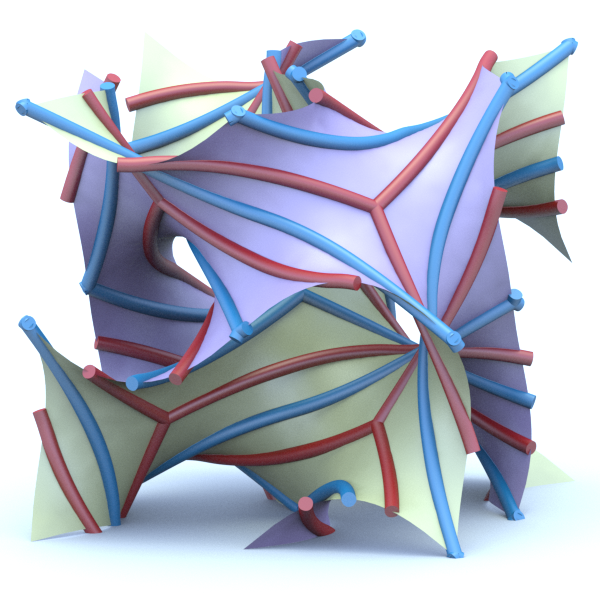}
\includegraphics[width=0.32\linewidth]{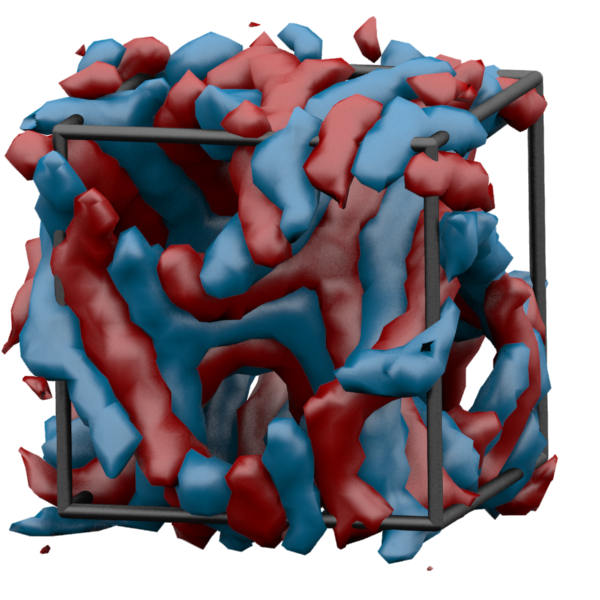}
\caption{A hyperbolic tiling in $\mathbb{H}^2$ with unbounded tiles and purely rotational symmetry (orbifold $2223$). The tiling is then shown in the gyroid triply-periodic minimal surface, which describes the domains formed in a numerical simulation of polymer self-assembly, shown on the right~\cite{Kirkensgaard2014}.}
\label{pnas}
\end{figure}

Beyond this enumeration, various sets of structures with purely rotational symmetry and various tile types have been explored in different contexts~\cite{Hyde2000,hyderams,evansper1,evansper2,evansper3,Pedersen2017}. A particular example relates to the complicated geometry of a self-assembled block copolymer system that has been recently characterised~\cite{Kirkensgaard2014}. The mutual repulsion of different polymer chain species in the star-shaped molecule causes an arrangement into complicated domains in $\mathbb{R}^3$. These domains can be described by hyperbolic tilings composed of unbounded tiles with a network-like boundary and rotational symmetry, decorating the gyroid TPMS~\cite{evansper1}, see Figure~\ref{pnas}. These mesostructures are among the most topologically complex morphologies identified in polymer sciences to date. More recently, hyperbolic honeycombs with rotational symmetry reticulated over TPMS have been used to construct new infinite deltahedra~\cite{Pedersen2018}.

The varied and physically relevant structures that arise from moving away from Coxeter symmetry groups of the hyperbolic tilings motivate the extension of the EPINET enumeration to a broader range of symmetries. So far, all approaches to enumerations for tilings with symmetry group not generated entirely by reflections have involved ad-hoc ideas for specific groups and have only met with success in a small number of restricted cases. Here, we develop a general framework to systematically enumerate isotopically distinct hyperbolic tilings with arbitrary symmetry groups of a given hyperbolic surface, which we consider as a finite-area Riemannian surface locally isometric to the hyperbolic plane. We specifically apply this enumerative technique to construct isotopically distinct hyperbolic tilings with purely rotational symmetry as candidates for reticulation over the candidate TPMS.


Delaney-Dress combinatorial tiling theory~\cite{DRESS1987,Huson1993,Delgado-Friedrichstilings} is an essential tool of the enumerative process. It deals with the classification of combinatorial classes of equivariant tilings, i.e. tilings with a specified symmetry group, of simply connected spaces, in which every tile is a bounded disk. Our recent generalization of combinatorial tiling theory allows, in theory, the enumeration of isotopy classes of tilings on hyperbolic surfaces~\cite{BenMyf1}. This generalization uses the notion of orbifolds and mapping class groups (MCGs), and interprets tilings as decorations, or piece-wise linear embeddings of graphs, on orbifolds, with Delaney-Dress (D-)symbols representing triangulations of orbifolds.

In this paper, we use this generalization to derive a methodology to completely enumerate all isotopy classes of equivariant tilings with symmetry group commensurate with some fixed hyperbolic surface, meaning that the symmetry group is isomorphic to a subgroup of the isometry group of that surface. In practice, our enumerative approach requires the conversion of the above theoretical framework (from \cite{BenMyf1}) into a practical setting, making use of topological ideas, computational group theory, braid theory, and presentations of MCGs. We describe this non-trivial implementation in detail, alongside limitations. Applying this result, we will show how to enumerate isotopy classes of equivariant tilings with a \emph{stellate} symmetry group, generated entirely by rotations (we also call any orbifold with a symmetry group generated entirely by rotations a \emph{stellate} orbifold). We focus on the enumeration of isotopy classes of tilings that are \emph{commensurate} with a genus-$3$ hyperbolic surface with symmetry group equal to $\star 246$ (see Section~\ref{sec:orbs} for notation), whereby we mean that the edge graph of the hyperbolic tiling projects to an edge graph of the hyperbolic surface homomorphically. This surface plays the role of the uniformized version of the surface within a unit cell for the gyroid, diamond, and primitive TPMS, so the enumerated tilings give rise to physically relevant structures of interest to the natural sciences.

This paper is structured into four sections which cumulatively build towards the enumeration process of isotopy classes of tilings, culminating in the application to explicit tilings of stellate symmetry. Section~\ref{sec:prelims} will be a recollection of orbifolds, combinatorial tiling theory, the framework from \cite{BenMyf1} for isotopic tiling theory as well as fundamental considerations for our enumeration. In section~\ref{sec:conciso}, we relate the geometry of orbifolds to the sets of generators of its symmetry group. The statement of our results depend heavily on the orbifold under consideration, so to avoid a lengthy presentation with many special cases, we instead separated the many different situations into different subsections. Section~\ref{sec:mixedbraids} then uses the connection of the semi-pure braid group on the sphere to the MCG to find appropriate presentations of MCGs. We subsequently introduce a natural data structure for tilings and move to an enumerative setting for tilings with stellate symmetry, and conclude by presenting a collection of examples of tilings in section~\ref{sec:steltileex}, illustrating the approach. 

\section{Preliminaries}\label{sec:prelims}
\subsection{The standard presentation of orbifolds}\label{sec:orbs}

Throughout this paper, let $\Gamma$ be a discrete group of isometries of $\mathbb{H}^2$. For simplicity, we restrict to the case where $\Gamma$ admits a compact fundamental domain, in which case $\Gamma$ is an \emph{NEC} (non-Euclidean crystallographic) group. Note that our results retain their validity in more generality for groups with a fundamental domain of finite-area, as we point out along the way. We identify the isomorphism class of an NEC group using Conway's \emph{orbifold symbol}~\cite{Conway2002,Macbeath1967}.  

\begin{definition} 
Let $\Gamma$ be an NEC group. A (hyperbolic) orbifold, $\mathcal{O}$, is topologically the quotient space $\mathbb{H}^2/\Gamma$ obtained by identifying points of $\mathbb{H}^2$ under the action of $\Gamma$. The orbifold structure retains the metric information carried by the particular isometries of $\Gamma$, by keeping track of the types of branching of the canonical projection map $p:\mathbb{H}^2\to \mathbb{H}^2/\Gamma,$ and an atlas of charts compatible with the $\Gamma$ action. The group $\Gamma$ is also known as the orbifold's symmetry group. 
\end{definition}
For more detailed definitions of the concepts involved, refer to \cite{ratcliffe2006foundations}. 
We distinguish between the underlying topological space and an orbifold with the additional structure by denoting them with $O$ and $\mathcal{O},$ respectively. We call $\mathcal{O}$ orientable if $O$ is orientable, and the nowhere dense set of branch points in $\mathbb{H}^2/\Gamma$ of the covering $p$ the \emph{singular locus} of the orbifold $\mathcal{O}.$  We refer to the isolated points of the singular locus as \emph{marked points}. 

Note that for more general finite-area orbifolds, cusps of the action of $\Gamma$ correspond to punctures in $O$ and can be treated as marked points. Moreover, for orbifolds with boundary components in $O$, the universal covering space $\tilde{U}$ is not $\mathbb{H}^2$ but a totally geodesic subspace thereof, i.e. where all geodesics in $\tilde{U}$ are geodesics in $\mathbb{H}^2$~\cite{katok1992fuchsian,Abikoff1980} and where the boundaries of $\tilde{U}$ are geodesics in $\mathbb{H}^2$.

The Conway symbol for $\Gamma$ has the form $A\cdots \star abc\cdots\times\cdots\circ$, where the different symbols give rise to special elements of $\Gamma$ that correspond to \emph{features} of $\mathcal{O},$ or transformations of $\mathbb{H}^2$ that generate $\Gamma$. For example, the symmetry group $2226$ is the symmetry group generated by four rotations of order $2$, $2$, $2$, and $6$, respectively. Using the Conway symbol, one can effectively compute the Euler-Poincar\' e characteristic of $\Gamma$, to ascertain whether it represents an NEC-group.

The symmetry group $\Gamma$ of an orbifold $\mathcal{O}$ is also called the \emph{fundamental group} of $\mathcal{O}$, denoted as usual by $\pi_1(\mathcal{O})$. Additionally to acting as deck transformations on the cover $p:\mathbb{H}^2\to\mathcal{O}$, elements of $\Gamma$ also have an interpretation as homotopy classes of based closed curves in $O$~\cite{Conway2002,ratcliffe2006foundations}. In this interpretation, we can picture the elements of $\Gamma$ as homotopy equivalence classes of closed curves in the underlying topological space $O$ of $\mathcal{O}$. To make sense of this, one needs to keep track of the types of branch points of $O$ and introduce further rules for homotopies of such curves, as for example the relation $\gamma^A=1$ for some $\gamma\in \Gamma$ means that the curve representing $\gamma$ is homotopically trivial when traversed $A$ times. Moreover, curves that touch a mirror boundary (corresponding to a substring of the Conway symbol of the form $\star abc...$) in $O$ transversally must lift to a path that crosses over the mirror in the covering space. 

\begin{figure}[!htbp]
  \centering
    \includegraphics[width=\textwidth]{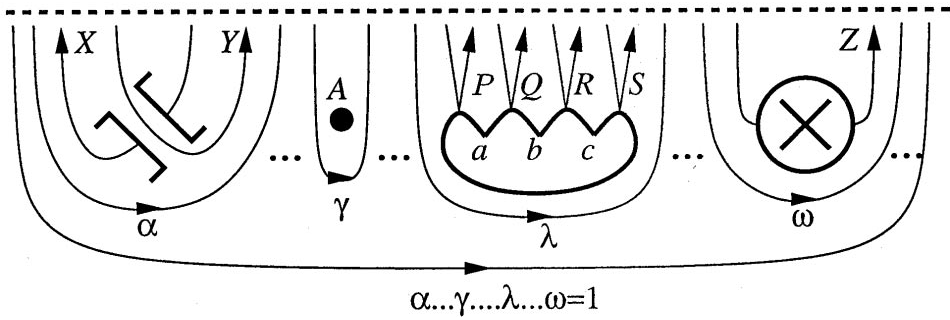}
  \caption{Paths in a compact orbifold that lift to generators of its symmetry group. Image courtesy of Daniel Huson~\cite{Conway2002}.} \label{fig:orbicurves}
\end{figure}

We briefly recall how the Conway symbol of an orbifold $\mathcal{O}$ gives rise to what we shall refer to as the \emph{standard presentation} of the orbifold symmetry group.  
There are two generators for the translations associated to each handle $\circ$, given by, say, $X$ and $Y$, and there is an oriented curve going around a handle that traces the commutator $[X,Y]=XYX^{-1}Y^{-1}=:\alpha.$ There are also generators for each gyration point of order $A$, and for a curve $\gamma$ in $\mathcal{O}$ going around the gyration point once we have $\gamma^A=1.$ Furthermore, for each mirror we have the usual Coxeter group relations, which depend on the angles of the intersecting mirrors. This means that each substring of the form $\star abc$ corresponds to a sequence of mirrors that together form a boundary component. Each of the letters $a,...$ corresponds to a mirror $\mathfrak{a}$ with relation $\mathfrak{a}^2=1,...$ and the mirror $\mathfrak{a}$ meets the mirror $\mathfrak{b}$ at an angle of $\pi/a$, which yields the relation $(\mathfrak{a}\mathfrak{b})^a=1$. As closed curves, each generator of a mirror corresponds to a closed curve from a base point that touches that mirror transversally and backtracks. In the case where the interior of the orbifold contains nontrivial features, we choose one mirror per boundary component that we give two generators, say $P$ and $S$, corresponding to the two ways of going around the mirror boundary component, and one generator $\lambda$ for the curve that goes around this boundary component once in positive orientation. We then add the relation $S=\lambda^{-1} P \lambda$. In the literature, $\lambda$ is known as a connecting generator. Figure~\ref{fig:orbicurves} illustrates the curves in $O$ representing the generators discussed here.

As a deck transformation, $\lambda$ corresponds to a hyperbolic transformation, whenever $\mathcal{O}$ has more than one nontrivial feature in its interior, or whenever there is another boundary component.
Next, going around a crosscap corresponds to a group element $\omega$ with $Z^2=\omega$, where $Z$ is a generator, corresponding to the curve passing through the crosscap once. In case $\mathcal{O}$ has (nonmirror) boundary components, every boundary curve gives rise to a generator $h_i$, a hyperbolic transformation, known as a \emph{boundary hyperbolic transformation}. For every puncture, we add a generator $p_j$, which is a parabolic deck transformation acting on $\mathbb{H}^2$~\cite{primermcgs} and corresponds to a curve around the puncture. There is one \emph{global relation} for an orbifold, namely, the product of the generating transformations discussed above has to be trivial: 
\begin{align}\label{eq:globalrel}
\gamma...h_i...p_j...\lambda...\omega...\alpha...=1.	
\end{align}
We refer to the presentation thus obtained as the \emph{standard presentation} of the fundamental group of $\mathcal{O}$. For example, the symmetry group $2226$ discussed previously has the natural relations for the orders of its elements and one extra relation, namely $r_{2_1}r_{2_2}r_{2_3}r_6=1$, where $r_{2_1}$ is a rotation about the first 2-fold rotation center and the subscripts are simply to keep track of the order. 

While a Conway symbol uniquely defines an isomorphism class of groups, this is not unambiguous as different symbols describe the same class. Most of these ambiguities are related to the order of the symbols, but there are also others~\cite{Conway2002}. In particular, for convenience, we can assume that in the presence of a crosscap, all handles are replaced by two crosscaps each~\cite{conwayzip}. From now on, when we talk about a set of generators for an orbifold, we usually mean a set of generators as provided by the Conway symbol like above. For emphasis, we will sometimes also write \emph{geometric generators}. Depending on the context, the generators may show up in a different order, which leads to a slightly different global relation \eqref{eq:globalrel}, with a permutation of the symbols. For our purposes, we need to fix an arbitrary ordering of the generators for the symmetry group under consideration. 

\subsection{Isotopic tiling theory}
The mathematical framework that serves for the enumeration of isotopy classes of tilings we present here is based on results on isotopic tiling theory~\cite{BenMyf1}. We recall the most important results for our purposes.
The elements of an orbifold fundamental group $\Gamma$ can be assigned types according to their algebraic properties and their action on the hyperbolic plane. For our purposes, we define the MCG $\Mod(\mathcal{O})$ of a hyperbolic orbifold $\mathcal{O}$ with symmetry group $\Gamma\subset \Iso(\mathbb{H}^2)$ as the isotopy classes of homeomorphisms $f$ of $\mathcal{O}$ that induce a homeomorphism of (a totally geodesic subspace of) $\mathbb{H}^2$ which induces an automorphism of $\Gamma.$ This corresponds to homeomorphisms of the underlying quotient space of $\mathcal{O}$ that preserve the types of branchings of the universal covering $U\to \mathcal{O}$. For example, for stellate orbifolds, this means that the homeomorphisms can only map one gyration point to another if they their orders of rotation agree. Recall also that the group of outer automorphisms of $\Gamma$ is defined as $\Out(\Gamma)=\Aut(\Gamma)/\Inn(\Gamma).$ Note that when the underlying topological space $O$ of $\mathcal{O}$ is orientable, we require the homeomorphisms in $\Mod(\mathcal{O})$ to be orientation-preserving. We shall refer to the following theorem as the MCG isomorphism in this paper.
\begin{theorem}{\cite{BenMyf1}}\label{thm:mcgout}
Let $\mathcal{O}$ be a nonorientable hyperbolic finite-area orbifold, with nonorientable underlying topological space. Then the MCG $\Mod(\mathcal{O})$ is isomorphic to $\Out_t(\pi_1(\mathcal{O}))$, the group of type-preserving outer automorphisms. If $\mathcal{O}$ is orientable, possibly containing mirrors, then the orientable MCG $\Mod(\mathcal{O})$ is isomorphic to $\Out^+(\pi_1(\mathcal{O}))$, the group of orientation and type-preserving automorphisms.
\end{theorem}
The enumeration of isotopy classes of tilings we present in the following sections is based primarily on this theorem and the standard presentation of the symmetry group discussed in the previous section. Given a set of geometric generators, by applying an outer automorphism to it, one obtains a conjugacy class of geometric generators. Using theorem~\ref{thm:mcgout}, one can show that different isotopy classes of tilings correspond to mutually non-conjugate sets of geometric generators for the symmetry group of the tiling that yield \emph{a priori} distinct isotopy classes of tilings~\cite{BenMyf1}. To ensure that non-conjugate sets of geometric generators yield distinct tilings, we need to assume that each edge orbit is given a different colour, for distinguishability, which gets rid of finitely many ambiguities. What is missing from this is more than a simple implementation. How an element of the MCG can be used to actually produce, effectively, the isotopy class of tiling it encodes still needs to be clarified, alongside a method to enumerate these MCG elements. Moreover, any such enumeration requires a presentation of the MCG, often for nonstandard surfaces, so we will derive presentations for some of these in section~\ref{sec:mixedbraids}. Our avenue of approach to an enumeration is to successively apply elements of the MCG to a set of geometric generators. 

\subsection{Enumerative aspects of isotopic tiling theory}\label{sec:ddtheory1}
A tiling is a partitioning of a metric space $\mathcal{X}$ into a locally finite collection of closed and bounded disks whose interiors are mutually disjoint. We call a point that is contained in at least $3$ tiles a \emph{vertex}, and the closures of connected components of the boundary of a tile with the vertices removed \emph{edges}.
Let $\mathcal{T}$ be a tiling of $\mathcal{X}$ and let $\Gamma$ be a discrete subgroup of $\Iso(\mathcal{X})$. If  $\mathcal{T}=\gamma \mathcal{T} :=\{\gamma T \;|\; T\in\mathcal{T}\}$ for all $\gamma \in \Gamma$ then we call the pair $(\mathcal{T},\Gamma)$ an \emph{equivariant tiling}.  
Two tiles $T_1,T_2\in \mathcal{T}$ are \emph{equivalent} or symmetry-related if there exists $\gamma\in\Gamma$ such that $\gamma T_1=T_2.$  
The \emph{orbit}, or equivalence class, of a tile is the subset of $\mathcal{T}$ given by images of $T$:  $\Gamma.T = \{ \gamma T$ for $ \gamma \in \Gamma \}$.
Given a particular tile $T\in\mathcal{T}$, the \emph{stabilizer subgroup} $\Gamma_T$ is the subgroup of $\Gamma$ that fixes $T$, i.e. $\Gamma_T = \{ \gamma \in \Gamma \;|\; \gamma T = T \}$.
Notice that the stabilizer subgroup of a tiling with closed disks is necessarily finite. 

A tile is called \emph{fundamental} if $\Gamma_T$ is trivial and we call the whole tiling fundamental if this is true for all tiles. 
An equivariant tiling is called \emph{tile}-$k$-\emph{transitive}, when $k$ is the number of equivalence classes of tiles under the action of $\Gamma$.  
Note that the above definitions do not require $\Gamma$ to be the maximal symmetry group for the tiling $\mathcal{T}$.
	
Two equivariant tilings $(\mathcal{T}_1, \Gamma_1)$ and $(\mathcal{T}_2, \Gamma_2)$ of $\mathcal{X}$ are \emph{equivariantly equivalent} if there is a homeomorphism, $\phi$, of $\mathcal{X}$ such that $\phi(T_1) \in \mathcal{T}_2$ for all $T_1 \in \mathcal{T}_1$ and such that $\phi$ induces a group isomorphism of $\Gamma_1$ onto $\Gamma_2$ by $\Gamma_2=\phi \Gamma_1 \phi^{-1}.$   
Dress~\cite{DRESS1987} shows that there is a complete invariant that detects when two tilings are equivariantly equivalent, for tilings of simply connected manifolds. This invariant is called the D-symbol and consists of a coloured graph that records adjacencies between tiles and their faces, augmented by weights that encode the group action of $\Gamma$ on $\mathcal{T}$. 

The enumeration of D-symbols up to a certain complexity is computationally tractable and was investigated in \cite{Huson1993,Delgado-Friedrichstilings}, allowing for efficient enumerations.

We call two equivariantly equivalent tilings on a hyperbolic surface $S$ isotopic if their edge graphs are isotopic in $S$~\cite{BenMyf1}.

Once one specifies a set of geometric generators of $\Gamma$ as transformations of $\mathbb{H}^2$, from a D-symbol one can find a (nonunique) corresponding decoration in terms of a combinatorial description of vertices and their connectivity that gives rise to a tiling of $\mathbb{H}^2$ with that D-symbol~\cite{BenMyf1}. Other isotopy classes of tilings with the same D-symbol are a result of changing the set of geometric generators according to the MCG isomorphism and using the same combinatorial description of the decoration. See figure~\ref{fig:tiling2224ht} for an illustration of different tilings resulting from such a combinatorial description. 

Every equivariant equivalence class of tilings can be built by successive applications of certain operations on tile-$1$-transitive tilings or their D-symbols, called the GLUE and SPLIT operations~\cite{Huson1993,BenMyf1}. A fundamental tile-$1$-transitive tiling can be built from any one chosen tile, \emph{the starting tile}, by applying elements of the symmetry group, resulting in a tessellation of the whole space. We therefore focus on the enumeration of tile-$1$-transitive tilings by constructing the decoration of a fundamental domain that gives rise to the tiling by enforcing the symmetries of the symmetry group. 
 
There are essentially two choices involved in the enumeration. The first is a choice of starting generators of a symmetry group that the MCG acts on to produce all other admissible sets of generators. This choice directly relates to the choice of method to produce a fundamental domain with given combinatorial structure in $\mathbb{H}^2$ from a given set of generators as a combinatorial description. Note that there are methods to produce any isotopy class of tilings from a given set of generators, just as it is possible to construct every possible tiling with given symmetry group by decorating a single fixed fundamental domain, so the choice of starting generators only makes sense after having fixed the method of construction. As a first step, we restrict ourselves to methods that yield a minimal total word length description of the Wilkie generators~\cite{Wilkie1966} of the symmetry group corresponding to edge traversals across boundaries of the fundamental starting tile in terms of the given geometric generators. Moreover, we can simplify the situation by requiring that the given generators act on a subset of the boundary of the starting tile or its neighbours sharing an edge by mapping some part of it to itself. This subset of the boundary can consist of one edge, two edges, or a point, depending on the generators. These choices leave at most a finite number of methods to produce isotopy classes of fundamental domains from given geometric generators. In practice, we will make an ad-hoc choice among these to produce a fundamental domain from a D-symbol that satisfy these conditions. We will discuss these in more detail for stellate symmetry groups below. 

 The second choice is that of an action of the MCG on sets of geometric generators on the symmetry group and, closely related, the presentation of the MCG that gives rise to the enumeration of its elements. The enumeration of MCG elements uses the fact that MCGs of orbifolds have solvable word problem~\cite{BenMyf1}. We will discuss well-known generators of MCGs in section~\ref{sec:conciso} and derive presentations of MCGs in section~\ref{sec:mixedbraids}. In section~\ref{sec:conciso}, we also explain how to enumerate the possible realizations of the geometric generators of $\Gamma=\pi_1(\mathcal{O})$ as isometries of $\mathbb{H}^2$, subsequently also referred to as the \emph{placements} of the generators, starting from a given one by deriving an action of $\Mod(\mathcal{O})$ on the sets of generators. 
	
Given a closed hyperbolic surface $S$ (such as the uniformized genus $3$ Riemannian surface that gives rise to the unit cell of the gyroid surface family), to enumerate all symmetric embeddings of graphs that lead to tilings with nontrivial symmetry group in $S$, one first has to identify all of the hyperbolic symmetry groups of $S$. Since $S$ is hyperbolic, by Hurwitz' theorem, there is a biggest discrete symmetry group $\Gamma_S$ that contains every symmetry of $S$ and whose fundamental domains in $\mathbb{H}^2$ have the smallest hyperbolic area amongst all symmetry groups of $S$. 

In practice, we find all admissible symmetry groups of $S$ by finding all of the finitely many groups $G$ with $T\subset G\subset \Gamma_S\subset \Iso(\mathbb{H}^2)$, where $T$ is the group of deck transformations of the universal Riemannian cover of $S$ by $\mathbb{H}^2$. For the purposes of EPINET, where $S$ is embedded in $\mathbb{R}^3$ into a unit cell such that it gives rise to a triply-periodic surface, $T$ is normal in $\Gamma_S$.  This can be readily checked by using the property of the unit cell that all translations under which the triply-periodic surface is invariant are a result of composing translations of the unit cell. Whence, we find all admissible subgroups by enumerating the finitely many subgroups of the finite group $\Gamma_S/T$ and adding generators for $T$.

For the purposes of enumerations of structures in $\mathbb{R}^3$ through isotopy classes of tilings it is important to check when two subgroups of $\Gamma_S$ are conjugate, since the symmetries of $\Gamma_S$ either lift to isometries of $\mathbb{R}^3$, or there is an index $2$ subgroup of $\Gamma_S$ that lifts~\cite{Perez2002,meeksphd}. When two symmetry groups are related by a symmetry that lifts to an ambient isometry in $\mathbb{R}^3$, then the resulting structures in $\mathbb{R}^3$ will also be related by a symmetry, so are either continuously deformable into one another, or are mirror images of each other. Using results on the existence of lifts of intrinsic symmetries of the diamond family of TPMS~\cite{meeksphd}, with symmetry group $\star 246$, we see that we can disregard subgroups of $\star 246$ that are conjugate by an element of $\star 246$ in the case of the primitive and diamond surfaces, as these symmetries lift to isometries of $\mathbb{R}^3$. In case of the gyroid, we need to construct all tilings that are related by some mirror symmetry of $\star 246$, as such a symmetry does not lift to an ambient isometry of $\mathbb{R}^3$. On the other hand, for the enumeration of isotopy classes of equivariant tilings on TPMS, one also needs to account for conjugate but distinct subgroups of symmetry groups $\star 246$. 

For the diamond, primitive and gyroid TPMS, the $131$ conjugacy classes of subgroups of the smallest (area-wise) symmetry group $\Gamma_S=\star 246$ have been listed~\cite{Robins2004} and the fundamental tile-$1$-transitive tilings for the Coxeter groups, generated entirely by reflections, have been enumerated and projected onto the diamond, primitive and gyroid TPMS~\cite{Ramsden2009}. Note that there is only one combinatorial class of fundamental tile-$1$-transitive tilings for Coxeter groups~\cite[Theorem $1.1$] {Lucic1990}, with a unique isotopy class~\cite{BenMyf1}. 

For a general hyperbolic symmetry group $G$, there are only a finite number of fundamental tile-$1$-transitive equivariant equivalence classes of tilings~\cite[Proposition $4$]{BenMyf1}, which can be enumerated as D-symbols for tilings~\cite{Dresshu}. Fixing the methods to produce a combinatorial class of tiling from each of the finitely many D-symbols, this translates to there only being a finite number of isotopically distinct fundamental tile-$1$-transitive tilings that arise from a given set of generators given as isometries of a hyperbolic group in $\Iso(\mathbb{H}^2)$. In the remainder of this section, we describe a method to find a nice starting set of generators for the symmetry group of a tiling that is related to minimizing the total edge length of the tiling on the surface.

Intuitively, we want to start with generators that act on subsets of tile edges and lead to a tiling on the surface $S$ with minimal total edge length, or minimal shearing. Given a fundamental domain for $\pi_1(S)$ in $\mathbb{H}^2$ with minimal total edge length, this condition translates to starting with a minimally sheared tiling in $\mathbb{H}^2$. The minimization procedures described here, including that outlined above for finding a method of construction of a combinatorial class of fundamental domain from given geometric generators as well as that for finding a starting set of generators described below cut down the possible choices of starting point for an enumeration to a finite number of equally suitable choices. The approach is based on the the fact that there are only finitely many isotopy classes of curves on a finite-area hyperbolic orbifold of length smaller than a given bound. In theory, it is therefore possible to choose a minimal representative out of the collection of curves that cut the orbifold into a disk that gives rise to a tile-$1$-transitive tiling. 

Let 
\begin{align}\label{eq:subgrouprels}
T\subset G\subset \Gamma_S\subset \Iso(\mathbb{H}^2)
\end{align}
be NEC groups, with $T$ the realization in $\Iso(\mathbb{H}^2)$ as deck transformations of the fundamental group of a closed hyperbolic surface $S$ and $\Gamma_S$ its smallest (area-wise) symmetry supergroup. Then there is a fundamental domain $F_G$ for $G$ produced by gluing copies of a fundamental domain of $\Gamma_S$, where for $F_{\Gamma_S}$, we choose a least-sheared representative. For $T$, we take a geodesical fundamental polygon $p$ in $\mathbb{H}^2$ that is least sheared. 

 We obtain a placement $Q$ of a set of geometric generators for any $G$ satisfying \eqref{eq:subgrouprels}, induced by a set of geometric generators for $\Gamma_S$ by requiring that elements of $Q$ act on the boundary of $F_G$ (or, depending on the context, a neighbouring tile) by mapping a subset of it to some subset of the boundary. One such set of generators will be our starting point and serve as a reference frame. See figures~\ref{fig:tiling3232a} and \ref{fig:tiling2224b} for an illustration of a fundamental polygon for a genus $3$ hyperbolic surface along with such a starting set of generators for the symmetry group $G=2233$ and $G=2224$, respectively. The fundamental domains are a result of gluing together copies of the uniquely determined (up to isometry) fundamental triangle for the $\star 246$ symmetry group, which plays the role of $\Gamma_S$ for the genus $3$ hyperbolic surface depicted. In both cases, the fundamental domain admits a mirror symmetry, which maps the generators to other generators of the group. For $2233$, the index $2$ supergroup $2\star 33$ is generated by the elements of $2233$ and the reflection symmetry of the fundamental domain. For $2224,$ it is the index $2$ supergroup $\star 2224$. 
 
 In more complicated cases with multiple such choices, we choose some starting set of generators minimizing the infimum of the length of the circumference of a fixed fundamental tile in $\mathbb{H}^2$ constructed in a fixed way, again with generators acting on subsets of tile edges that are incident to the fundamental tile. If we choose a minimally sheared version of $p$ in $\mathbb{H}^2$ to begin with, these two methods of finding a starting set roughly agree where they overlap. Note that an alternate natural choice for a starting set of generators can also be obtained by choosing one of the finitely many possibilities of geometric generators that act, as before, on tile edges incident to a Dirichlet fundamental domain. 
Summarizing, the idea is simply to use a least-sheared version of a tiling on $S$ resulting from a fixed construction as a starting point and choose a set of generators that acts on subsets of tile edges incident to one of the tiles. For stellate groups, using generators acting on the boundary of the index $2$ Coxeter supergroup gives a natural starting point, see Figures~\ref{fig:tiling2222276} and \ref{fig:tiling2222277}, where the indicated generators of $22222$ act on the boundary of a fundamental tile that is the result of doubling the fundamental domain for the Coxeter supergroup $\star 22222$. 

A phenomenon related to the ambiguity of the starting set of generators and illustrated by the examples above occurs when some equivariant tiling of $\mathbb{H}^2$ with symmetry group $G$ exhibits symmetries that are not elements of $G.$ For example, when a surface $S$ admits a symmetry group $R$ generated solely by rotations. A fundamental domain that admits a reflectional symmetry for $R$ is obtained by mirroring \emph{the} fundamental domain of the Coxeter supergroup of $R$ that is generated solely by reflections. The emergence of such `accidental' symmetries is a result of combinatorial tiling theory being cast in terms of equivariant tilings, where some symmetry group of a tiling must be specified, possibly without being the full group of symmetries of the tiling. These accidental symmetries may easily be broken by adding extra edges to the tiles, but because of their existence, some tilings are invariant w.r.t. a change of generators and associated decoration of the orbifold. We will see further examples of this in section~\ref{sec:steltileex} below. Notice that the emerging symmetries do not have to be part of the isometry group of the tiled surface $S$. The possible ambiguities in the description of isotopy classes of classical tilings with the MCG of their symmetry group relate to finite subgroups of the MCG~\cite[Section $7$]{BenMyf1}. It is because of these finitely many ambiguities that the enumeration that results from isotopic tiling theory as presented here is, strictly speaking, an enumeration of coloured tilings, where each edge orbit is given a different colour, for distinguishability.

The following three sections successively build towards the enumeration of isotopy classes of tilings. In section~\ref{sec:conciso}, we relate the structure of orbifolds to the sets of generators of its group. In particular, we derive an action of the mapping class group of the orbifold associated to the symmetry group of a tiling on the set of sets of geometric generators. Section~\ref{sec:mixedbraids} then uses the connection of the semi-pure braid group on the sphere to its MCG to move to an enumerative setting for tilings with stellate symmetry. We also explicitly give descriptions of certain subgroups of mapping class groups and of tilings as decorations on orbifolds. This culminates, in section~\ref{sec:steltileex}, in an array of examples of isotopically distinct tilings of the hyperbolic plane with symmetries that are commensurate with the primitive, diamond and gyroid triply-periodic minimal surfaces.

\section{Relating the structure of orbifolds to the sets of generators of its group}\label{sec:conciso}
The aim of this section is to explain how the MCG acts on the sets of geometric generators for the symmetry group of a tiling. Together with a presentation of the MCG in question, this will give rise to enumerations of isotopy classes of tilings with given symmetry group. 
\subsection{The MCG isomorphism}
The set of sets of generators of $\Gamma$ that is relevant for producing isotopically distinct tilings of a (finite-area) hyperbolic surface $M$ with symmetry group $\Gamma$ is in one-to-one correspondence to the orbifold MCG by theorem~\ref{thm:mcgout}. Just knowing that such an abstract correspondence exists is insufficient for applications though, hence we begin by making the definition of the MCG isomorphism more suitable to our purposes. 

We will describe the MCG isomorphism from theorem~\ref{thm:mcgout} by relating the action of certain generators of the MCG on curves to the resulting group elements in the orbifold fundamental group $\Gamma$, using the interpretation of group elements in $\Gamma$ as homotopy classes of curves in the associated orbifold $\mathcal{O}$. The main point is that one can interpret elements of the symmetry group $\Gamma$ of the orbifold as homotopy classes of closed curves in the labelled quotient space $O=\mathbb{H}^2/\Gamma$ that avoids the singular locus of $\mathcal{O}$. 

The following well-known result for covering spaces holds in greater generality for orbifold covering spaces.

Let $f$ be a homeomorphism of a surface $S$ with base point $p$ not in the singular locus, which lifts to a homeomorphism $\tilde{f}$ of the universal covering space and fix a base point $\tilde{p}$ in the fiber above $p$. Denote the deck transformation corresponding to the closed curve $\alpha$, with homotopy class $[\alpha]\in \pi_1(S,p)$, by $\delta_{[\alpha],p}$, and the homomorphism $f$ induced on the fundamental group by $f_*$. We then have the relation 
\begin{align}\label{eq:decklifts}
\tilde{f}\circ \delta_{[\alpha],p}\circ \tilde{f}^{-1}=\delta_{f_*([\alpha]),f(p)}.
\end{align} 
The relation \eqref{eq:decklifts} can be checked directly by applying both sides to $\tilde{f}(\tilde{p}),$ where $p$ is not in the singular locus. Then note that deck transformations are uniquely determined by where they map a single point (not in the singular locus)~\cite[Proposition $1.34$]{Hatcher}. One can use equation \eqref{eq:decklifts} to interpret the form the MCG isomorphism takes in \cite{BenMyf1}, where it matches the left hand side of the equation, to match the classical formulation of the Dehn-Nielsen-Baer theorem in \cite[Theorem $8.1$]{primermcgs}, where it corresponds to the right hand side of \eqref{eq:decklifts}.

The above suggests a method to construct the automorphisms of $\Gamma$ that are the images of homeomorphisms of $\mathcal{O}$ explicitly. First interpret a set of geometric generators of $\Gamma$ as homotopy classes of closed paths in the orbifold. One can now simply draw pictures of these generators as curves in $O$, look at how these are changed by a homeomorphism, and read off the new word representing the resulting path. Doing this for all curves representing generators of $\Gamma$ defines an induced automorphism of $\Gamma$ and, in particular, a new generating set. We treat the generating set of $\Gamma$ that corresponds to the set of curves such that cutting $O$ along these yields the standard global group relation \eqref{eq:globalrel}, which we shall write as $\Pi=1$ in the following. See figure~\ref{fig:cutsurf} for an illustration of such curves on orientable surfaces. Compare this picture with figure~\ref{fig:disksurface}(a), which shows the surface cut open along the set of curves in figure~\ref{fig:cutsurf} for a genus $3$ surface, where the $S_i$ can represent boundary components, punctures, or mirrors. Different presentations of the fundamental group correspond to a different set of generators in the orbifold and different actions of automorphisms on the sets of generators. In case one is interested in other presentations, the results can be translated to the presentation of interest, by cutting and gluing the orbifold suitably to obtain the standard presentation~\cite{zieschang1980surfaces}. Figure~\ref{fig:disksurface}(b) shows an example of another way of cutting a genus $3$ surface into a disk, corresponding to natural generators of opposite edge identifications for the fundamental dodecagon that glues up to become the uniformized version of the surface giving rise to the primitive, diamond and gyroid surface family. 
\begin{figure}[!htbp]
   \begin{subfigure}[t]{\textwidth}
\begin{center}
\begin{tikzpicture}
\tikzset{->-/.style={decoration={
  markings,
  mark=at position #1 with {\arrow{>}}},postaction={decorate}}}
	\tikzstyle{vertex}=[circle,fill,scale=0.3]
	
	\newcommand\Ra{2}
	\newcommand\An{115}
	\coordinate (M) at (-1,0);
    \draw (-3,0) arc (-180:-180+\An:2);
    \draw (-3,0) arc (-180:-180-\An:2);
    \coordinate (A) at ($(M)+({-180+\An}:\Ra)$);
    \coordinate (B) at ($(A)+(1.5,0)$);
    \draw[-] (A) to[bend left] (B);
    \path let \p1 = (A) in coordinate (A1) at (\x1,-\y1);
    \path let \p1 = (B) in coordinate (B1) at (\x1,-\y1);
    \draw[-] (A1) to[bend right] (B1);

    \draw[shorten >=0.01cm,shorten <=0.01cm,-] ($(M)+(-0.5,0)$) to[bend left] ($(M)+(0.5,0)$);
    \draw[shorten >=-0.3cm, shorten <=-0.3cm,-] ($(M)+(-0.5,0)$) to[bend right] ($(M)+(0.5,0)$);
    
    \coordinate (hole2) at  ($(M)+(3.2,0)$);
    \draw[shorten >=0.01cm,shorten <=0.01cm,-] ($(hole2)+(-0.5,0)$) to[bend left] ($(hole2)+(0.5,0)$);
    \draw[shorten >=-0.3cm, shorten <=-0.3cm,-] ($(hole2)+(-0.5,0)$) to[bend right] ($(hole2)+(0.5,0)$);
	\draw ($(hole2)+(0,2)$) arc (90:115:2);	
	\draw ($(hole2)+(0,2)$) arc (90:65:2);	
	\draw ($(hole2)-(0,2)$) arc (-90:-115:2);	
	\draw ($(hole2)-(0,2)$) arc (-90:-65:2);
	
   \coordinate (A) at ($(hole2)+({-180+\An}:\Ra)$);
    \coordinate (B) at ($(A)+(1.5,0)$);
    \draw[-] (A) to[bend left] (B);
    \path let \p1 = (A) in coordinate (A1) at (\x1,-\y1);
    \path let \p1 = (B) in coordinate (B1) at (\x1,-\y1);
    \draw[-] (A1) to[bend right] (B1);

   \coordinate (hole3) at  ($(hole2)+(3.2,0)$);
    \draw[shorten >=0.01cm,shorten <=0.01cm,-] ($(hole3)+(-0.5,0)$) to[bend left] ($(hole3)+(0.5,0)$);
    \draw[shorten >=-0.3cm, shorten <=-0.3cm,-] ($(hole3)+(-0.5,0)$) to[bend right] ($(hole3)+(0.5,0)$);
	\draw ($(hole3)+(0,2)$) arc (90:115:2);	
	\draw ($(hole3)+(0,2)$) arc (90:65:2);	
	\draw ($(hole3)-(0,2)$) arc (-90:-115:2);	
	\draw ($(hole3)-(0,2)$) arc (-90:-65:2);

    \node at ($(hole3)+(1.6,0)$) {\ldots};
    
   \coordinate (A) at ($(hole3)+({-180+\An}:\Ra)$);
    \coordinate (B) at ($(A)+(1.5,0)$);
    \draw[-] (A) to[bend left] (B);
    \path let \p1 = (A) in coordinate (A1) at (\x1,-\y1);
    \path let \p1 = (B) in coordinate (B1) at (\x1,-\y1);
    \draw[-] (A1) to[bend right] (B1);

\coordinate (hole4) at  ($(hole3)+(3.2,0)$);
    \draw[shorten >=0.01cm,shorten <=0.01cm,-] ($(hole4)+(-0.5,0)$) to[bend left] ($(hole4)+(0.5,0)$);
    \draw[shorten >=-0.3cm, shorten <=-0.3cm,-] ($(hole4)+(-0.5,0)$) to[bend right] ($(hole4)+(0.5,0)$);
	\draw ($(hole4)+(0,2)$) arc (90:115:2);	
	\draw ($(hole4)+(0,2)$) arc (90:0:2);	
	\draw ($(hole4)-(0,2)$) arc (-90:-115:2);	
	\draw ($(hole4)-(0,2)$) arc (-90:0:2);

\coordinate (Base) at ($(M)+(1.5,-0.5)$);
\draw  [fill] (Base) circle [radius=0.025];

\draw  [fill] ($(Base)-(-0.6,0.5)$) circle [radius=0.05] ;
\node at ($(Base)-(-0.1,0.5)$) {\ldots};
\draw  [fill] ($(Base)-(0.4,0.5)$) circle [radius=0.05];
\coordinate (S) at ($(Base)-(0.4,0.5)$);

\draw [green,->] (Base) to[out=90,in=0] ($(M)+(0,1)$) node[above] {$Y_1$};
\draw [green] ($(M)+(0,1)$) to[out=180,in=90] ($(M)-(1.5,0)$) to[out=-90,in=200] (Base);

\draw [->,blue] (Base) to[out=250,in=0] ($(Base)-(0.75,0.8)$) node[below] {$S_n$};
\draw [blue] ($(Base)-(0.75,0.8)$)  to[out=180,in=230] (Base);

\draw [blue,->] (Base) to[out=-22,in=-10] ($(Base)+(0.85,-0.8)$) node[below] {$S_1$} ;
\draw [blue] ($(Base)+(0.85,-0.8)$) to[out=170,in=-60] (Base);

\coordinate (1) at ($(M)+({225}:\Ra)$);

\draw [green] (Base) to[out=150,in=0] ($(M)+(0.5,0)$) ;
\draw [green,dashed] ($(M)+(0.5,0)$) to[out=-90,in=135] (1);
\draw [green,->] (Base)  to[out=-140,in=20] ($(1)+(1.5,0.3)$) node[below ] {$X_1$};
\draw [green]  ($(1)+(1.5,0.3)$)  to[out=200,in=-45] ($(1)$);

\draw [violet,->] (Base) to[out=90,in=180] ($(hole2)+(0,1)$) node[above] {$X_2$};
\draw [violet] ($(hole2)+(0,1)$) to[out=0,in=90] ($(hole2)+(1,0)$) to[out=-90,in=20] (Base);

\coordinate (2) at ($(hole2)+({60}:\Ra)$);

\draw [violet] (Base) to[out=50,in=220] ($(hole2)+(0.5,0)$)  ;
\draw [violet,dashed] ($(hole2)+(0.5,0)$) to[out=80,in=150] (2)  node[below ] {$Y_2$};
\draw [violet,->] (Base)  to[out=10,in=260] ($(hole2)+(1.5,0)$);

\draw [violet] ($(hole2)+(1.5,0)$) to[out=80,in=-30] ($(2)$);

\coordinate (3) at ($(hole3)+({60}:\Ra)$);
\coordinate (4) at ($(hole3)-(1.5,0.5)$);
\draw [cyan,->] (Base) to[out=-5,in=200] (4) to[out=20,in=180] ($(hole3)+(0,1)$)  node[above] {$X_3$};
\draw [cyan] ($(hole3)+(0,1)$) to[out=0,in=60] ($(hole3)+(1,0)$) to[out=-120,in=-8] (Base);

\draw [cyan,dashed] ($(hole3)+(0.5,0)$) to[out=70,in=150] (3)  node[below] {$Y_3$};
\draw [cyan,->] (Base)  to[out=-10,in=260] ($(hole3)+(1.3,0)$);

\draw [cyan] ($(hole3)+(1.3,0)$) to[out=80,in=-30] ($(3)$);

\draw [cyan] (Base) to[out=-6,in=-110] ($(hole3)+(0.5,0)$);

\coordinate (5) at ($(hole4)+({60}:\Ra)$);
\coordinate (6) at ($(hole4)-(1.5,0.4)$);
\draw [red,->] (Base) to[out=-15,in=240] (6) to[out=60,in=180] ($(hole4)+(0,1)$)  node[above] {$X_g$};
\draw [red] ($(hole4)+(0,1)$) to[out=0,in=60] ($(hole4)+(1,0)$) to[out=-120,in=0] ($(hole4)-(0,0.5)$) to[out=180,in=40] (6);

\draw [red,dashed] ($(hole4)+(0.5,0)$) to[out=70,in=90] ($(hole4)+(2,0)$)  node[left] {$Y_g$};
\draw [red] ($(6)+(0.5,-0.9)$)  to[out=20,in=270] ($(hole4)+(2,0)$);
\draw [red,->] (Base) to[out=-17,in=200] ($(6)+(0.5,-0.9)$); 
\draw [red] (6)  to[out=50,in=-120] ($(hole4)+(0.5,0)$);

\draw  [fill] (Base) circle [radius=0.025];


\end{tikzpicture}
\end{center}
  \end{subfigure}
  \caption{Closed curves on surfaces $S_{g,n}$ of genus $g\le 3$ with $n$ marked points that generate $\pi_1(S_{g,n})$ such that cutting along the curves produces a disk with marked points, yielding the global relation \eqref{eq:globalrel}.}\label{fig:cutsurf}
\end{figure}

\begin{figure}[!htbp]
\centering
   \begin{subfigure}[t]{0.45\textwidth}
   \vskip 0pt
\begin{center}
\begin{tikzpicture}[scale=1.5,
    tlabel/.style={pos=0.4,right=-1pt},
    baseline=(current bounding box.center)
    ]
\newdimen\R
\R=2cm
\begin{scope}[very thick,decoration={
    markings,
    mark=at position 0.4 with {\arrow{>}}}
    ] 
\draw[postaction={decorate}, green] (0:\R) --node[right]{a} (30:\R)  ;
\draw[postaction={decorate}, green] (30:\R) --node[above right]{b}  (60:\R);
\draw[postaction={decorate}, green]  (90:\R) --node[above]{a}  (60:\R)  ;
\draw[postaction={decorate}, green] (120:\R) --node[above left]{b} (90:\R);
\draw[postaction={decorate}, violet] (120:\R)  --node[above left]{c} (150:\R) ;

\draw[postaction={decorate}, violet] (150:\R) --node[left]{d} (180:\R) ;
\draw[postaction={decorate}, violet] (210:\R) --node[left]{c} (180:\R);
\draw[postaction={decorate}, violet] (240:\R) --node[below left]{d} (210:\R);
\draw[postaction={decorate}, cyan] (240:\R) --node[below left]{e} (270:\R);

\draw[postaction={decorate}, cyan] (270:\R) --node[below]{f} (300:\R);
\draw[postaction={decorate}, cyan] (330:\R) --node[below right]{e} (300:\R);

\draw[postaction={decorate}, cyan] (360:\R) --node[below right]{f} (330:\R);

\draw  [fill, blue ] (0,0) circle [radius=0.05] ;
\node[blue] at (0,0.5) {$S_2$};
\draw [->,blue] (2,0) to[out=150,in=20] (-0.1,-0.2) to[out=210,in=-90] (-0.3,0) to[out=90,in=180] (0,0.2);
\draw [blue](0,0.2) to[out=0,in=140] (2,0);

\coordinate (S2) at ($(1.3,-0.2)$); 

\draw  [fill, blue] (S2) circle [radius=0.05] node[below right ]{$S_1$};
\draw [->,blue] (2,0) to[out=170,in=20] ($(S2)+(-0.1,-0.2)$) to[out=210,in=-90] ($(S2)+(-0.3,0)$) to[out=90,in=180] ($(S2)+(0,0.2)$);
\draw [blue] ($(S2)+(0,0.2)$) to[out=0,in=140] (2,0);
\end{scope}
\end{tikzpicture}
\end{center}
\caption{}
  \end{subfigure}
 \begin{subfigure}[t]{0.45\textwidth}
 \vskip 0pt
\begin{center}
\includegraphics[width=\textwidth, height=\textwidth]{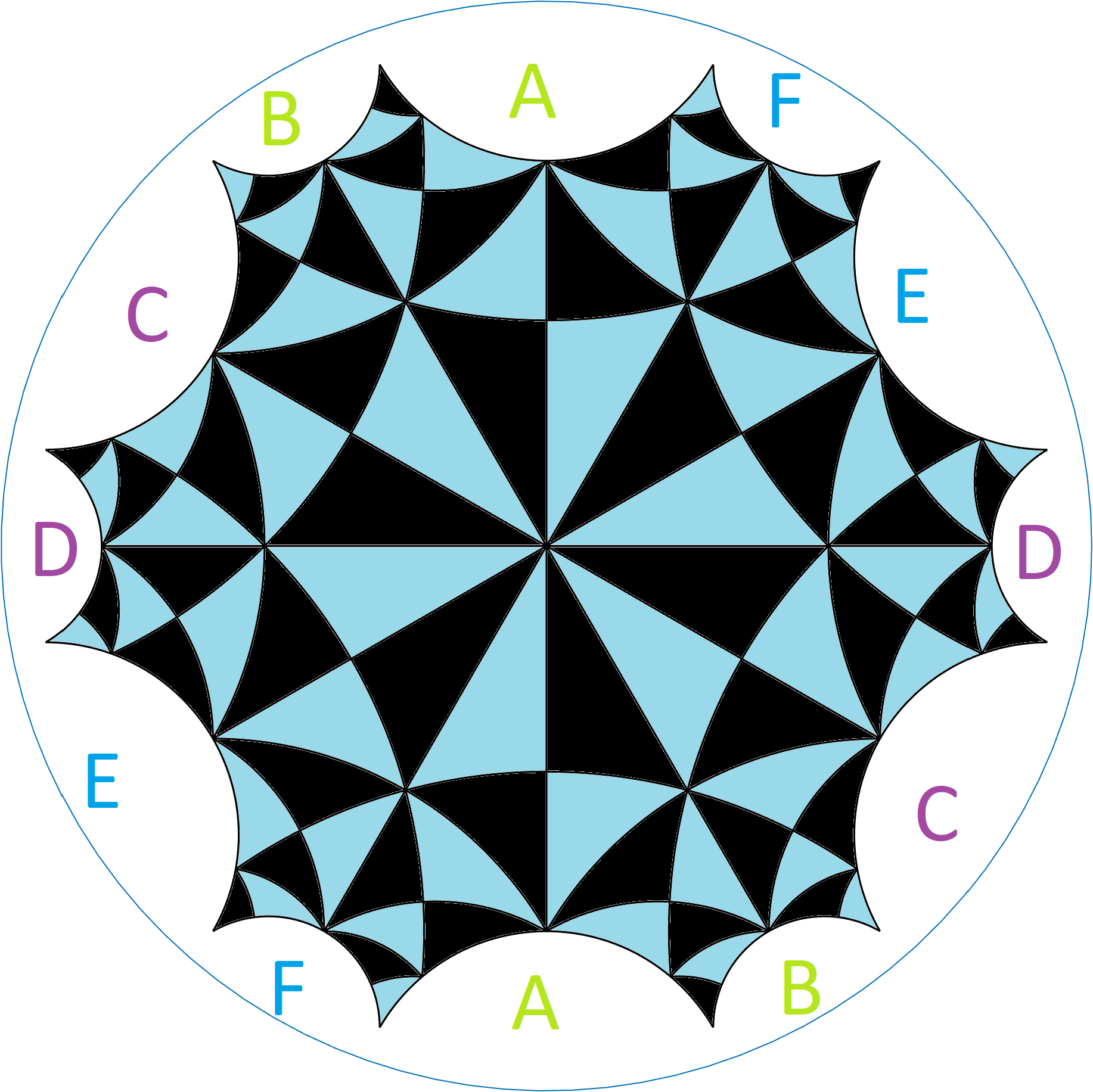}
\end{center}
\caption{}
  \end{subfigure}
  \caption{(a) A surface of genus $3$, cut open into a disk along the curves of figure~\ref{fig:cutsurf}. (b) Another way of cutting a genus $3$ surface into a disk, corresponding to the presentation of the fundamental dodecagon giving rise to the gyroid family of TPMS as the embedding of a genus $3$ surface into the three-torus. }\label{fig:disksurface}
\end{figure}

It is worth pointing out that one could, in theory, also attempt to use the method described in \cite{mccool75} to try to find an algebraic representation of some MCGs as a set of operations on the sets of generators of $\Gamma$. Moreover, the related but simpler algorithm in \cite{Dahmani2011} finds a generating set of algebraic operations, but without relations. However, we are interested in a presentation that also captures our intuition of what a complexity ordering for decorations on surfaces might look like. For this, generators for the MCG with a geometric interpretation are important. For example, the MCG of a classical surface can be generated by as few as two torsion elements~\cite{Korkmaz2004}, which is not very illuminating. Also, the method in \cite{mccool75} is very involved and most probably results in a presentation that is too complicated for most computational algorithms to handle. As far as we are aware, no explicit presentation has been derived from these methods~\cite[p. $130$]{primermcgs}. Note that computational group theory packages in computer algebra programs available today such as GAP, a programming language well-suited to studying problems in group theory~\cite{GAP4}, can only solve relatively simple problems. For example, even the word problem for the classical MCG for genus-$3$ surfaces using the well-known Gervais presentation~\cite{GERVAIS2001} is computationally too involved to comprehensively solve for the Knuth-Bendix program for GAP, despite it being well-known that the word problem is solvable~\cite[Theorem $4.2$]{primermcgs}. 

Last but not least, the geometric generators of the MCG that we use relate to the twisting of the decoration around handles of the hyperbolic surface covering, making it easier to relate them to the entanglement of the graphs produced. 

\subsection{The action of the MCG on sets of generators}
The MCGs of interest to us are extensions of classical MCGs in the general case~\cite[Section $5$]{BenMyf1}, but even when they are subgroups of classical MCGs, when the orbifold in question does not have mirror symmetries, their presentations are not treated in the literature. We present an example illustrating the situation, which will be the main focus of this paper. This is the case where the orbifold is \emph{stellate}, the only features it contains being gyration points. This restriction means that topologically, $\mathcal{O}$ is a sphere with some marked points. Note, on the other hand, that once a presentation of the symmetry group is fixed, the cyclic ordering of the gyration points is fixed and gyration points have \emph{neighbours}, corresponding to the neighbouring factors of the cyclic word $\Pi$, the left hand side of the global group relation \eqref{eq:globalrel} of $\Gamma$. 

\subsubsection{Half-twists and Dehn twists around separating curves}
First assume that $\Gamma$ is generated by rotations of the same order. For any instance of $x_ix_k$ in $\Pi$ define $\Phi(x_i):=x_ix_kx_i^{-1}$, $\Phi(x_k):=x_i$ and $\Phi(x_j):=x_j$ for $j\neq k,i$ and $\Pi_\Phi$ as the word where all $x_j$ are replaced by $\Phi(x_j)$. We see that $\Pi_\Phi=\Pi$. So, as long as the involved rotations $x_i,x_k$ are of the same order, $\Phi$ induces an automorphism of $\Gamma$. This automorphism is called a \emph{half-twist}. It is well-known that this automorphism corresponds directly to the action of the standard generators of the Braid group $B_n$ on the word $\Pi$, which also serves as a definition for $B_n$~\cite{Lando2004}. The presentation of the braid group with these generators is given by~\cite{artinbraids,Bohnenblust1947}
 \begin{align}\label{eq:braidpresentation}
\begin{split}
B_n=\langle\sigma_1,...,\sigma_{n-1}|& [\sigma_i,\sigma_j]=1  \quad |i-j|>1,\\
										& \sigma_i\sigma_{i+1}\sigma_i=\sigma_{i+1}\sigma_i\sigma_{i+1}
										\rangle.		
\end{split}
\end{align}
Note that there are three other major ways to define the braid group~\cite[Chapter $9$]{primermcgs}, and the different perspectives yield different insights into their structure.

 Stellate orbifolds are topologically spheres with marked points, leading to phenomena such as the existence of Dirac braids~\cite{Murasugi1999} for braids on the sphere. The presentation of the orientable MCG $\Mod(S_{n})$ of the sphere $S_n$ with $n$ marked points is again classical and given as~\cite[p. $128$]{primermcgs}
\begin{align}\label{eq:stelmcgpresentation}
\begin{split}
\Mod(S_{n})=\langle\sigma_1,...,\sigma_{n-1}|& [\sigma_i,\sigma_j]=1  \quad |i-j|>1,\\
										& \sigma_i\sigma_{i+1}\sigma_i=\sigma_{i+1}\sigma_i\sigma_{i+1},\\
										&(\sigma_1\cdots\sigma_{n-1})^n=1,\\
										&(\sigma_1\cdots\sigma_{n-1}\sigma_{n-1}\cdots\sigma_1)=1\rangle.		
\end{split}
\end{align} 
We will come back to discuss the relationship between the braid group and MCGs in more detail below in section~\ref{sec:mixedbraids}, where we derive group presentations of certain subgroups. The link between braid groups and MCGs on surfaces extends to more general surface braid groups and MCGs of more complicated surfaces~\cite{Birman1974}, see also the short exact sequence \eqref{eq:birmanexseq} below.

We now explain how to find automorphisms realizing general permutations using half-twists, which is useful for applying the presentation \eqref{eq:mixedbraidrep} below to stellate symmetry groups.
Assume that not all gyration points are of the same order. The above action of the generators of the braid group does not make sense for elements $x_i,x_k$ that are not of the same order. Instead, we have to consider appropriate subgroups of the braid group that only account for allowed permutations, and these correspond to subgroups of the MCG of finite index. Indeed, there is a well-known~\cite{primermcgs} short exact sequence induced by the action of $\Mod(\mathcal{O})$ on marked points of $O$:
\begin{align}\label{eq:puremcgsymm}
1\to \PMod(\mathcal{O})\to \Mod(\mathcal{O})\to \Sigma_n\to 1,
\end{align}
where $\Sigma_n$ is the symmetric group on $n$ elements and $\PMod$ is the \emph{pure mapping class group} of transformations that fix all marked points. Note that \eqref{eq:puremcgsymm} is valid for both only the orientation preserving transformations as well as the full MCG. The problem of finding appropriate subgroups realizing only allowed permutations could be solved geometrically, starting with the fact that Dehn twists around essential\footnote{A closed curve in a surface is called essential if it is not homotopic into the neighbourhood of a single, possibly marked, point.} closed loops generate the orientable pure MCG of closed orientable surfaces. Then a decomposition of a permutation into transpositions of only neighbouring gyration points yields half twists for exchanging two arbitrary marked points. In cyclic notation for permutations, after having decomposed cycles into a product of transpositions, one readily verifies that $(a,b)=(a,a+1)^{-1}(a+1,b)(a,a+1)$ gives an inductive procedure to decompose transpositions in this way. Squaring the corresponding half twist yields the Dehn twist along a simple closed loop around the two points involved. Taking those half twists that correspond to allowed permutations along with all Dehn twists generated by the remaining half twists yields a generating set for such subgroups of the MCG. With these generators, the Reidemeister-Schreier process~\cite[Section $2.4$]{Lyndon2001} can be used to find a presentation. 

The situation can be further simplified by using the known presentation of the \emph{pure braid group} $PB_n$, which we can define here as the subgroup of $B_n$ consisting of those elements that maps every generator in $\Pi$ to a conjugate of itself. This subgroup of the braid group is generated by the elements~\cite{artinbraids} 
\begin{align}\label{eq:purebraids}
a_{i,j}=(\sigma_{j-1}\sigma_{j-2}\cdots\sigma_{i+1})\sigma_i^2(\sigma_{j-1}\cdots\sigma_{i+1})^{-1}, \quad 1\le i <j\le n.
\end{align}
Regardless of the subgroup of the MCG of interest, the pure mapping class group that the pure braid group corresponds to will always be contained in it as a finite index subgroup. 
  
The action of the braid group on stellate orbifold symmetry groups discussed above can also be derived in a different way that illustrates our approach, by analysing the effect of a geometrically defined half-twist on curves in $O$ that represent generators of $\Gamma$. From a more geometric point of view, half-twists can be defined as elements of the MCG that interchange the positions of two marked points and square to what is known as a Dehn twist along a curve around the two marked points~\cite{primermcgs}. Figure~\ref{fig:halftwists} shows the derivation of the representation in $\Out(\pi_1(\mathcal{O}))$ of a half-twist around two neighbouring marked points, showing its effect on two curves corresponding to adjacent group elements of the symmetry group. It is easy to see that geometrically, two half twists around the same two gyration points yields a Dehn twist, around a curve encircling the two gyration points. The ordering of the gyration points in figure~\ref{fig:halftwistsa} is a result of requiring the standard presentation of the orbifold; see also figure~\ref{fig:cutsurf}.

  \begin{figure}[!htbp]
  \imagewidth=0.45\textwidth
\captionsetup[subfigure]{width=0.9\imagewidth,justification=raggedright}
  \begin{subfigure}[t]{0.45\textwidth}
  \centering
\begin{tikzpicture}
\coordinate (1) at (0,0);
\draw  [fill] (-1,0) circle [radius=0.05] node[above] (-1,0) {$S_2$};
\draw  [fill] (-0.3,0) circle [radius=0.05] node[above] (0,0) {$S_1$};
\node at (-2,0) {\ldots};
\draw  [fill] (-2.5,0) circle [radius=0.05] node[above] (-2.5,0) {$S_n$};
\draw [fill] (-2,1) circle [radius=0.025];
\draw [->,green] (-2,1) to[out=0,in=90] (-0.7,0) to[out=-90,in=0] (-1,-0.5);
\draw [green] (-1,-0.5) to[out=170,in=-85] (-1.5,0.8) to[out=95,in=0] (-2,1);
\draw [->,blue] (-2,1) to[out=0,in=100] (0,0) to[out=-80,in=0] (-0.3,-0.5);
\draw [blue] (-0.3,-0.5) to[out=-180,in=-85] (-0.6,0.4) to[out=95,in=0] (-2,1);
\draw [->,red] (-0.9,0.05) to[out=30,in=150] (-0.4,0.05);
\draw [<-,red] (-0.9,-0.05) to[out=-30,in=-150] (-0.4,-0.05);
\end{tikzpicture}
    \caption{}\label{fig:halftwistsa}
  \end{subfigure}
  \begin{subfigure}[t]{0.45\textwidth}
  \centering
  \begin{tikzpicture}
\draw  [fill] (-1,0) circle [radius=0.05] node[above] (-1,0) {$S_1$};
\draw  [fill] (-0.3,0) circle [radius=0.05] node[above] (0,0) {$S_2$};
\node at (-2,0) {\ldots};
\draw  [fill] (-2.5,0) circle [radius=0.05] node[above] (-2.5,0) {$S_n$};
\draw [fill] (-2,1) circle [radius=0.025];

\draw [->,blue] (-2,1) to[out=0,in=90] (0.2,0.1) to[out=-90,in=0] (-0.3,-0.7) to[out=180,in=-90] (-1.5,0);
\draw [blue] (-1.5,0) to[out=90,in=180] (-1,0.7) to[out=0,in=90] (-0.7,0)
to[out=-90,in=180]  (-0.3,-0.7);

\draw [->,green] (-2,1) to[out=0,in=100] (0,0) to[out=-80,in=0] (-0.3,-0.5);
\draw [green] (-0.3,-0.5) to[out=-180,in=-85] (-0.6,0.4) to[out=95,in=0] (-2,1);
\end{tikzpicture}
    \caption{}\label{fig:halftwistsb}
  \end{subfigure}
  \caption{(a) Curves representing group elements and the half-twist, indicated in red. (b) The resulting curves after twisting.}\label{fig:halftwists}
\end{figure}

More precisely, the half twists of figure~\ref{fig:halftwists} around $S_1$ and $S_2$ take the form 
\begin{align}\label{eq:half-twists}
	 \begin{cases} 
	S_1\mapsto S_1S_2S_1^{-1},\\
	S_2\mapsto S_1.
   \end{cases}
\end{align}

\begin{figure}[!htbp]
\imagewidth=0.31\textwidth
\captionsetup[subfigure]{width=0.9\imagewidth,justification=raggedright}
  \begin{subfigure}[t]{0.31\textwidth}
  \centering
    \includegraphics[width=\textwidth, height=\textwidth]{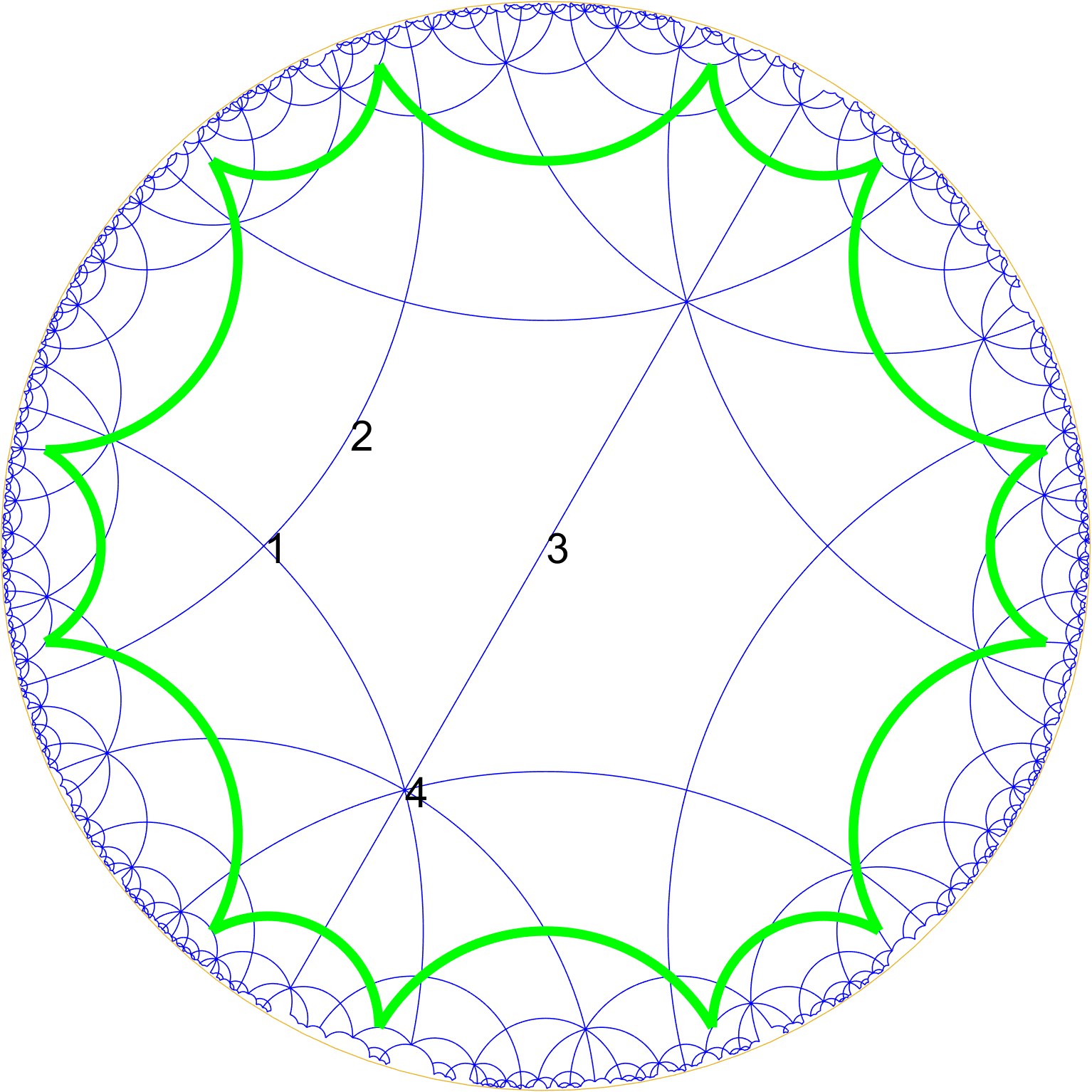}
    \caption{}\label{fig:tiling2224ht1}
  \end{subfigure}
  \begin{subfigure}[t]{0.31\textwidth}
  \centering
    \includegraphics[width=\textwidth, height=\textwidth]{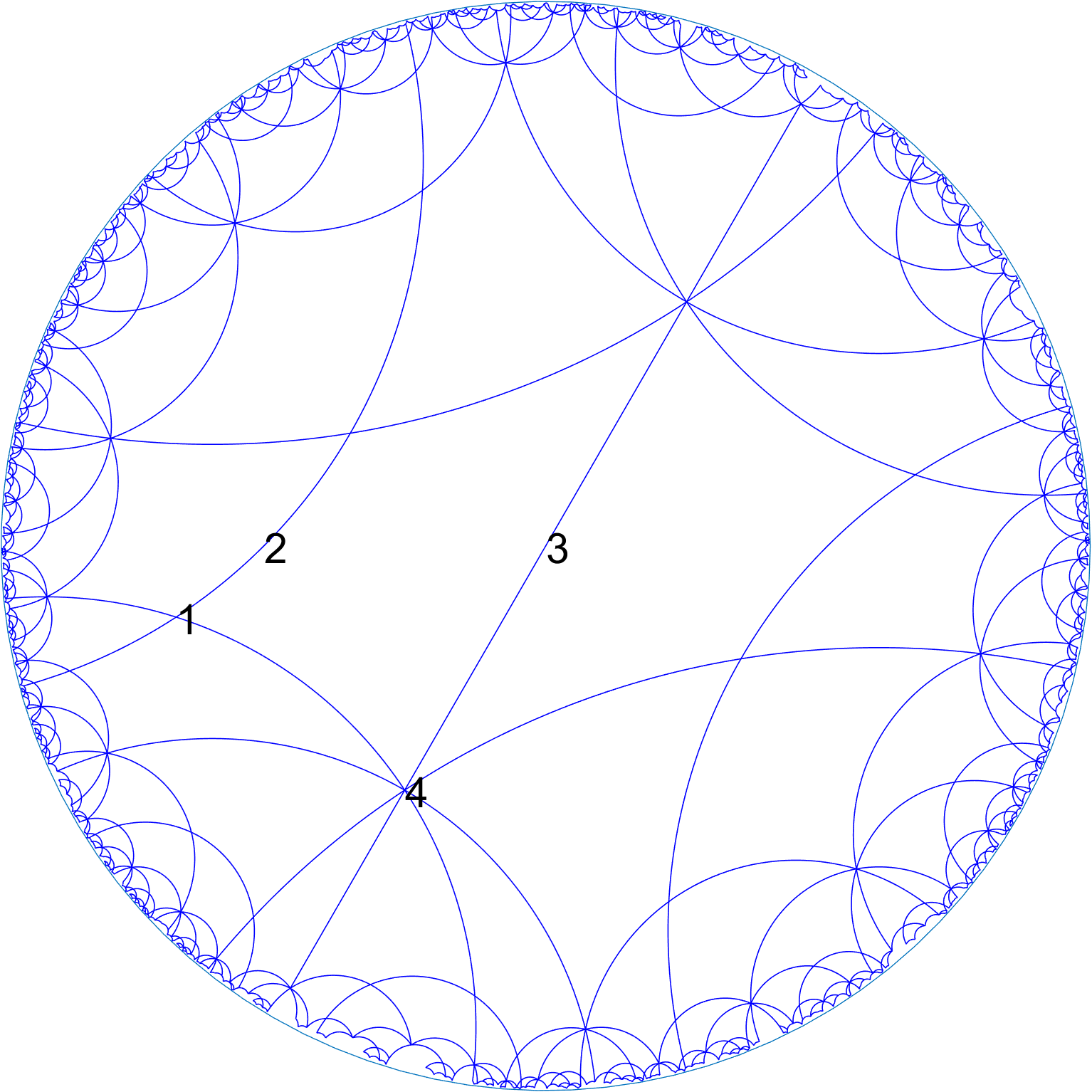}
    \caption{}\label{fig:tiling2224ht2}
  \end{subfigure}
  \begin{subfigure}[t]{0.31\textwidth}
  \centering
    \includegraphics[width=\textwidth, height=\textwidth]{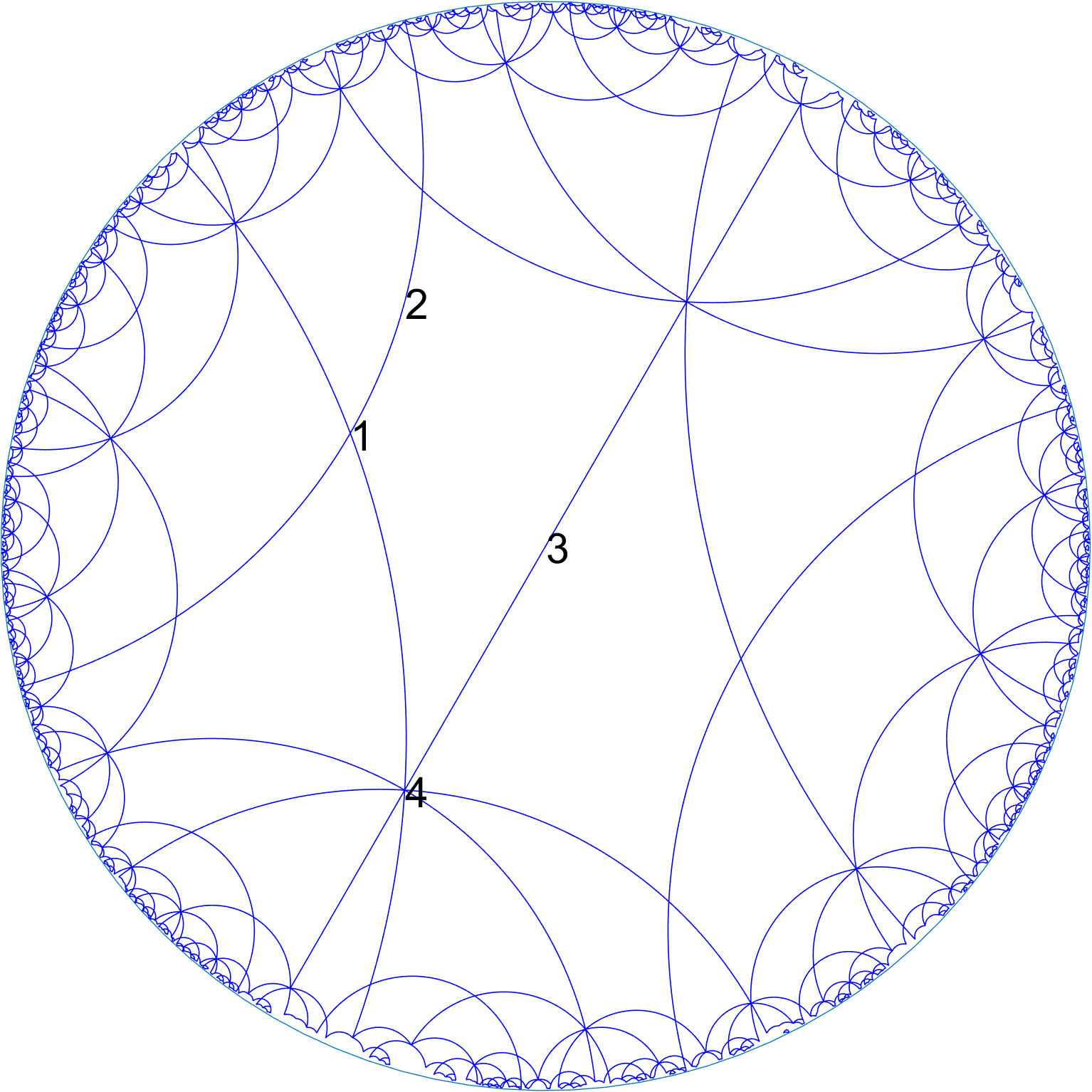}
    \caption{}\label{fig:tiling2224ht3}
  \end{subfigure}
  \caption{The effect of two inverse half twists on the isotopy class of tilings with symmetry group $2224$ (a) Starting set of generators with a tiling commensurate with the fundamental dodecagon for the P, D and gyroid TPMS shown in green. (b) A twisted tiling, commensurate with the same dodecagon. (c) The tiling with inverse twist.}\label{fig:tiling2224ht}
\end{figure}

Figure~\ref{fig:tiling2224ht} shows the effect of a half twist and its inverse on the starting set of generators with a fixed D-symbol for a tiling with symmetry group \orb{2224}.

We now analyse how Dehn twists and some other prominent elements of MCGs act on sets of generators. We will use right Dehn twists in this paper. 
\begin{example}
We first consider an example. Let $S_5$ be the five-times punctured sphere and consider the pure MCG $\PMod(S_5)\subset \Mod(S_5)$ of classes of orientable homeomorphisms that fix all punctures, with $\pi_1(S_5)$ generated by the transformations $\{e_j\}_{i=1}^5$ associated to the punctures. It is well known that $\PMod(S_5)$ is generated by Dehn twists around some essential simple closed curves~\cite{primermcgs}. After stereographic projection, we picture all punctures as lying in an ordered row and consider the set $\{t_i\}_{i=1}^s$ of Dehn twists around the set of curves $\{c_i\}_{i=1}^s$. We take the curves $c_i$ to be given by those that enclose exactly two consecutive points, and $s=4$. These generate the pure MCG~\cite{primermcgs}. The curves $c_i$ now enclose punctures $(p_k)_{k\in K_i}$ with $K_i\subset\{1,...,5\}$, which in $\pi_1(S_5)$ corresponds to $A_i:=\Pi_{j\in K_i}e_j$. Under the MCG isomorphism we obtain $\varphi:\{t_i\}_{i=1}^s \to \Out(\pi_1(S_5))$ with $\varphi(t_i)$ the automorphism that fixes each $e_j$ for $j\notin K_i$ and sends $e_j$ to $A_ie_jA_i^{-1}$ for $j\in K_i$. 
\end{example}
We are now in a position to define the action on sets of generators of some common elements and generators of MCGs. It is well-known that the pure MCG of a surface without marked points is generated by Dehn twists along essential closed curves in the orientable case~\cite{primermcgs}. In the nonorientable case, one additionally needs what are called boundary or crosscap \emph{slides}~\cite{lickorish_1965,Korkmaz2002}.

Dehn twists along essential separating simple closed curves are the most straightforward to handle and it is helpful to have an expression for such a general class of MCG elements. Similar to the above example we treat features that are adjacent in the left hand side of \eqref{eq:globalrel} differently than ones that are not. We then have a natural geometric interpretation of Dehn twists along separating simple curves that go around a chain of adjacent elements. A Dehn twist $t_c$ along a curve $c$ that encloses a chain of neighbouring features, say $\beta_1,...,\beta_k$, in that order, corresponds to the automorphism of $\Gamma$ induced by 
\begin{align}\label{eq:dehntwistauto}
\gamma\mapsto
	 \begin{cases} 
      (\beta_1\cdots\beta_k)\gamma(\beta_1\cdots\beta_k)^{-1} \quad \gamma \text{ is inside } c,\\
      \gamma \quad \text{ otherwise,}
   \end{cases}
\end{align}
where $\gamma$ is any of the generators of $\Gamma$ described in section~\ref{sec:orbs}. There is some leeway in the choice of which side is the inside of a given curve, but it is easy to see that this choice corresponds to a conjugation in $\Gamma$ of the resulting automorphism. One readily checks that the algebraic transformation in \eqref{eq:half-twists} squares to the right Dehn twist \eqref{eq:dehntwistauto} around the two involved marked points. 

\subsubsection{Orientation reversing automorphisms by reflections}
We now know how to represent simple Dehn twists and half twists around neighbouring gyration points (or punctures, or boundary components) of the same type as automorphisms of $\pi_1(\mathcal{O})$, and also how to combine these to obtain transformations  corresponding to more complicated permutations. The rest of the section deals with more complicated generators for more complicated orbifolds. 

We start with the action of a general orientation reversing homeomorphism on orbifolds whose underlying topological space is orientable. Figure~\ref{fig:reflectionaut} shows the situation for a genus $2$ surface with two marked points, or boundary components. From this, one derives the automorphism in \eqref{eq:reflectionaut}, which gives the corresponding automorphism for a general orientable orbifold with global relation of the form $S_1S_2\cdots S_n [a_1,b_1]\cdots [a_g,b_g]=1$, where both $g=0$ and $n=0$ are allowed. It is straightforward to verify that the resulting transformation $\mathcal{R}$ given by 
\begin{align}\label{eq:reflectionaut}
	 \begin{cases} 
	S_i\mapsto S_n^{-1}S_{n-1}^{-1}\cdots S_{i+1}^{-1}S_i^{-1}S_{i+1}\cdots S_n,\\
	a_j\mapsto b_{g-j+1},\\
	b_j\mapsto a_{g-j+1}.
   \end{cases}
\end{align}
defines an automorphism. Note that $\mathcal{R}$ has order two, maps each boundary component to a conjugate of its inverse, and thus defines an orientation reversing homeomorphism and automorphism for all orientable orbifolds. 

It is easy to see that the tile-$1$-transitive tiling corresponding to the set of closed curves on a classical surface that gives rise to the standard presentation of its fundamental group is invariant under this automorphism. The existence of an orientation reversing homeomorphism of order $2$ means that, for any orientable orbifold, the short exact sequence 
\begin{align}\label{eq:orientshortseq}
1\to \Mod(\mathcal{O})\to \Mod^\pm(\mathcal{O})\to \mathbb{Z}/2\mathbb{Z}\to 1
\end{align}
splits, where $\Mod^\pm(\mathcal{O})$ denotes the full MCG. This implies that the presentation of a MCG of an orientable orbifold is essentially determined by a presentation of the classes of orientation preserving transformations. To obtain a presentation of $\Mod^\pm(\mathcal{O})$ from one of $\Mod(\mathcal{O})$ first add a generator for the transformation $\mathcal{R}$ in \eqref{eq:reflectionaut} with order $2$ to the generators and relators of $\Mod(\mathcal{O})$.  Then, add relators corresponding to conjugation by $\mathcal{R}$ of each other generator, where the action of conjugation can be derived from \eqref{eq:reflectionaut}.

In case one of the $S_i$ generators in \eqref{eq:reflectionaut} represents the connecting generator of a mirror boundary, the transformation has to be adjusted slightly for each such mirror, since the curves touching the mirror boundaries are also transformed under the orientation reversing homeomorphism. We present an example illustrating the general situation. The homeomorphism depicted in figure~\ref{fig:reflectionaut} changes the direction of travel around the mirror boundary of the curves in $O$ corresponding to the mirror generators. For a boundary component $\star abc$ with connecting generator $\lambda$, we obtain, locally,   
\begin{align}\label{eq:reflectionautmirrors}
	 \begin{cases} 
	 \mathfrak{a}\mapsto \lambda^{-1}\mathfrak{a}\lambda,\\
	 \mathfrak{b}\mapsto \lambda^{-1}\mathfrak{b}\lambda,\\
	 \mathfrak{c}\mapsto \lambda^{-1}\mathfrak{c}\lambda,\\
	\tilde{\mathfrak{c}}\mapsto \mathfrak{c}=\lambda\mathfrak{\tilde{c}}\lambda^{-1},\\
		\lambda\mapsto \lambda^{-1}.
   \end{cases}
\end{align}
As in section~\ref{sec:orbs}, $\mathfrak{\tilde{c}}$ corresponds to the curve that bounces off the same mirror as $\mathfrak{c}$, but travels around the mirror boundary the other way. In case the mirror boundary component admits an extra symmetry, as for example in the case $\star 24 42$, the reflection above can appear slightly differently, since parts of a mirror might be mapped to another one, with straightforward adjustments.
Note that any homeomorphism corresponding to the automorphism defined by \eqref{eq:reflectionautmirrors} generally cannot just be supported in a small neighbourhood of the mirror boundary component, as it necessarily reverses the orientation of the boundary curve that encircles the mirror. In order for \eqref{eq:reflectionautmirrors} to yield an automorphism of the symmetry group in the general case, it needs to be combined with the transformation of the connecting generator in \eqref{eq:reflectionaut}.

\begin{figure}[!htbp]
\centering
\begin{tikzpicture}[
    tlabel/.style={pos=0.4,right=-1pt},
    baseline=(current bounding box.center)
    ]
\newdimen\R
\R=2cm
\begin{scope}[very thick,decoration={
    markings,
    mark=at position 0.4 with {\arrow{>}}}
    ] 
\draw[postaction={decorate}] (45:\R)  --node[above right]{a} (0:\R);
\draw[postaction={decorate}] (90:\R) --node[above]{b} (45:\R);
\draw[postaction={decorate}] (90:\R) --node[above]{c} (135:\R);
\draw[postaction={decorate}] (135:\R) --node[above left]{d} (180:\R);

\draw[postaction={decorate}] (225:\R) --node[below left]{c} (180:\R) ;
\draw[postaction={decorate}] (270:\R) --node[below]{d} (225:\R);
\draw[postaction={decorate}] (270:\R) --node[below]{a} (315:\R);
\draw[postaction={decorate}] (315:\R) --node[below right]{b} (360:\R);
\draw[dashed, red] (90:\R) -- (270:\R);

\draw  [fill] (0,0.5) circle [radius=0.05] ;
\node[blue] at (0.3,0.5) {$S_1$};
\draw [->,blue] (0,-2) to[out=130,in=180] (0,0.7);
\draw [blue] (0,0.7) to[out=0,in=40] (0,0.3) to[out=220,in=120] (0,-2);
\draw [->,green] (0,-2) to[out=100,in=200] (0,-0.02);
\draw [green] (0,-0.02) to[out=20,in=40] (0,-0.4) to[out=220,in=95] (0,-2);
\draw  [fill] (0,-0.2) circle [radius=0.05] ;
\node[green] at (0.3,-0.2) {$S_2$};
\end{scope}
\end{tikzpicture}
$\longrightarrow$
\begin{tikzpicture}[
    tlabel/.style={pos=0.4,right=-1pt},
    baseline=(current bounding box.center)
    ]
\newdimen\R
\R=2cm
\begin{scope}[xshift=6cm,very thick,decoration={
    markings,
    mark=at position 0.4 with {\arrow{>}}}
    ] 
\draw[postaction={decorate}] (45:\R)  --node[above right]{d} (0:\R);
\draw[postaction={decorate}] (90:\R) --node[above]{c} (45:\R);
\draw[postaction={decorate}] (90:\R) --node[above]{b} (135:\R);
\draw[postaction={decorate}] (135:\R) --node[above left]{a} (180:\R);

\draw[postaction={decorate}] (225:\R) --node[below left]{b} (180:\R) ;
\draw[postaction={decorate}] (270:\R) --node[below]{a} (225:\R);
\draw[postaction={decorate}] (270:\R) --node[below]{d} (315:\R);
\draw[postaction={decorate}] (315:\R) --node[below right]{c} (360:\R);
\draw[dashed, red] (90:\R) -- (270:\R);

\draw  [fill] (0,0.5) circle [radius=0.05] ;
\node[blue] at (-0.3,0.5) {$S_1$};
\draw [->,blue] (0,-2) to[out=50,in=0] (0,0.7);
\draw [blue] (0,0.7) to[out=180,in=130] (0,0.3) to[out=-40,in=60] (0,-2);
\draw [->,green] (0,-2) to[out=80,in=340] (0,-0.02);
\draw [green] (0,-0.02) to[out=160,in=140] (0,-0.4) to[out=320,in=85] (0,-2);
\draw  [fill] (0,-0.2) circle [radius=0.05] ;
\node[green] at (-0.3,-0.2) {$S_2$};
\end{scope}
\end{tikzpicture}
\caption{An orientation reversing homeomorphism on an orientable orbifold of genus $2$, with two marked points, corresponding to a reflection across the red line, illustrating the general situation.}\label{fig:reflectionaut}
\end{figure}
\subsubsection{Crosscap transpositions and boundary slides}
In the case of an orbifold with a nonorientable underlying quotient space, we do not need an orientation reversing transformation like above as a generator, by theorem~\ref{thm:mcgout}. We now turn to finding the action of the MCG on the geometric generators of $\pi_1(\mathcal{O})$ for a nonorientable orbifold $\mathcal{O}$ with crosscaps. In addition to Dehn twists along two-sided curves (for which Dehn-twists are well-defined), \emph{crosscap slides} and \emph{boundary slides} are needed to generate the MCG of a nonorientable surface~\cite{Korkmaz2002}. MCGs of nonorientable surfaces have been, like their orientable counterparts, extensively studied in the literature, but are somewhat less well understood~\cite{Lickorish1964,nonorientmcgbounds,Parlak2017,Lesniak2017}. 

\begin{figure}[!htbp]
\captionsetup[subfigure]{width=0.9\imagewidth,justification=raggedright}
  \begin{subfigure}[t]{\textwidth}
  \centering
  \includegraphics[width=\textwidth]{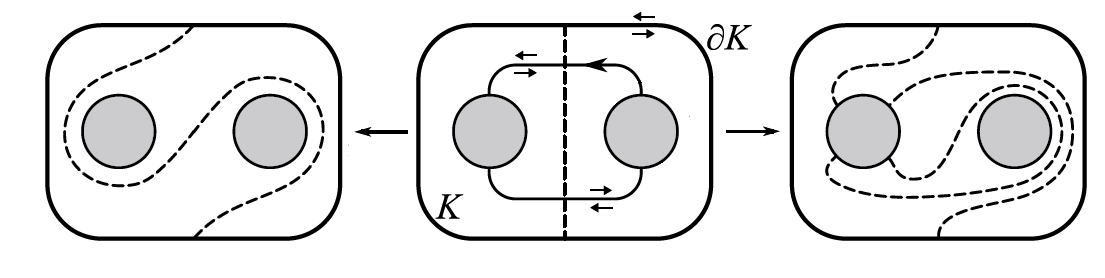}
  \end{subfigure}
  \caption{The effect on the dashed curve of a crosscap transposition on the left and a crosscap slide on the right. The homeomorphism is supported in the Kleinbottle $K$ with one boundary component $\partial K$.}\label{fig:crosscapslide}
\end{figure}

Despite there being to date no known presentation for the MCG of a general nonorientable surface in terms of these generators, they are very natural, cannot be easily decomposd any further into simpler constituent geometric transformations of other types and are supported in very small subsurfaces. We also look at \emph{crosscap transpositions}~\cite{Parlak2017}. See figure~\ref{fig:crosscapslide} for an illustration of both a crosscap transposition and crosscap slide. Roughly speaking, crosscap and boundary slides correspond to sliding the crosscap, resp. boundary along a one-sided loop~\cite{Korkmaz2002}. 
 The form all of these transformations take as automorphisms are derived, like before, from the geometric picture of the surface $O$ with a set of curves such that cutting the surface along these curves yields a planar polygon with identifications of edges to produce the global relation of the standard presentation of $\pi_1(\mathcal{O})$ in section~\ref{sec:orbs}. 

Let $\zeta$ be a two sided loop that passes through two neighbouring crosscaps represented by $A$ and $B$, and encircles no other features. Let $\alpha$ be the curve associated with $A$. A crosscap transposition $U_{\alpha,\zeta}$ is supported in a neighbourhood of $\zeta\cup\alpha$, which is a Klein bottle $K$ with one boundary component, see figure~\ref{fig:crosscapslidecurves}, and has the following representation as an automorphism, where all other generators are constant. 
\begin{align}\label{eq:crosscaptransrep}
U_{\alpha,\zeta}:(A,B)\mapsto (B,(B^2)^{-1}AB^2)
\end{align}
The Dehn twist $T_\zeta$ has the representation
\begin{align}\label{eq:twosideddehn}
T_\zeta : (A,B)\mapsto (AB^{-1}A^{-1},AB^2).
\end{align}
The crosscap slide $Y_{\alpha,\zeta}=U_{\alpha,\zeta}T_\zeta$ is represented by 
\begin{align}\label{eq:crossslide}
Y_{\alpha,\zeta} : (A,B)\mapsto (AB^2,(AB^2)^{-1}B^{-1}(AB^2)).
\end{align}
We exemplify the derivations of the above equations by drawing the picture for \eqref{eq:twosideddehn} in figure~\ref{fig:twistsnonorient}, with figure~\ref{fig:crosscapslidecurves1} showing the result of the twisting, from which \eqref{eq:twosideddehn} can be read off. Outside the shown neighbourhood, the transformation is the identity.
It is straightforward to check that these transformation yield automorphisms of $\Gamma$ and that they satisfy 
\begin{align*}
Y_{\alpha,\zeta}^2=U_{\alpha,\zeta}^2=T_{\partial K},
\end{align*}
which constitutes an algebraic proof of this well-known geometric relation~\cite{lickorish_1965}. Note that \eqref{eq:dehntwistauto} can be used to compute $T_{\partial K}$.
A boundary slide $S_A$ in a nonorientable surface that has a boundary or puncture $P$ neighbouring the cross cap $A$ is represented by 
\begin{align}\label{eq:boundslide}
(P,A)\mapsto((PA)P^{-1}(PA)^{-1},PA),
\end{align}
which again is readily seen to induce an automorphism and again comes from drawing a picture with curves. Topologically, boundaries behave exactly like punctures and like marked (gyration) points. In particular, the half-twists discussed in the context of gyration points and punctures give rise to similar half-twists for interchanging (mirror) boundaries. Note, though, that individual mirror components are also transformed, in the same way as the corresponding connecting generators given by curves encircling the mirror boundary components, as \eqref{eq:reflectionautmirrorsslide} below illustrates. Also, the boundary slide of a nonorientable orbifold for the case of a mirror boundary component also reverses the orientation of the curve around the boundary and therefore uses a version of the transformation \eqref{eq:reflectionautmirrors} for the transformation of the individual mirrors of the boundary component, as we now explain. We again use the example that led to \eqref{eq:reflectionautmirrors} to illustrate the general situation. If $P$ in \eqref{eq:boundslide} represents a curve around a mirror boundary component $\star abc$, then $S_A$ acts on generators as 
\begin{align}\label{eq:reflectionautmirrorsslide}
	 \begin{cases} 
	 \mathfrak{a}\mapsto (PA)P^{-1}\mathfrak{a}P(PA)^{-1}\lambda,\\
	 \mathfrak{b}\mapsto (PA)P^{-1}\mathfrak{b}P(PA)^{-1},\\
	 \mathfrak{c}\mapsto (PA)P^{-1}\mathfrak{c}P(PA)^{-1},\\
	\tilde{\mathfrak{c}}\mapsto (PA)\mathfrak{c}(PA)^{-1}=S_A(P)^{-1}S_A(\mathfrak{c})S_A(P),\\
	(P,A)\mapsto((PA)P^{-1}(PA)^{-1},PA).
   \end{cases}
\end{align}

 \begin{figure}[!htbp]
 \imagewidth=0.45\textwidth
\captionsetup[subfigure]{width=0.9\imagewidth,justification=raggedright}
  \begin{subfigure}[t]{0.45\textwidth}
  \centering
\begin{tikzpicture}[scale=2]
\tikzset{->-/.style={decoration={
  markings,
  mark=at position #1 with {\arrow{>}}},postaction={decorate}}}
\begin{scope}
\draw[clip] [rounded corners] (0,0) rectangle (2,2);
\draw [rounded corners, thick] (0,0) rectangle (2,2);

\draw[purple]  (0.5,1) to[out=0,in=0]  (0.5,3);
\draw[purple,->-=0.8,>=stealth] (1,3) to[out=180,in=180]  (0.5,1);

\draw[blue] (1,3) to[out=0,in=0]  (1.5,1);
\draw[blue,>=stealth,->-=.8] (2,3) to[out=180,in=180]  (1.5,1);
\draw [rounded corners] (0.5,0.5) rectangle (1.5,1.5);
\draw [->,>=stealth] (0.9,1.5) -- (0.7,1.5);
\draw [fill,lightgray] (0.5,1) circle [radius=0.2];
\draw (0.5,1) circle [radius=0.2];
\node[purple] at (0.5,1) {$A$};
\draw [fill,lightgray] (1.5,1) circle [radius=0.2];
\draw (1.5,1) circle [radius=0.2];
\node[blue] at (1.5,1) {$B$};
\draw[-{Latex[length=1mm, width=0.5mm]}] (0.9,0.45) -- (0.7,0.45);
\draw[-{Latex[length=1mm, width=0.5mm]}] (0.7,0.55) -- (0.9,0.55);

\draw[-{Latex[length=1mm, width=0.5mm]}] (1.2,1.45) -- (1.4,1.45);
\draw[-{Latex[length=1mm, width=0.5mm]}] (1.4,1.55) -- (1.2,1.55);

\node[below] at (1.5,0.5) {\tiny $\zeta$};
\end{scope}
\node[right] at (2,2) {$\partial K$};

\end{tikzpicture}
    \caption{}\label{fig:crosscapslidecurves}
  \end{subfigure}
  \begin{subfigure}[t]{0.45\textwidth}
  \centering
\begin{tikzpicture}[scale=2]
\tikzset{->-/.style={decoration={
  markings,
  mark=at position #1 with {\arrow{>}}},postaction={decorate}}}
\begin{scope}
\draw[clip] [rounded corners] (0,0) rectangle (2,2);
\draw [rounded corners, thick] (0,0) rectangle (2,2);

\draw[blue,>=stealth,->-=.7] (2.2,3) to[out=180,in=30]  (1,1.7) to[out=-150,in=30] (0.4,1.6)
 to[out=-150,in=170] (0.5,1) to[out=-30,in=180] (0.9,0.6); 
 \draw[blue,>=stealth,->-=.5] (0.9,0.6) to[out=0,in=-135] (1.5,1);
 \draw[blue,>=stealth,->-=.8]   (1.5,1) to[out=45,in=-80] (1.6,1.4)  to[out=100,in=50] (1.1,1.6) 
 to[out=-130,in=180] (1.5,1);

\draw[blue,>=stealth,->-=.3] (1.5,1) to[out=0,in=0]  (1,3);

\draw[purple,->-=0.8,>=stealth] (0.6,3) to[out=180,in=220]  (0.5,1);
\draw[purple,->-=0.8,>=stealth]  (0.5,1) to[out=-135,in=180]  (0.7,0.3) to[out=0,in=-125]  (1.6,0.6)
to[out=65,in=-45]  (1.5,1);
\draw[purple,->-=0.5,>=stealth]  (1.5,1) to[out=-180,in=0] (0.5,1);
\draw[purple,->-=0.5,>=stealth] (0.5,1) to[out=180,in=-90] (0.2,1.3) to[out=90,in=-100] (0.3,1.8)
to[out=80,in=-90] (0.6,3) ;

\draw [rounded corners] (0.5,0.5) rectangle (1.5,1.5);
\draw [->,>=stealth] (0.9,1.5) -- (0.7,1.5);
\draw [fill,lightgray] (0.5,1) circle [radius=0.2];
\draw (0.5,1) circle [radius=0.2];
\node[purple] at (0.5,1) {$A$};
\draw [fill,lightgray] (1.5,1) circle [radius=0.2];
\draw (1.5,1) circle [radius=0.2];
\node[blue] at (1.5,1) {$B$};
\draw[-{Latex[length=1mm, width=0.5mm]}] (0.9,0.45) -- (0.7,0.45);
\draw[-{Latex[length=1mm, width=0.5mm]}] (0.7,0.55) -- (0.9,0.55);

\draw[-{Latex[length=1mm, width=0.5mm]}] (1.2,1.45) -- (1.4,1.45);
\draw[-{Latex[length=1mm, width=0.5mm]}] (1.4,1.55) -- (1.2,1.55);

\node[below] at (1.5,0.5) {\tiny $\zeta$};
\end{scope}
\node[right] at (2,2) {$\partial K$};

\end{tikzpicture}
    \caption{}\label{fig:crosscapslidecurves1}
  \end{subfigure}
  \caption{Dehn twists around two sided curves in nonorientable surfaces. (a) The Kleinbottle $K$ with boundary and two-sided curve $\zeta$. (b) The resulting curves, in red and blue, after twisting right around $\zeta$.}\label{fig:twistsnonorient}
\end{figure}
Boundary and crosscap slides for non-neighbouring features are a result of combining transpositions of boundaries and crosscaps appropriately, similar to the situation for gyration points in the discussion above following \eqref{eq:puremcgsymm}. As far as we are aware, the only known presentations of MCGs of some classes of nonorientable surfaces with punctures make use of the generators discussed here~\cite{Paris2013,Szepietowski2008}. 

\subsubsection{Point-pushing, Dehn twists, Lickorish's generators}
In the orientable case, following \cite{birmanbraidsmcgs} and \cite{birmanbraids}, we use a similar approach of finding representations of generators corresponding to sliding singular points or boundaries around handles. The pictures one has to draw quickly become rather involved, and we again focus on geometric generators here, for which the standard relation \eqref{eq:globalrel} for the orbifold groups limits the ways of choosing the curves that represent the orbifold group elements. We work with the curves depicted in figure~\ref{fig:cutsurf}.
 
\begin{figure}[!htbp]
\captionsetup[subfigure]{width=0.9\imagewidth,justification=raggedright}
  \begin{subfigure}[t]{\textwidth}
  \centering
  \includegraphics[width=\textwidth]{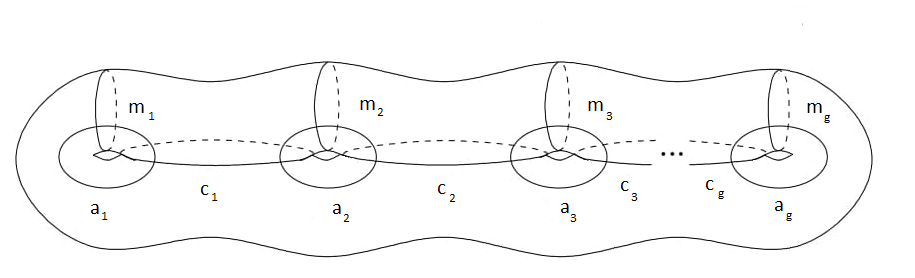}
  \end{subfigure}
  \caption{Curves corresponding to Lickorish's generators for the MCG of a closed surface of genus $g$.}\label{fig:licksgens}
\end{figure}

There are different prominent sets of generators for the MCG of a classical surface, with different properties. We will focus on the Dehn twists around the closed curves depicted in figure~\ref{fig:licksgens}, known as \emph{Lickorish generators}, along with point-pushes as defined in \cite{birmanbraids}. Gervais' presentation is another prominent presentation~\cite{GERVAIS2001}. After some experimentation with the two presentations with GAP and KBMAG~\cite{kbmag}, we found that the Lickorish set of generators seems to be easier to handle in practice. Figure~\ref{fig:licksgens} shows Lickorish's generators for the MCG of a classical surface~\cite{lickorish_1965}.

Assume that $\Pi$ has a subword of the form $S[X,Y]$. Then, moving $S$ around $Y$ corresponds to the situation of figure~\ref{fig:pointpush}. Figure~\ref{fig:pointpushb} shows the effect on the blue curve of pushing $S$ around the red curve. 

  \begin{figure}[!htbp]
\imagewidth=0.45\textwidth
\captionsetup[subfigure]{width=0.9\imagewidth,justification=raggedright}
  \begin{subfigure}[t]{0.45\textwidth}
  \centering
\begin{tikzpicture}
\tikzset{->-/.style={decoration={
  markings,
  mark=at position #1 with {\arrow{>}}},postaction={decorate}}}
	\tikzstyle{vertex}=[circle,fill,scale=0.3]
	
	\newcommand\Ra{3}
	\newcommand\An{115}
	\coordinate (M) at (-1,0);
    \draw (-4,0) arc (-180:-180+\An:3);
    \draw (-4,0) arc (-180:-180-\An:3);
    \coordinate (A) at ($(M)+({-180+\An}:\Ra)$);
    \coordinate (B) at ($(A)+(2.5,0)$);
    \draw[-] (A) to[bend left] (B);
    \path let \p1 = (A) in coordinate (A1) at (\x1,-\y1);
    \path let \p1 = (B) in coordinate (B1) at (\x1,-\y1);
    \draw[-] (A1) to[bend right] (B1);

    \draw[shorten >=0.01cm,shorten <=0.01cm,-] ($(M)+(-0.5,0)$) to[bend left] ($(M)+(0.5,0)$);
    \draw[shorten >=-0.3cm, shorten <=-0.3cm,-] ($(M)+(-0.5,0)$) to[bend right] ($(M)+(0.5,0)$);
    \node at ($(M)+(3,0)$) {\ldots};

\coordinate (Base) at ($(M)+(1.5,-0.5)$);

\coordinate (S) at ($(Base)-(0.5,1)$);

\draw [green,->] (Base) to[out=90,in=0] ($(M)+(0,1)$) node[above] {Y};
\draw [green] ($(M)+(0,1)$) to[out=180,in=90] ($(M)-(1.5,0)$) to[out=-90,in=200] (Base);

\draw [->,blue] (Base) to[out=270,in=0] ($(Base)-(0.85,1.3)$);
\draw [blue] ($(Base)-(0.85,1.3)$)  to[out=180,in=230] (Base);

\coordinate (1) at ($(M)+({225}:\Ra)$);

\draw [cyan] (Base) to[out=150,in=0] ($(M)+(0.5,0)$) ;
\draw [cyan,dashed] ($(M)+(0.5,0)$) to[out=-90,in=135] (1);
\draw [cyan,->] (Base)  to[out=-140,in=20] ($(1)+(1.5,0.3)$) node[above ] {X};
\draw [cyan]  ($(1)+(1.5,0.3)$)  to[out=200,in=-45] ($(1)$);

\draw [red] (S) to[out=245,in=-100] ($(-3.7,0)$)  to[out=90,in=180] ($(M)+(-0.2,0.9)$);
\draw [red] ($(M)+(-0.2,0.9)$)  to[out=0,in=100] ($(M)+(1,0)$) to[out=-80,in=70] (S);

\coordinate (1) at ($(M)+({260}:\Ra)$);

\draw [violet] (S) to[out=110,in=-60] ($(M)+(0.5,0)$) ;
\draw [violet,dashed] ($(M)+(0.5,0)$) to[out=-90,in=170] (1);
\draw [violet] (S)  to[out=-100,in=-10] ($(1)$);

\draw  [fill] (Base) circle circle [radius=0.025];

\draw  [fill] ($(Base)-(-0.7,1)$) circle [radius=0.05] ;
\node at ($(Base)+(0.2,-1)$) {\ldots};
\draw  [fill] ($(Base)-(0.5,1)$) circle [radius=0.05] node[left,blue] {$S$};


\end{tikzpicture}
    \caption{}\label{fig:pointpusha}
  \end{subfigure}
  \begin{subfigure}[t]{0.45\textwidth}
  \centering

\begin{tikzpicture}
\tikzset{->-/.style={decoration={
  markings,
  mark=at position #1 with {\arrow{>}}},postaction={decorate}}}
	\tikzstyle{vertex}=[circle,fill,scale=0.3]
	
	\newcommand\Ra{3}
	\newcommand\An{115}
	\coordinate (M) at (-1,0);
    \draw (-4,0) arc (-180:-180+\An:3);
    \draw (-4,0) arc (-180:-180-\An:3);
    \coordinate (A) at ($(M)+({-180+\An}:\Ra)$);
    \coordinate (B) at ($(A)+(2.5,0)$);
    \draw[-] (A) to[bend left] (B);
    \path let \p1 = (A) in coordinate (A1) at (\x1,-\y1);
    \path let \p1 = (B) in coordinate (B1) at (\x1,-\y1);
    \draw[-] (A1) to[bend right] (B1);

    \draw[shorten >=0.01cm,shorten <=0.01cm,-] ($(M)+(-0.5,0)$) to[bend left] ($(M)+(0.5,0)$);
    \draw[shorten >=-0.3cm, shorten <=-0.3cm,-] ($(M)+(-0.5,0)$) to[bend right] ($(M)+(0.5,0)$);
    \node at ($(M)+(3.5,0)$) {\ldots};

\coordinate (Base) at ($(M)+(1.5,-0.5)$);

\coordinate (S) at ($(Base)-(0.5,1)$);
\coordinate (1) at (-2.5,-1.5);

\draw [blue,->] (Base)   to[out=-110,in=70] ($(S)+(0.2,0)$) to[out=-110,in=0] ($(M)-(0,2.5)$) to[out=180,in=-90]
($(M)-(2,0)$) to[out=90,in=180] ($(M)+(0,1.5)$); 
\draw [blue] ($(M)+(0,1.5)$)  to[out=0,in=90] ($(S)+(0.1,0)$) to[out=-90,in=0] ($(S)-(0,0.1)$) to[out=180,in=-90] ($(S)-(0.1,0)$) to[out=90,in=0] ($(M)+(0,0.9)$); 
\draw [blue] ($(M)+(0,0.9)$) to[out=180,in=90] ($(M)+(-1,0)$) to[out=-90,in=225] ($(S)+(0,-0.2)$) to[out=45,in=-120] (Base);

\draw  [fill] (Base) circle circle [radius=0.025];

\draw  [fill] ($(Base)-(-0.7,1)$) circle [radius=0.05] ;
\node at ($(Base)+(0.2,-1)$) {\ldots};
\draw  [fill] ($(Base)-(0.5,1)$) circle [radius=0.05] node[left,blue] {$S$};


\end{tikzpicture}
    \caption{}\label{fig:pointpushb}
  \end{subfigure}
  \caption{Point pushes in orientable surfaces. (a) Curves representing group elements (blue, green, cyan), and both available pushes, indicated in red and violet. (b) The effect on $S$ of pushing around the red curve.}\label{fig:pointpush}
\end{figure}
Ultimately, we find that we can write the point push around of $S$ around $Y$ as
{\small
\begin{align}\label{eq:pointtwistY}
\begin{cases} 
S\mapsto SY^{-1}SYS^{-1}\\
X\mapsto SY^{-1}S^{-1}YXS^{-1}\\
Y\mapsto SYS^{-1}.
 \end{cases}
\end{align}}

A similar picture for a push around $X$ yields
{\small
\begin{align}\label{eq:pointtwistX}
\begin{cases} 
S\mapsto SXSX^{-1}S^{-1}\\
X\mapsto SXS^{-1}\\
Y\mapsto YS^{-1}.
 \end{cases}
\end{align}}
We see that $Y$ and $X$ are mapped to hyperbolic transformations in both cases. Notice that in both \eqref{eq:pointtwistY} and \eqref{eq:pointtwistX} `forgetting' that $S$ is a feature yields the trivial automorphism. 

Assume now that $\Pi$ has a subword of the form $[a_1,b_1][a_2,b_2]$. Note that the free homotopy classes of $X_i$ in figure~\ref{fig:cutsurf} and those of $a_i$ in figure~\ref{fig:licksgens} coincide.
 In Lickorish's generators for the MCG of a classical surface, there are two twists per handle that leave each commutator relation invariant on its own, corresponding to the twists around $a_i,m_i$ in figure~\ref{fig:licksgens}. In figure~\ref{fig:pointpusha}, these are the twist around the curves around $X$ and $Y$, respectively. As automorphisms, these are respectively given by
\begin{align}\label{eq:singhandledehns1}
T_{a_i}:\begin{cases} 
a_i\mapsto a_i\\
b_i\mapsto b_ia_i,
 \end{cases}
\end{align}
and
\begin{align}\label{eq:singhandledehns2}
T_{m_i}:\begin{cases} 
a_i\mapsto a_ib_i^{-1}\\
b_i\mapsto b_i,
 \end{cases}
\end{align}
for $i=1,2$. The more complicated twist among the Lickorish generators, namely the $c_k$ in figure~\ref{fig:licksgens}, mingles two adjacent handles. We have to treat the case $k=1$ separately from the others, as illustrated by figure~\ref{fig:cutsurf}. Figure~\ref{fig:lickc12} shows the curves $c_1$ and $c_2$ on a genus $3$ surface, cut open into a disk along the curves of the standard presentation, shown in figure~\ref{fig:cutsurf}. We obtain
\begin{align}\label{eq:twohandledehnsc1}
T_{c_1}:\begin{cases} 
a_1\mapsto a_1,\\
b_1\mapsto a_2b_2^{-1}a_2^{-1}b_1a_1,\\
a_2\mapsto a_2b_2^{-1}a_2^{-1}b_1a_1b_1^{-1}a_2,\\
b_2\mapsto b_2,
 \end{cases}
\end{align}
for the twist along $c_1.$
\begin{figure}
\centering
  \begin{subfigure}[t]{0.45\textwidth}
  \centering
\begin{tikzpicture}[
    tlabel/.style={pos=0.4,right=-1pt},
    baseline=(current bounding box.center)
    ]
\newdimen\R
\R=2cm
\begin{scope}[very thick,decoration={
    markings,
    mark=at position 0.4 with {\arrow{>}}}
    ] 
\draw[postaction={decorate}, green] (0:\R) --node[right]{a} (30:\R)  ;
\draw[postaction={decorate}, green] (30:\R) --node[above right]{b}  (60:\R);
\draw[postaction={decorate}, green]  (90:\R) --node[above]{a}  (60:\R)  ;
\draw[postaction={decorate}, green] (120:\R) --node[above left]{b} (90:\R);
\draw[postaction={decorate}, violet] (120:\R)  --node[above left]{c} (150:\R) ;

\draw[postaction={decorate}, violet] (150:\R) --node[left]{d} (180:\R) ;
\draw[postaction={decorate}, violet] (210:\R) --node[left]{c} (180:\R);
\draw[postaction={decorate}, violet] (240:\R) --node[below left]{d} (210:\R);
\draw[postaction={decorate}, cyan] (240:\R) --node[below left]{e} (270:\R);

\draw[postaction={decorate}, cyan] (270:\R) --node[below]{f} (300:\R);
\draw[postaction={decorate}, cyan] (330:\R) --node[below right]{e} (300:\R);

\draw[postaction={decorate}, cyan] (360:\R) --node[below right]{f} (330:\R);

\draw (60:\R) --node[below right]{$c_1$} (180:\R);

\draw  [fill] (1.5,0) circle [radius=0.05] ;
\end{scope}
\end{tikzpicture}
  \caption{}
  \end{subfigure}
  \begin{subfigure}[t]{0.45\textwidth}
  \centering
\begin{tikzpicture}[
    tlabel/.style={pos=0.4,right=-1pt},
    baseline=(current bounding box.center)
    ]
\newdimen\R
\R=2cm
\begin{scope}[very thick,decoration={
    markings,
    mark=at position 0.4 with {\arrow{>}}}
    ] 
\draw[postaction={decorate}, green] (0:\R) --node[right]{a} (30:\R)  ;
\draw[postaction={decorate}, green] (30:\R) --node[above right]{b}  (60:\R);
\draw[postaction={decorate}, green]  (90:\R) --node[above]{a}  (60:\R)  ;
\draw[postaction={decorate}, green] (120:\R) --node[above left]{b} (90:\R);
\draw[postaction={decorate}, violet] (120:\R)  --node[above left]{c} (150:\R) ;

\draw[postaction={decorate}, violet] (150:\R) --node[left]{d} (180:\R) ;
\draw[postaction={decorate}, violet] (210:\R) --node[left]{c} (180:\R);
\draw[postaction={decorate}, violet] (240:\R) --node[below left]{d} (210:\R);
\draw[postaction={decorate}, cyan] (240:\R) --node[below left]{e} (270:\R);

\draw[postaction={decorate}, cyan] (270:\R) --node[below]{f} (300:\R);
\draw[postaction={decorate}, cyan] (330:\R) --node[below right]{e} (300:\R);

\draw[postaction={decorate}, cyan] (360:\R) --node[below right]{f} (330:\R);

\draw (330:\R) --node[below]{$c_2$} (180:\R);

\draw  [fill] (1.5,0) circle [radius=0.05] ;
\end{scope}
\end{tikzpicture}
  \caption{}
  \end{subfigure}
    \caption{Curves correpsonding to Lickorish's generators in figure~\ref{fig:licksgens}, for a genus $3$ surface, cut along curves of the standard presentation of figure~\ref{fig:cutsurf}. (a) The curve $c_1$ of Lickorish's generators. (b) The curve $c_2$ of Lickorish's generators.}\label{fig:lickc12}
\end{figure}
Assuming that $\Pi$ has a subword of the form $[a_i,b_i][a_{i+1},b_{i+1}]$ for $i>1,$ we define $\gamma_i:=a_{i+1}b_{i+1}^{-1}a_{i+1}^{-1}b_ia_i$ and find that the twist along $c_i$ for $i>1$ takes the form
\begin{align}\label{eq:twohandledehnsc2}
T_{c_i}:\begin{cases} 
a_i\mapsto a_i\gamma_i^{-1},\\
b_i\mapsto \gamma_i b_i,\\
a_{i+1}\mapsto \gamma_i a_{i+1},\\
b_{i+1}\mapsto b_{i+1}.
 \end{cases}
\end{align}

Lastly, we need a way to represent a point push around a handle that is not neighbouring the pushed point. Half-twists can be applied to change the ordering of the singular points and boundary components, crosscap transpositions achieve the same for crosscaps, so we are only missing a `transposition of handles.' For two neighbouring handles $[a_1,b_1][a_2,b_2]$ like above, such a transformation can easily be verified to be given by
\begin{align}\label{eq:handletrans}
\begin{cases} 
a_1\mapsto [a_1,b_1] a_2 [a_1,b_1]^{-1},\\
b_1\mapsto [a_1,b_1] b_2 [a_1,b_1]^{-1},\\
a_2\mapsto a_1, \\
b_2\mapsto b_1.
 \end{cases}
\end{align} 
Since this is a type-preserving automorphism of the orbifold fundamental group, theorem~\ref{thm:mcgout} implies that there is a homeomorphism inducing this transformation on the fundamental group. For classical surfaces, this is also a consequence of the well-known Dehn--Nielsen--Baer theorem~\cite{primermcgs}.

With the generators and their action on sets of generators of an orbifold fundamental group derived in this section, the presentation of any MCG can be readily translated into the form they take as automorphisms of $\pi_1(\mathcal{O})$. We will see how to accomplish this for stellate orbifolds in the next section, where we derive a group presentation of a class of MCGs in terms of Dehn twists and half-twists.
\section{The semi-pure mapping class group on the sphere}\label{sec:mixedbraids}
In this section we derive presentations for MCGs of stellate orbifolds to produce the tilings in the next section by exploiting a connection to braid groups $B_n$, where we assume $n\ge 3$. Note that braid groups in the literature do not follow functional notation and if interpreted as maps need to be read from left to right. We follow this tradition, so implementations of the results in this section must be translated to functional notation. We will use both the notion of the braid group as the fundamental group of a configuration space as well as that on the MCG of homeomorphisms of a disk with marked points fixing the boundary in the following~\cite{primermcgs}. 
\begin{definition}
Given a partition $\mathcal{P}$ of $n\in \mathbb{N}$, we define the \emph{semi-pure braid group} $SPB_{\mathcal{P}}$ of type $\mathcal{P}$ to be the subgroup of the braid group $B_n$ that under the canonical morphism to the symmetric group $\Sigma_n$ yields only permutations that respect the partition of $n$, i.e. only permute elements within each set in the partition. We similarly define \emph{semi-pure surface braid groups} and the \emph{semi-pure mapping class group} $\Mod(S_{\mathcal{P}})$ of a surface $S$ with $n$ marked points as the subgroup of the MCG that fixes the partition $\mathcal{P}$ of the marked points. 
\end{definition}
The semi-pure braid group of the plane was studied in~\cite{Manfredini1997}. There is a well-known short exact sequence that relates the braid group of an orientable surface $S$ to the MCG of $S$ and an $n$-times punctured version of $S$, which we denote with $S^*$. This sequence is known as the Birman exact sequence \cite[Theorem $9.1$]{primermcgs} and, as long as $\pi_1(\Hom^+(S))=1$, reads
\begin{align}\label{eq:birmanexseq}
1\to \pi_1(C(S,n)) \to \Mod^+(S^*) \to \Mod^+(S)\to 1,
\end{align}
where $C(S,n)$ is the configuration space of $n$ distinct, unordered points in $S$ and the MCGs involved stem from orientable homeomorphisms. The braid group $B_n$ is isomorphic to $\pi_1(C(\mathbb{C},n))$ \cite[Section $9.1.2$]{primermcgs} and the short exact sequence \eqref{eq:birmanexseq} generalizes to the situation of semi-pure braid groups and the semi-pure MCG, with the same proof, where $\pi_1(C(S,n))$ is replaced by the semi-pure braid group of interest, and $\Mod^+(S^*)$ is replaced by the corresponding semi-pure MCG. One can use \eqref{eq:birmanexseq} to obtain a presentation of $\Mod(S_{\mathcal{P}})$ from one of the braid group for a large class of surfaces. There is also a similar short exact sequence for nonorientable surfaces~\cite[Section $7$]{Szepietowski2008}.
From now on, we restrict to the case where $\Mod(S_{\mathcal{P}})$ denotes the semi-pure MCG on the sphere, as this is the relevant case for stellate orbifolds.

Classical surfaces with negative Euler characteristic satisfy $\pi_1(\Hom^+(S))=1$~\cite{Hamstrom1,hamstrom1965,hamstrom1966}. However, stellate orbifolds do not obey \eqref{eq:birmanexseq}, as the orbifolds underlying topological space is the sphere $S^2$ and therefore $\pi_1(\Hom^+(S^2))=\mathbb{Z}/2\mathbb{Z}$~\cite{Smale1959}.

The standard presentation \eqref{eq:braidpresentation} of $B_n$ is phrased in terms of the generators $\sigma_i$, which translate to half-twist generators in the MCG.  Each $\sigma_i$ crosses the strand in position $i$ in front of the strand in position $i + 1$. In order to find a presentation of the spherical braid group $B_n(S^2)$, we need to append the relation 
\begin{align}\label{eq:Xrel}
X:=\sigma_1...\sigma_{n-1}\sigma_{n-1}...\sigma_1=\id
\end{align}
to the presentation of $B_n$ in terms of the generators $\sigma_i$~\cite{fadell1962,Zariski1936}. Notice that $X$ represents a pure braid.

Consider the element 
\begin{align}
z:=(\sigma_1\cdots \sigma_{n-1})^n\in B_n.
\end{align} 
The infinite cyclic center of the pure braid group $PB_n$ \emph{and} of $B_n$ is generated by $z$~\cite{artinbraids}. Taking the interpretation of $B_n$ as the MCG of the $n$-times punctured disk with homeomorphisms that fix the boundary pointwise, one sees that geometrically, $z$ corresponds to the Dehn twist around the boundary. Following equation $(9.1)$ and figure $9.6$ of \cite{primermcgs}, adding the relation $z=\id$ to the relations of $B_n(S^2)$ turns this group into the corresponding MCG of the sphere $\Mod(S_n)$ with $n$ identical singularities. 

For the semi-pure MCG on the sphere, it is essential to find an expression for the nontrivial extra relations resulting from the topology of the sphere in terms of generators of the pure Braid group. We start with how to express $z$ in terms of the $a_{i,j}$ from \eqref{eq:purebraids}. Using the interpretation of the braid group as the fundamental group of a configuration space on $n$ points, we have~\cite[Chapter $9$]{primermcgs}
\begin{align}\label{eq:relstel1}
z=(a_{1,2}a_{1,3}\cdots a_{1,n})\cdots(a_{n-2,n-1}a_{n-2,n})(a_{n-1,n}).
\end{align} 
In this interpretation, $a_{i,j}$ represents the element of the braid group where the $i$th point is pushed around the $(i+1)$st point. We use the same ideas for $X=\id.$ The original proof in \cite{fadell1962} that $X=\id$ suffices as an extra relation for the full braid group on the sphere hinges on the existence of twists that transforms the analogous relations for the other braids into this one. In terms of the pure braid group, however, we must add back in the missing relations and express them in terms of the pure generators. There are $n$ such relations, and each has an interpretation as the braid where the $i$-th strand passes over all $n-i$ strands to its right, turning back to pass underneath all strands until the first one, and back over the first $i-1$ back to where it started. See \cite[p.\ $194$]{Murasugi1999} for an illustration. The first factor in the expression is $a_{i,i+1}...a_{i,n}$. The second factor, for similar reasons, is given by $a_{i-1,i}a_{i-2,i}...a_{1,i}$ and we obtain 
\begin{align}\label{eq:relstel2}
a_{i,i+1}...a_{i,n}a_{i-1,i}a_{i-2,i}...a_{1,i}, \quad 1\le i\le n
\end{align} 
as a set of equations corresponding to the relation for $X$ in \eqref{eq:Xrel}. This explains how to write the relation corresponding to $X$ in terms of $n$ relations on the generators of the pure braid group. 

We now list all isomorphism classes of stellate groups that are subgroups of $\star 246$ and contain the translations associated to the genus $3$ polygon, which represents the unit cells of the gyroid, primitive and diamond TPMS using the list of commensurate subgroups in \cite{Robins2004}. From now on, we restrict our attention to the orientable MCG for simplicity. A further simplification comes from the fact that the thrice-punctured sphere $S_{3}$ has trivial pure MCG~\cite{primermcgs}, so the only nontrivial elements stem from permutations of the points and every permutation is clearly realized as an orientation preserving homeomorphism. We therefore have
\begin{align}
\Mod(S_{3})\cong \Sigma_3.
\end{align} 

\begin{table}[!htbp]
\centering
\begin{tabular}{cc}
\toprule
Orbifold $\mathcal{O}$& Mapping Class Group $\Mod(\mathcal{O})$\\
\midrule
$246$& Trivial\\
$266$ & $\mathbb{Z}_2$\\
$344$ & $\mathbb{Z}_2$\\
$2223$ & $\Mod(S_{3,1})$\\
$2224$ &$\Mod(S_{3,1})$\\
$2226$ & $\Mod(S_{3,1})$\\
$2233$& $\Mod(S_{2,2})$\\
$2244$& $\Mod(S_{2,2})$\\
$2266$& $\Mod(S_{2,2})$\\
$4444$ & $\Mod(S_4)$\\
$22222$ & $\Mod(S_5)$\\
$22223$ & $\Mod(S_{4,1})$\\
$222222$ & $\Mod(S_6)$\\
$22222222$ & $\Mod(S_8)$\\
\bottomrule
\end{tabular}
\caption{The $14$ mapping class groups of stellate orbifolds of symmetry groups of the gyroid, primitive, and diamond surfaces.}\label{table:stellatemcgs}
\end{table}	

We obtain the list in table~\ref{table:stellatemcgs} of relevant MCGs. Note that the generator of $\mathbb{Z}_2$ in the table is in both cases a half-twist about the two gyration points of the same order.  

Choosing generators pins down one particular presentation, which determines the ordering of the gyration points up to cyclic permutations. The representation of a half twist involving two non-adjacent points is more complicated than \eqref{eq:half-twists}, but can be readily computed using the decomposition of a permutation into adjacent transpositions from section~\ref{sec:conciso}.

A presentation of the subgroup $SPB_{\mathcal{P}}$ of $B_n$ consisting of exactly those elements that leave invariant a partition $\mathcal{P}$ of $n$ of the form $(1...h_1)(h_1+1...h_2)...(h_{m-1}+1...h_m)$ was derived in \cite{Manfredini1997}, following an application of the Reidemeister-Schreier method to glean presentations for subgroups from a presentation of the group. After reexamining theorem $4,$ the most general theorem in the paper, we find that there, $i$ can equal $j$ for the generators $A_{h_ih_j+1}$. Also, the braid relations in equation $(2.2)$ have to be replaced by the usual ones as applicable. Furthermore, in the indexing of $(2.2)$, $i$ can equal $h_{j+1}-2$. Lastly, equation $(2.9)$ needs to be replaced by the relevant equation in \eqref{eq:mixedbraidrep} below. The rather lengthy presentation of the subgroup $\Mod(S_{\mathcal{P}})$ of the MCG that leave invariant the partition $\mathcal{P}$ of $n$ in terms of the generators $\sigma_i$ and the artin generators $a_{i,j}$ of the pure braid group from \eqref{eq:purebraids} is a direct consequence of the above discussion and is as follows.
\begin{minipage}{\textwidth}
\begin{theorem}
Let $\mathcal{P}$ be a partition of $n$ of the form $(1...h_1)(h_1+1...h_2)...(h_{m-1}+1...h_m)$. Then, adding the extra relations \eqref{eq:relstel1} and \eqref{eq:relstel2} to the presentation \eqref{eq:mixedbraidrep} below yields a presentation of the semi-pure MCG $\Mod(S_\mathcal{P})$ for stellate orbifolds.
\rule{\textwidth}{.2em}
\begin{flalign}\label{eq:mixedbraidrep}
&\textbf{Generators} & \nonumber\\
&\sigma_i, \quad i\neq h_t, \quad 1\le t < m, \qquad a_{h_i,h_j+1}, \quad 1\le i\le j< m &\nonumber\\
&\textbf{Relations} &\nonumber\\
&\sigma_i\sigma_j=\sigma_j\sigma_i & \text{if } |i-j|\ge 2,\nonumber\\
&\sigma_i\sigma_{i+1}\sigma_{i}=\sigma_{i+1}\sigma_i\sigma_{i+1},& \text{if } h_j<i\le h_{j+1}-2,\nonumber\\
&a_{h_i,h_j+1}\sigma_k=\sigma_ka_{h_i,h_j+1}& \text{if } h_i-1\neq k\neq h_j+1,\nonumber\\
&a_{h_i,h_j+1}a_{h_k,h_m+1}=a_{h_k,h_m+1}a_{h_i,h_j+1} & \nonumber\\
&\qquad \text{if }h_j<h_k \text{ or if }h_i<h_k<h_m<h_j, &\nonumber\\
&a_{h_i,h_j+1}(\sigma_{h_j+1}a_{h_i,h_j+1}\sigma_{h_j+1})=(\sigma_{h_j+1}a_{h_i,h_j+1}\sigma_{h_j+1})a_{h_i,h_j+1}& \text{if }h_j+1\neq h_{j+1},\nonumber\\ 
&a_{h_i,h_j+1}(\sigma_{h_i-1}a_{h_i,h_j+1}\sigma_{h_i-1})=(\sigma_{h_i-1}a_{h_i,h_j+1}\sigma_{h_i-1})a_{h_i,h_j+1}& \text{if }h_i-1\neq h_{i-1},\nonumber\\ 
&(\sigma_{h_i-1}a_{h_i,h_j+1}\sigma^{-1}_{h_i-1})(\sigma_{h_j+1}a_{h_i,h_j+1}\sigma^{-1}_{h_j+1})= &\nonumber\\
&(\sigma_{h_j+1}a_{h_i,h_j+1}\sigma^{-1}_{h_j+1})(\sigma_{h_i-1}a_{h_i,h_j+1}\sigma^{-1}_{h_i-1})&\\
&\qquad\text{if both } h_i-1\neq h_{i-1} \text{ and } h_j+1\neq h_{j+1},&\nonumber\\
&(\sigma_{h_i-1} a_{h_i,h_j+1}\sigma^{-1}_{h_i-1})a_{h_i,h_k+1}=a_{h_i,h_k+1}(\sigma_{h_i-1}a_{h_i,h_j+1}\sigma^{-1}_{h_i-1})& \nonumber\\
&\qquad \text{if both }h_j<h_k \text{ and } h_i-1\neq h_{i-1}, &\nonumber\\
&a_{h_k,h_j+1}(\sigma_{h_j+1} a_{h_i,h_j+1}\sigma^{-1}_{h_j+1})=(\sigma_{h_j+1}a_{h_i,h_j+1}\sigma^{-1}_{h_j+1})a_{h_k,h_j+1}&\nonumber\\
&\qquad\text{if both }h_i<h_k \text{ and } h_j+1\neq h_{j+1}, &\nonumber\\
\hline
&a_{h_i,h_j}a_{h_i,h_k+1}a_{h_j,h_k+1}=a_{h_i,h_k+1}a_{h_j,h_k+1}a_{h_i,h_j},&\nonumber\\
&a_{h_i,h_j}a_{h_i,h_k+1}a_{h_j,h_k+1}=a_{h_j,h_k+1}a_{h_i,h_j}a_{h_i,h_k+1},&\nonumber\\
&(a_{h_i,h_j}a_{h_i,h_j+1}a^{-1}_{h_i,h_j})(\sigma_{h_j+1}a_{h_k,h_j+1}\sigma^{-1}_{h_j+1})=&\nonumber\\
&(\sigma_{h_j+1}a_{h_k,h_j+1}\sigma^{-1}_{h_j+1})(a_{h_i,h_j}a_{h_i,h_j+1}a^{-1}_{h_i,h_j}) & \nonumber\\
&\qquad \text{if both } h_i<h_k \text{ and }h_j+1\neq h_{j+1}&\nonumber\\
&(a_{h_i,h_j}a_{h_i,h_k+1}a^{-1}_{h_i,h_j})a_{h_j,h_l+1}=a_{h_j,h_l+1}(a_{h_i,h_j}a_{h_i,h_k+1}a^{-1}_{h_i,h_j})&\text{if }h_j\le h_k<h_l.\nonumber 
\end{flalign}
\rule{\textwidth}{.2em}
\end{theorem}
\vspace*{5pt}
\end{minipage}

In the last four equations, whenever $h_j\neq h_{j-1}+1$, i.e., when the symbols $a_{h_i,h_j}$ are not generators of $B_{n,\mathcal{P}}$, we must replace the $a_{h_i,h_j}$ by their expression in terms of the $\sigma_i$ and $a_{h_i,h_j+1}.$ We find the following relations for the pure generators $a_{i,j}.$ Here, $a_{i,j}$ corresponds to the Dehn twist around the simple curve that encircles the marked points $i$ and $j$ (from the point of view of viewing the braid group as the MCG of a disk, with boundary fixed pointwise). We therefore observe that
\begin{align}
a_{i,j}=\sigma_{j-1}...\sigma_xa_{i,x}\sigma_{x}^{-1}...\sigma_{j-1}^{-1}
\end{align}
for $x<j$. Similarly, we establish 
\begin{align}
a_{x,j}=\sigma_{x}^{-1}...\sigma_{z-1}^{-1}a_{z,j}\sigma_{z-1}...\sigma_{x}
\end{align}
for $j>z>x$.
Note, again, that in contrast to the algebraic braid group given in terms of abstract group elements, whose multiplication is read from left to right, MCG elements act from right to left. This leads to $(\sigma_1...\sigma_i)^{-1}$ acting as $T_{\sigma_1}^{-1}...T_{\sigma_i}^{-1}$, where the $T_{\sigma_j}$ are the associated twists in the MCG of a disk corresponding to the twists $\sigma_j$ in the braid group. 

Lastly, to find a presentation for a semi-pure Braid group we assumed that the elements of the same type are grouped together as neighbours as in the partition of $n$ above in \eqref{eq:mixedbraidrep}. If this is not the case, then we must conjugate all elements of the semi-pure braid group with a braid $b$ that maps to an element in the symmetric group realizing the appropriate partition from the one at hand. The ordering of the points within a part in the partition is arbitrary. The element $b$ is then decomposed into neighbouring transpositions for which \eqref{eq:half-twists} gives the representation as an automorphism. 

\section{Tilings with Stellate Symmetry Groups in $\star 246$}\label{sec:steltileex}
\begin{figure}[!htbp]
\imagewidth=0.31\textwidth
\captionsetup[subfigure]{width=0.9\imagewidth,justification=raggedright}
  \begin{subfigure}[t]{0.31\textwidth}
  \centering
    \includegraphics[width=\textwidth, height=\textwidth]{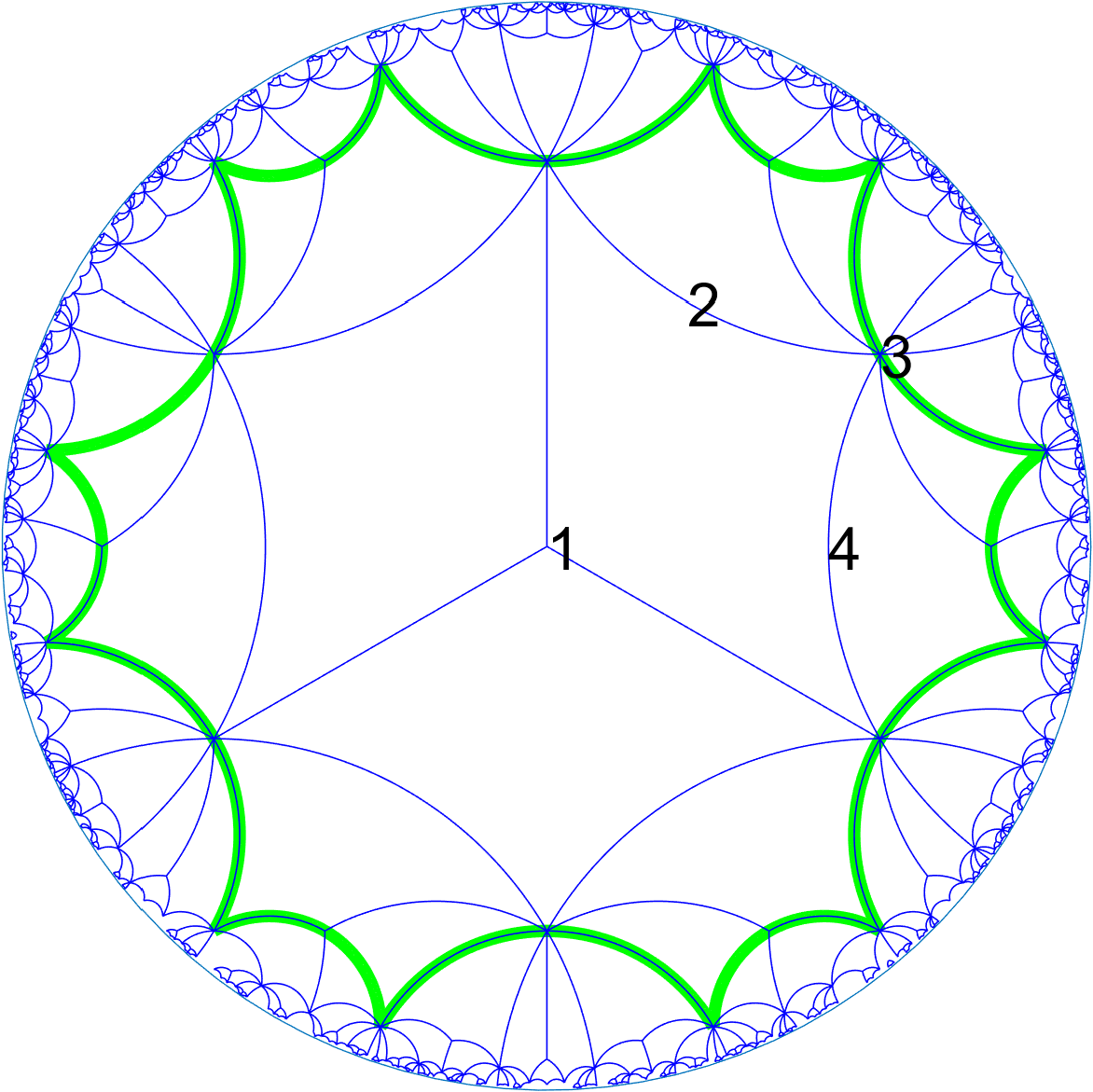}
    \caption{Tiling number $1$.}\label{fig:tiling3232a}
  \end{subfigure}
  \begin{subfigure}[t]{0.31\textwidth}
  \centering
    \includegraphics[width=\textwidth, height=\textwidth]{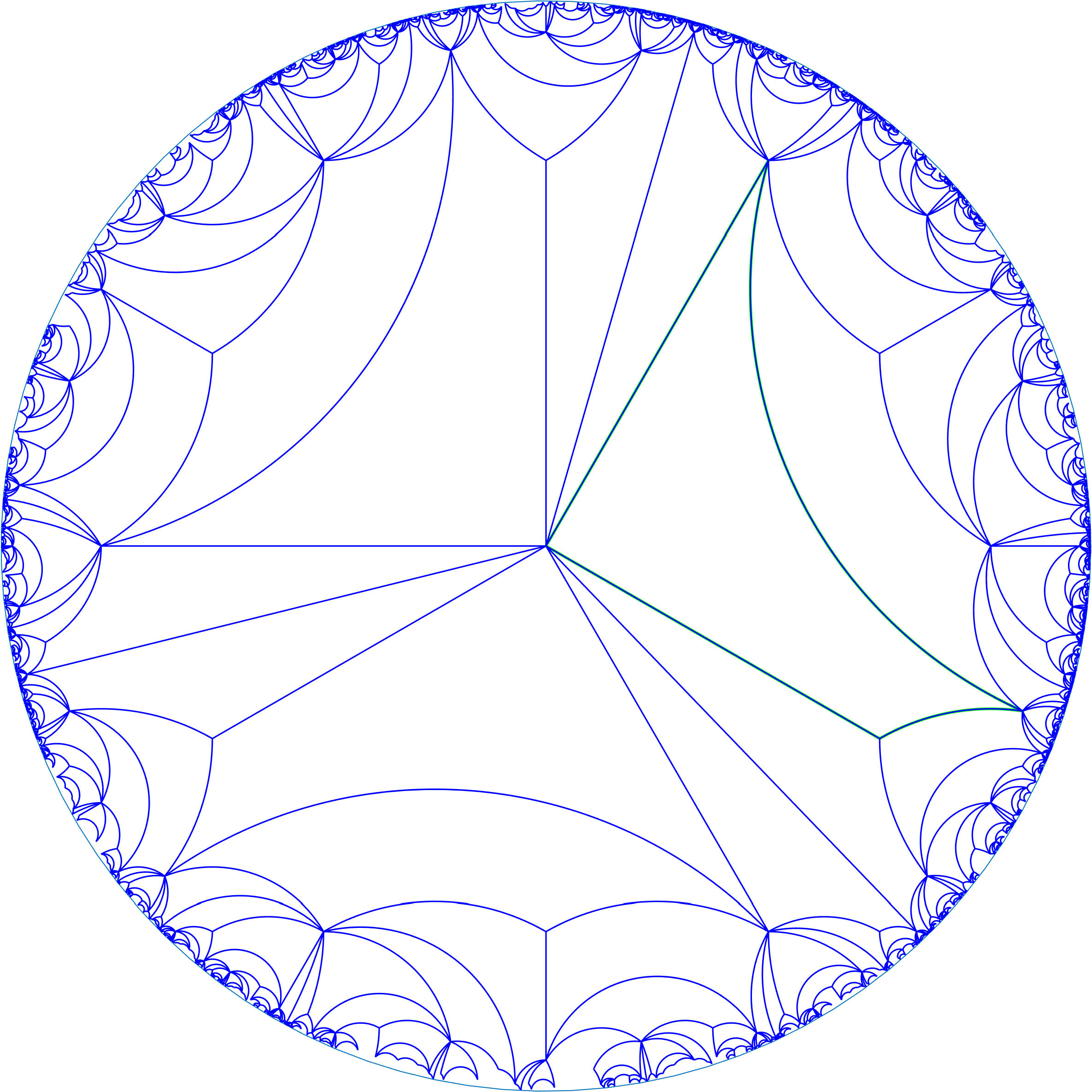}
  \caption{Tiling number $2$.}
  \label{fig:tiling3232b}
  \end{subfigure}
  \begin{subfigure}[t]{0.31\textwidth}
  \centering
    \includegraphics[width=\textwidth, height=\textwidth]{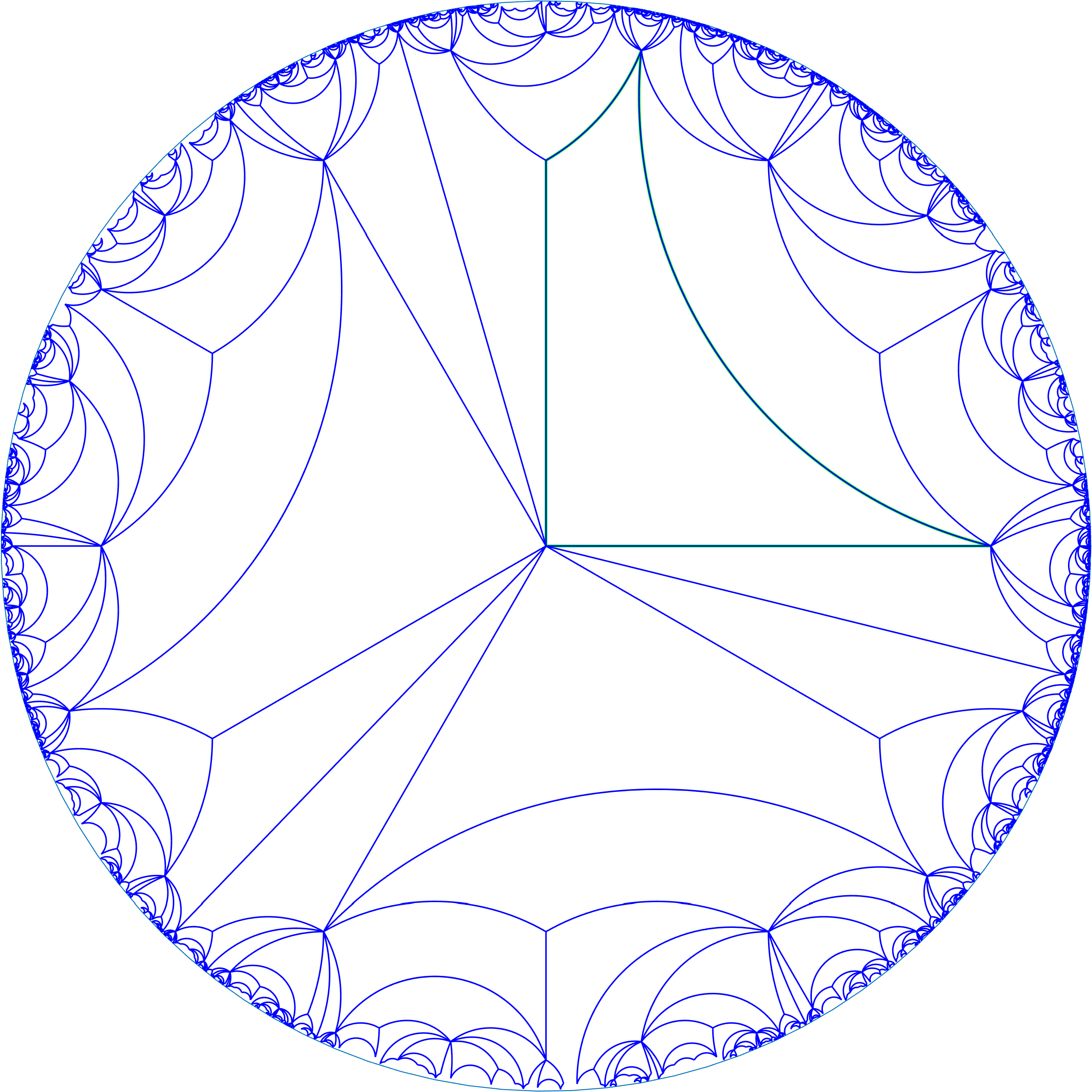}
    \caption{Tiling number $3$.}
  \end{subfigure}
    \\
  \begin{subfigure}[t]{0.31\textwidth}
  \centering
    \includegraphics[width=\textwidth, height=\textwidth]{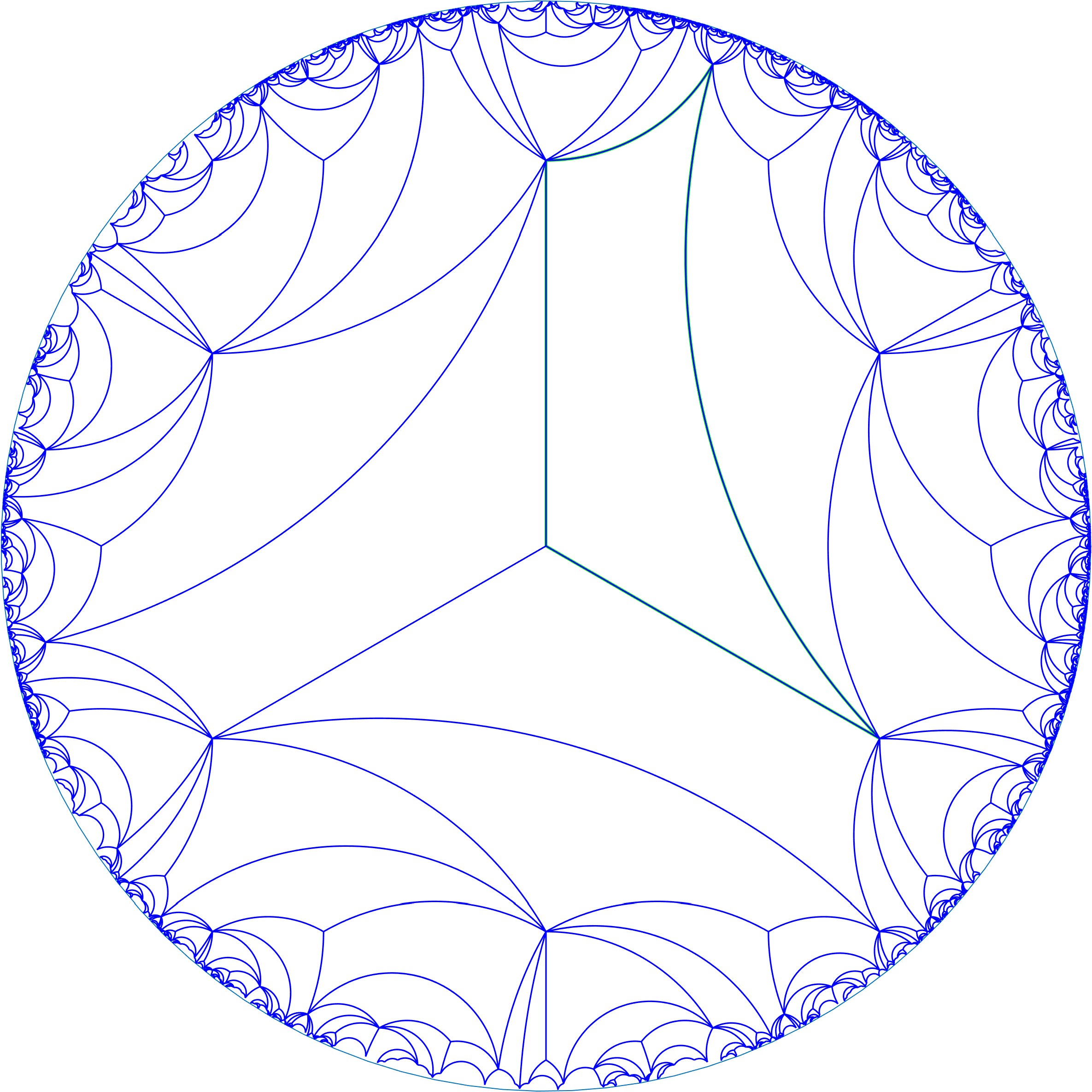}
    \caption{Tiling number $4$.}
  \end{subfigure}
  \begin{subfigure}[t]{0.31\textwidth}
  \centering
    \includegraphics[width=\textwidth, height=\textwidth]{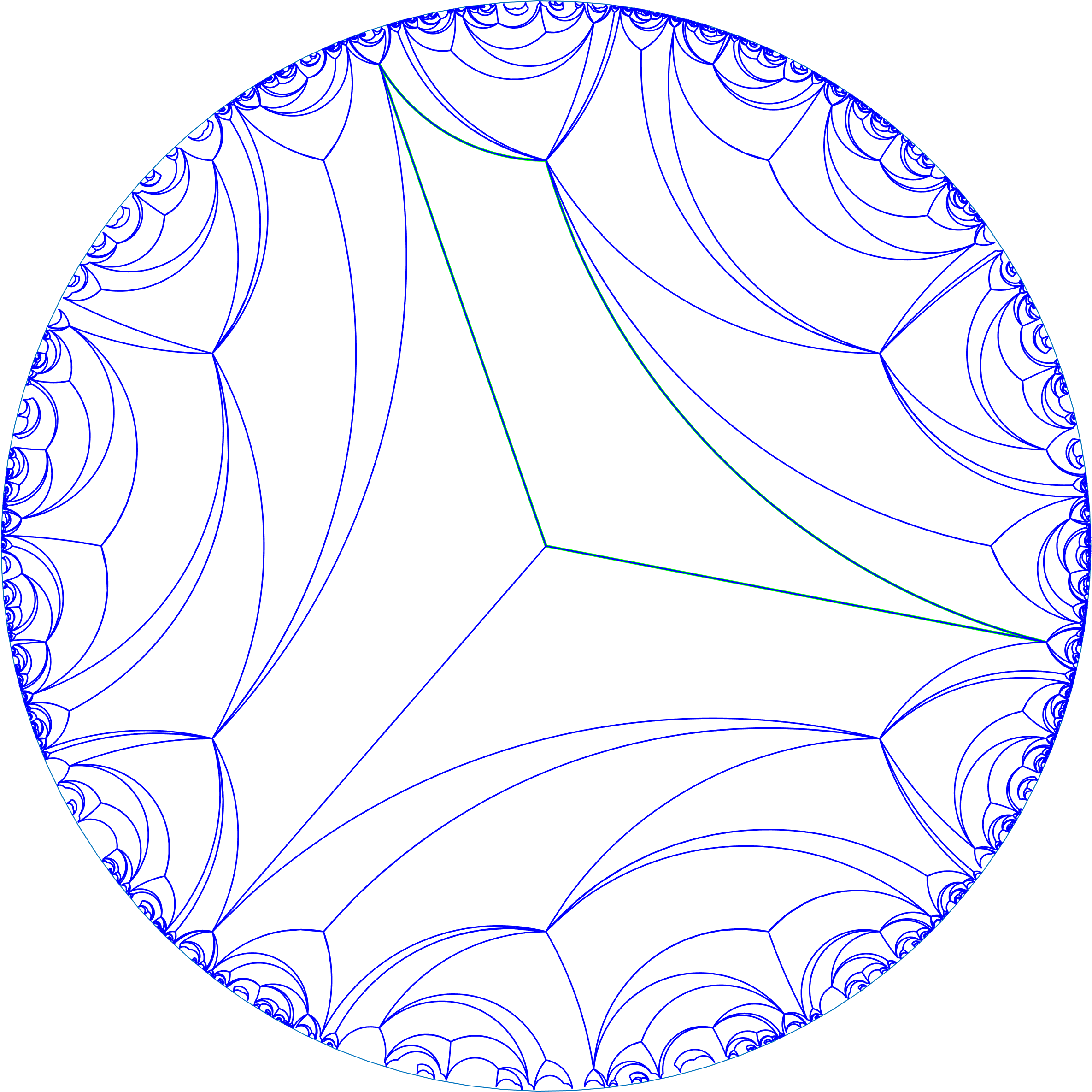}
  \caption{Tiling number $5$.} \label{fig:tiling3232e}
  \end{subfigure}
  \begin{subfigure}[t]{0.31\textwidth}
  \centering
    \includegraphics[width=\textwidth, height=\textwidth]{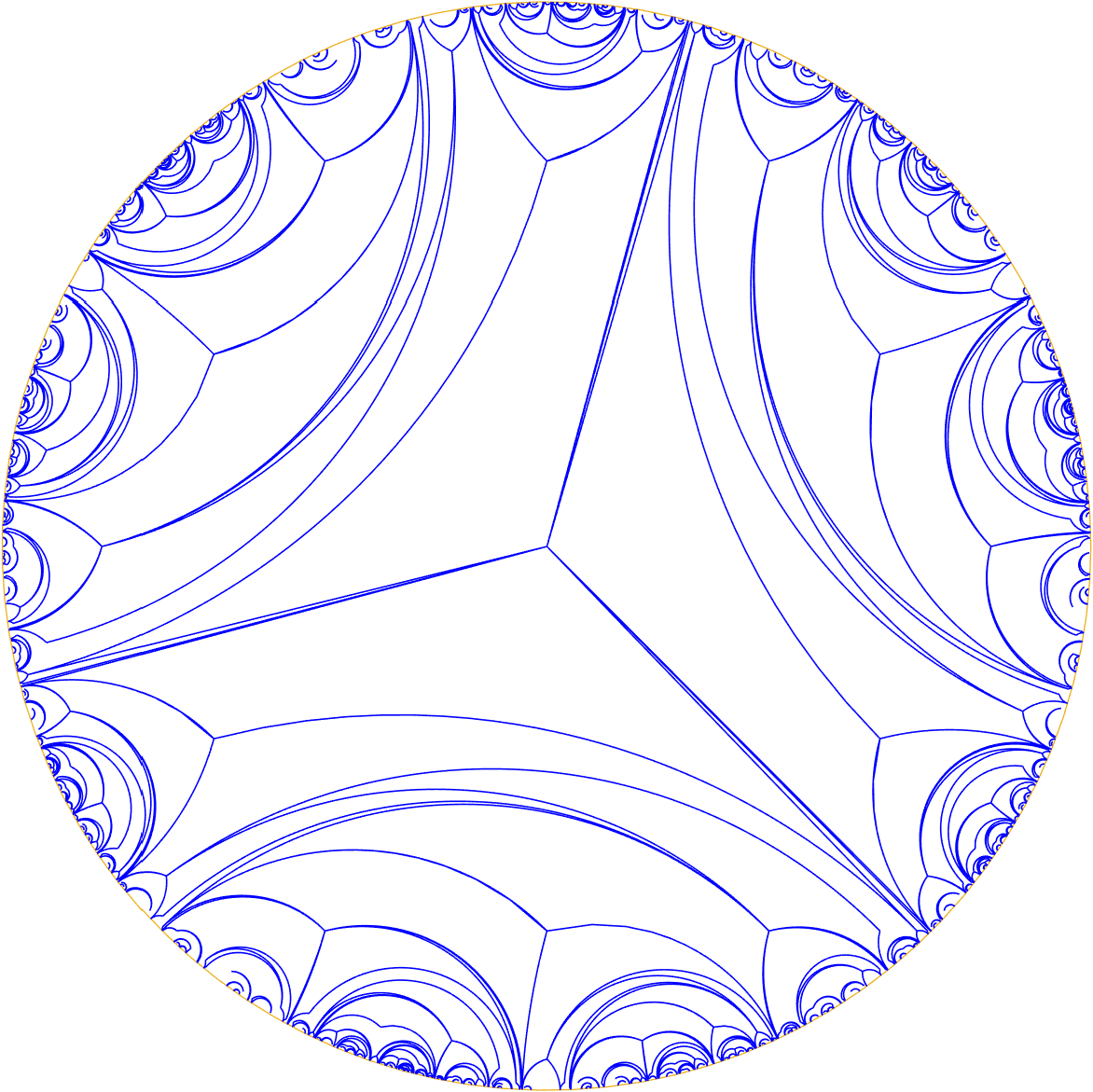}
    \caption{Tiling number $40$.}
  \end{subfigure}    
  \caption{Isotopically distinct fundamental tilings with symmetry group $2233$ and the same D-symbol. All tilings fit onto the genus-$3$ surface obtained by identifying opposite edges of the dodecagon in green in figure~\ref{fig:tiling3232a}. Starting from the top left, which shows a simplest starting tiling, the tilings to the right are a result of `twisting' the tiling, and numbers $1$ to $5$ and $40$ in our enumeration, respectively. Figure~\ref{fig:tiling3232a} shows the placements of the generators, figures~\ref{fig:tiling3232b} to~\ref{fig:tiling3232e} show the boundary of the fundamental tile on which the generators act, in green.}\label{fig:tiling3232}
\end{figure}

This section is dedicated to the explicit construction of tilings of $\mathbb{H}^2$ that are commensurate with the candidate TPMS, the gyroid, primitive, and diamond surfaces.
Using the computer programming language GAP and in particular the Knuth-Bendix package KBMAG~\cite{kbmag}, we have successfully implemented the results of the previous sections for general semi-pure MCGs for stellate orbifolds to produce a list of MCG elements ordered by word length from a presentation of the MCG. Together with how the generators of the MCG act on the generators of the symmetry group as derived in section~\ref{sec:conciso}, we produce a list of sets of generators of the symmetry group in $\Iso(\mathbb{H}^2)$, starting from a chosen starting set of generators. Then, from the placements of generators, we produced hyperbolic tilings from decorations of the corresponding orbifold derived from D-symbols encoding their equivariant equivalence class. 
The sequence of pictures of tilings were produced using MATLAB and are part of an exhaustive enumeration of isotopy classes of tilings with a given compatible symmetry group. We concentrate on the most challenging cases to highlight the approach. Every tiling presented here is drawn on the hyperbolic plane and is commensurate with the genus-$3$ hyperbolic surface $S$ obtained by identifying opposite edges of the dodecagon illustrated in figures~\ref{fig:tiling3232a} and \ref{fig:disksurface}(b). 

One of the biggest challenges in working with the semi-algorithms provided in KBMAG is that sometimes small changes in parameters can lead to very different results~\cite{kbmag,Holt2017}. For the full MCG of the genus $3$ surface, there is a presentation known as the Gervais presentation~\cite{GERVAIS2001}, appealing for applications because its generators are symmetric and supported in small subsurfaces. Unfortunately, we did not succeed in choosing a set of parameters for KBMAG to solve the word problem with this presentation.  Nevertheless, keeping in mind that for an unambiguous enumeration it is sometimes necessary to enumerate cosets of the MCG w.r.t.\ finite groups, the Knuth-Bendix package seems like the most comprehensive option available to date. 

We focus exclusively on decorations given by piece-wise geodesics on the orbifold $\mathcal{O}$ in its induced metric from $\mathbb{H}^2$. Any other tiling can be isotoped to produce such a decoration of the orbifold. In general, for more complicated tilings, one has to insert additional points not at vertices of the tiling at which we allow breaks of geodesics. Viewing tilings as piece-wise geodesic decorations on metric orbifolds lets us specify any tiling by a finite number of points $P\subset \mathbb{H}^2$ along with an adjacency matrix $A$ specifying which points are neighbours. This leads to a data structure similar to other known structures for symmetric tilings, such as that for periodic structures in CGAL~\cite{cgal}. 
The idea of the ensuing data structure hinges on finding $P$ for a piece-wise geodesic tiling $\mathcal{T}$ invariant under $\Gamma=\pi_1(\mathcal{O})$ so that when connecting points in $P$ according to $A$, one obtains a graph $G$ such that $\Gamma G$ encompasses all edges in $\mathcal{T}.$ To do so, find a fundamental domain $D$ for $\Gamma$ and a set $P_0$ containing a unique copy of each vertex and corner of the tiling $\mathcal{T}$ in $D$. Then, for every point $p\in P_0$, find the set $E_p$ of all edges incident to $p$ and find its other endpoint $p'$, and mark each edge by the unique element $f\in\Iso(\mathbb{H}^2)$ s.t. there exists a unique point $p_0\in P_0$ with $f(p_0)=p'$. An edge in $E_p$ can then be been as having a start and an endpoint $p,p_0\in P_0$, resp., and the marking $f$. Then, the set $P$ can be found as the (not necessarily disjoint) union of $P_0$ and, for all $p\in P_0$, the isometries of the edges in $E_p$ applied to their endpoint.

\begin{figure}[!htbp]
\imagewidth=0.31\textwidth
\captionsetup[subfigure]{width=0.9\imagewidth,justification=raggedright}
  \begin{subfigure}[t]{0.31\textwidth}
  \centering
    \includegraphics[width=\textwidth, height=\textwidth]{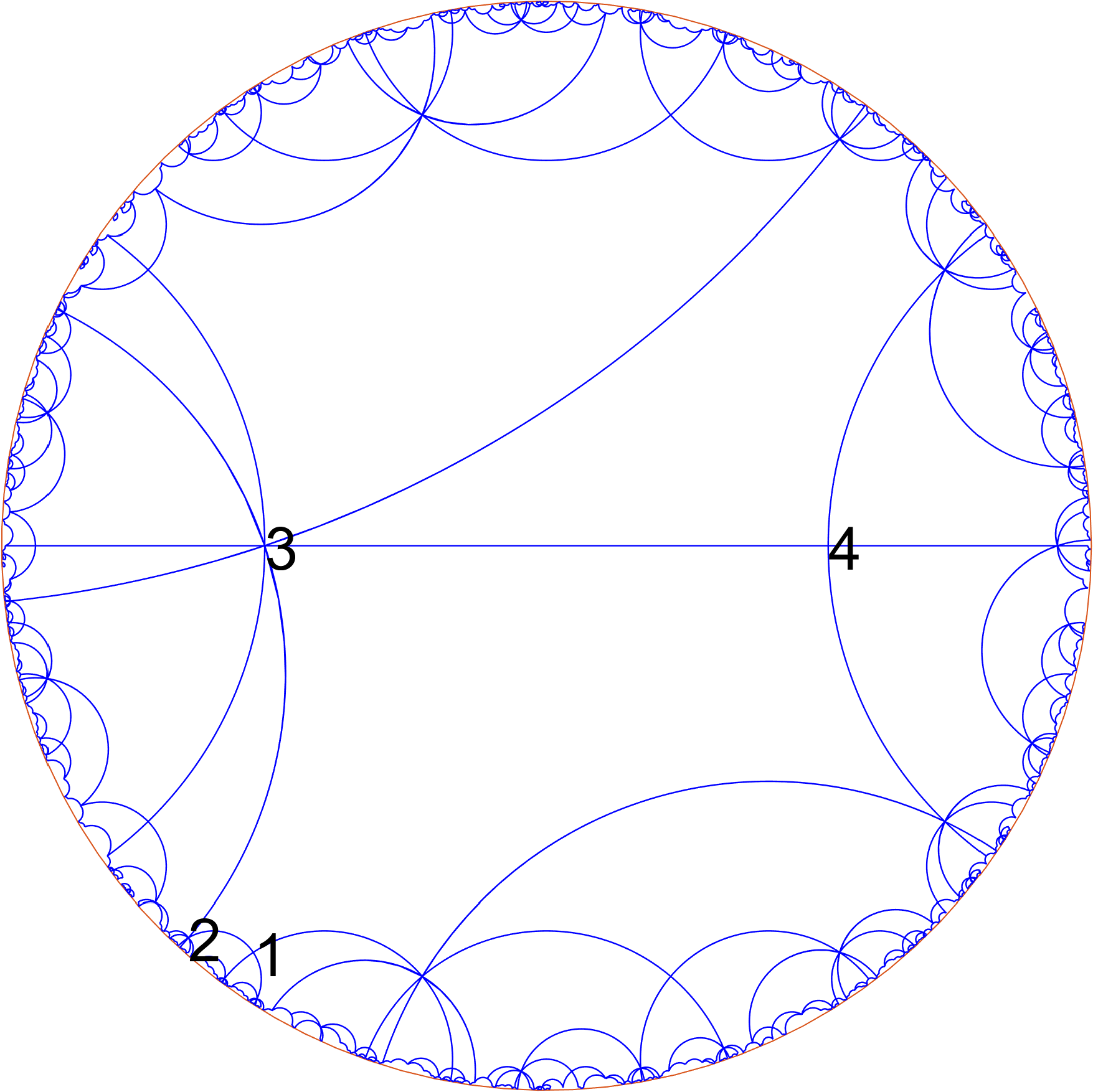}
    \caption{Tiling number $20$.}\label{fig:tiling4444b}
  \end{subfigure}
  \begin{subfigure}[t]{0.31\textwidth}
  \centering
    \includegraphics[width=\textwidth, height=\textwidth]{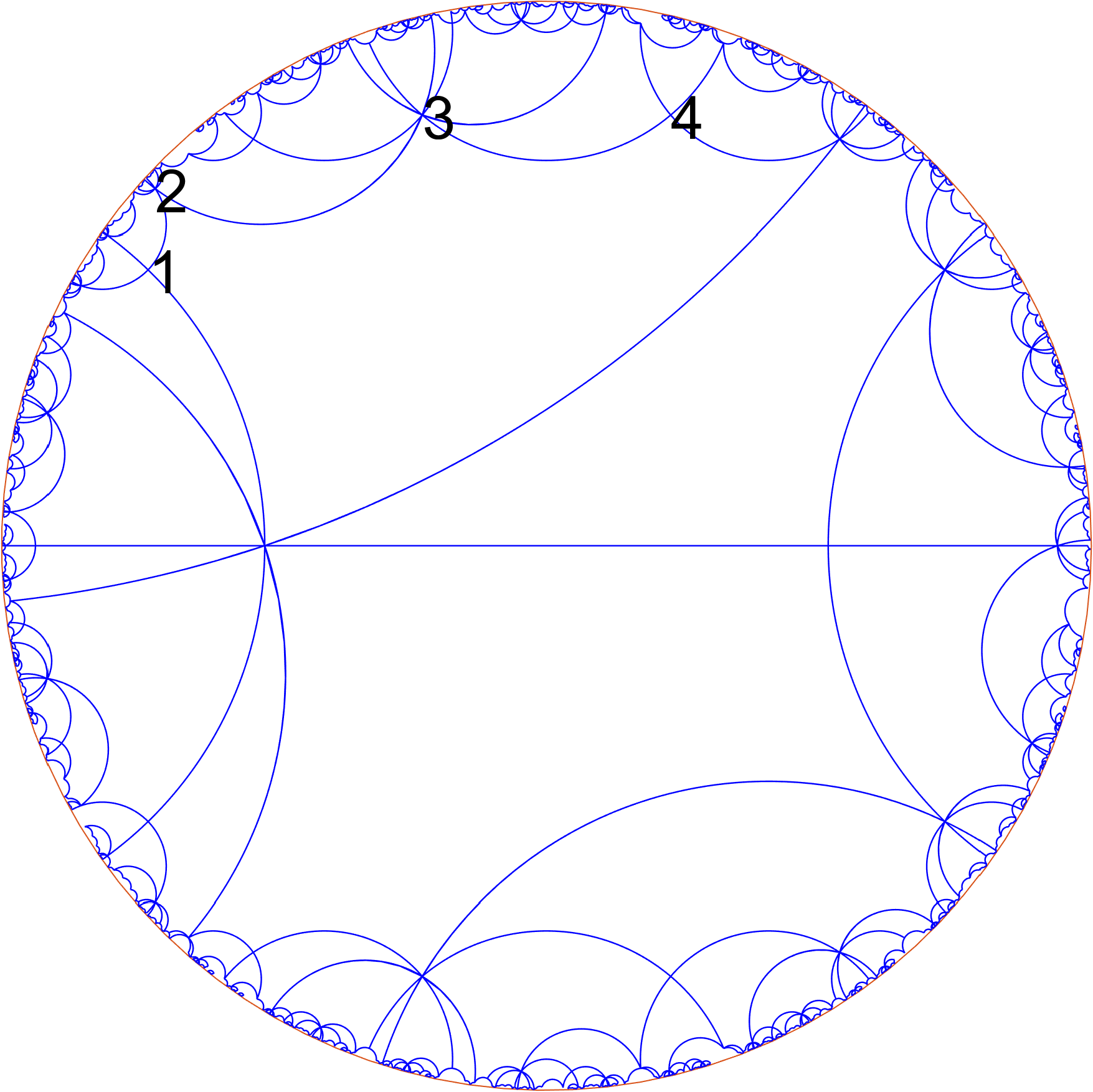}
    \caption{Tiling number $40$.}\label{fig:tiling4444vertsa}
  \end{subfigure}
  \begin{subfigure}[t]{0.31\textwidth}
  \centering
    \includegraphics[width=\textwidth, height=\textwidth]{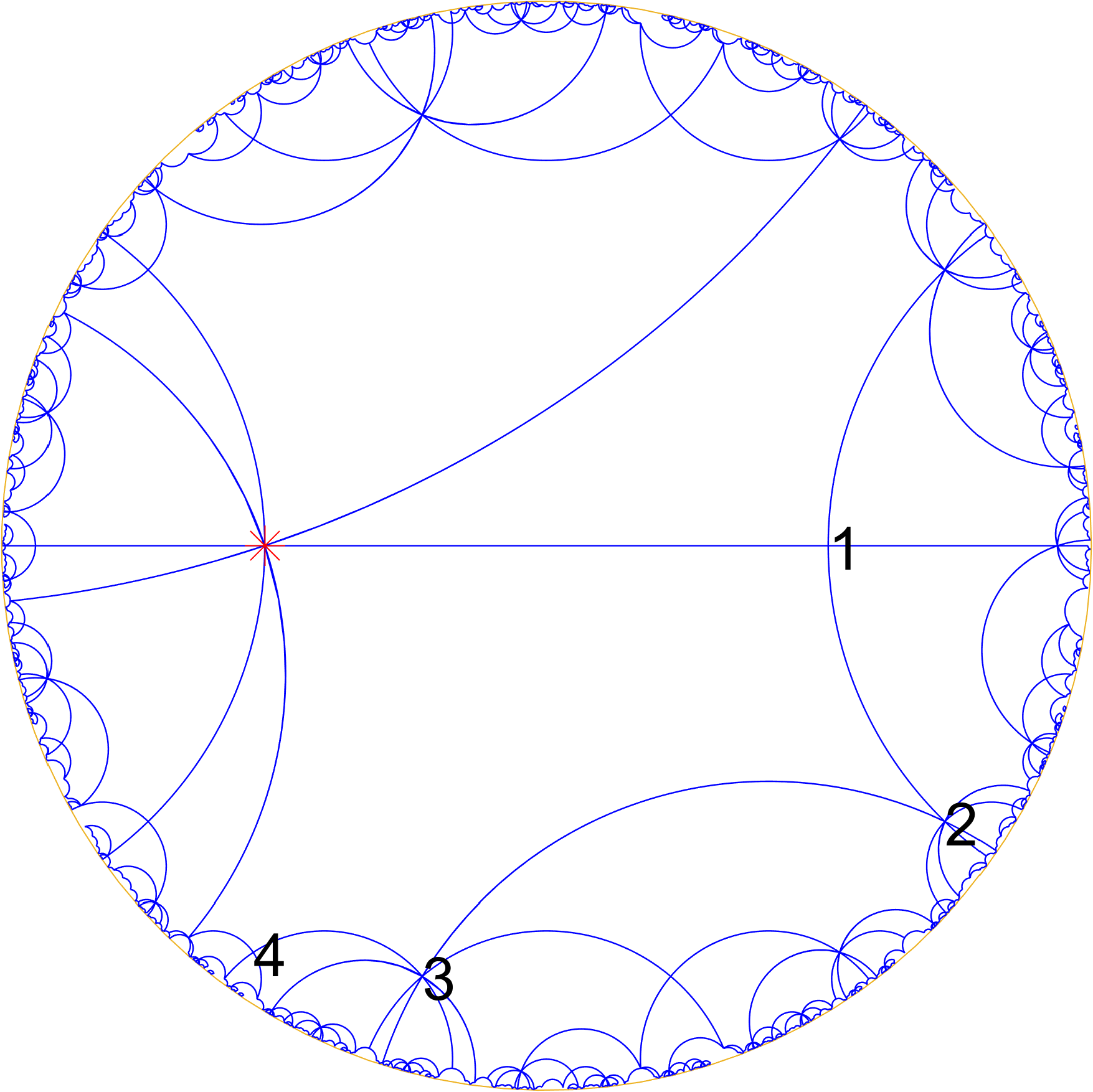}
    \caption{Tiling number $40$.}\label{fig:tiling4444vertsb}
  \end{subfigure}
  \caption{Tiling no. 40 with symmetry $4444$. Image (a) shows a placement of generators different from (b) and both yield identical tilings. The placements of the generators here are a permutation of the generators in the other tiling. The numbering of symmetry elements on the same tile is shown in (c), after conjugation by rotation around the indicated symmetry point. The decoration is too simple to detect the differences, so the tilings show up as the same.}\label{fig:tiling4444verts}
\end{figure}

The enumeration presented here is, as mentioned above, an enumeration of coloured tilings, where each edge orbit is given a different colour. Here, it is possible that two tilings appear the same because we have not gotten rid of ambiguities as a result of internal symmetries of the tilings, i.e. realizations of the equivariant equivalence classes that exhibit more symmetries than shown here~\cite[Proposition $1$]{BenMyf1}. Consider, for example, figure~\ref{fig:tiling4444b} and the identical tiling shown in figure~\ref{fig:tiling4444vertsa}. Observe that there are $6$ points of increased symmetry on the boundary of each tile in the associated combinatorial class of tilings. A different colouring of the edge connecting the vertices labelled $1$ and $2$ in each tiling would distinguish the two, as comparing the placements of generators in figures~\ref{fig:tiling4444b} and \ref{fig:tiling4444vertsb} illustrates. In general, all ambiguities in the enumeration are a result of an automorphism of the graph representing the D-symbol, which can be efficiently checked~\cite[Section $7$]{BenMyf1}. Note that, for EPINET, it may well be of interest to chemists to not eliminate this ambiguity when looking at structures with symmetrically distinct constituent parts, as illustrating by two different edge colours in figure~\ref{pnas}.

\begin{figure}[!htbp]
\imagewidth=0.31\textwidth
\captionsetup[subfigure]{width=0.9\imagewidth,justification=raggedright}
  \begin{subfigure}[t]{0.31\textwidth}
  \centering
    \includegraphics[width=\textwidth, height=\textwidth]{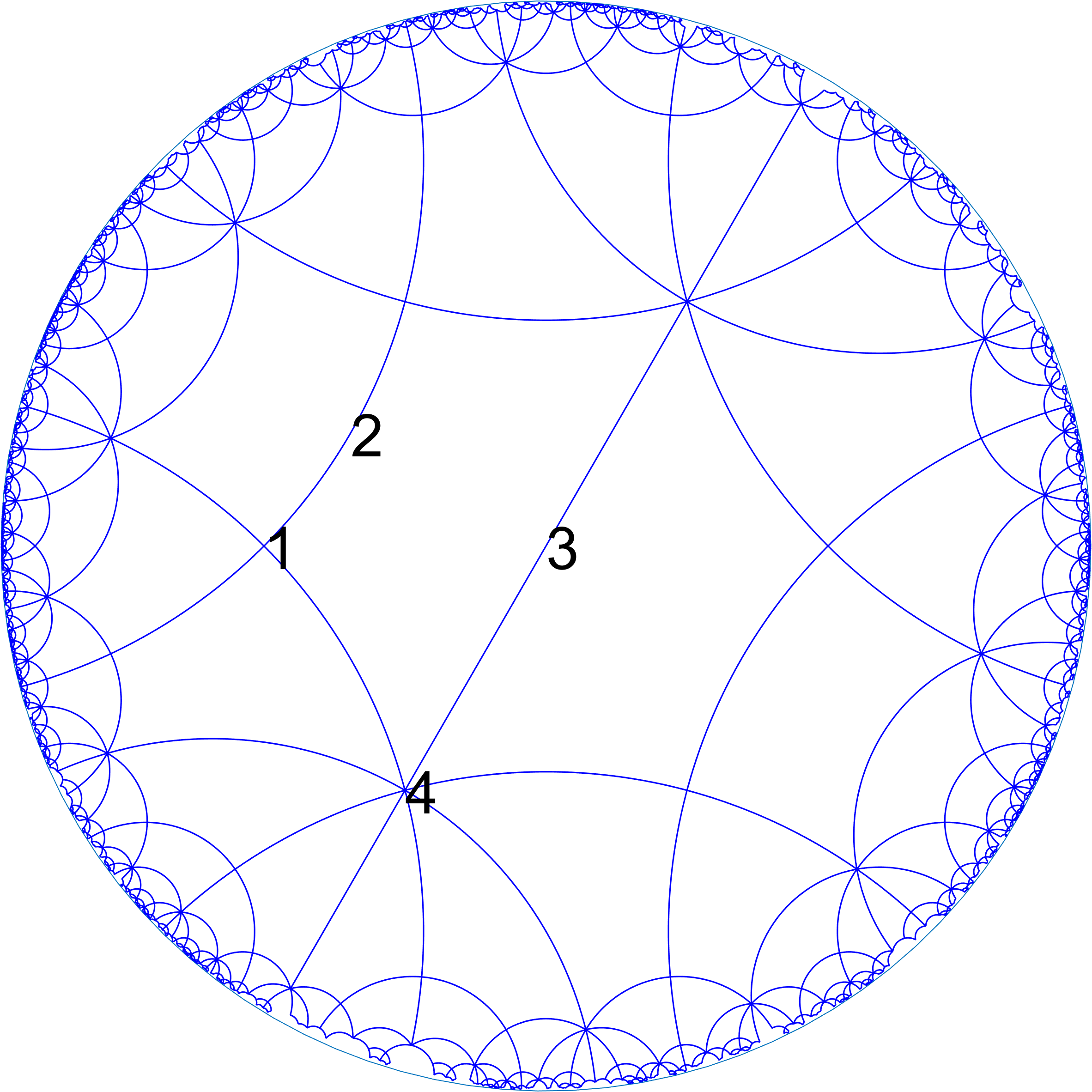}
    \caption{Tiling number $1$.}\label{fig:tiling2224b}
  \end{subfigure}
  \begin{subfigure}[t]{0.31\textwidth}
  \centering
    \includegraphics[width=\textwidth, height=\textwidth]{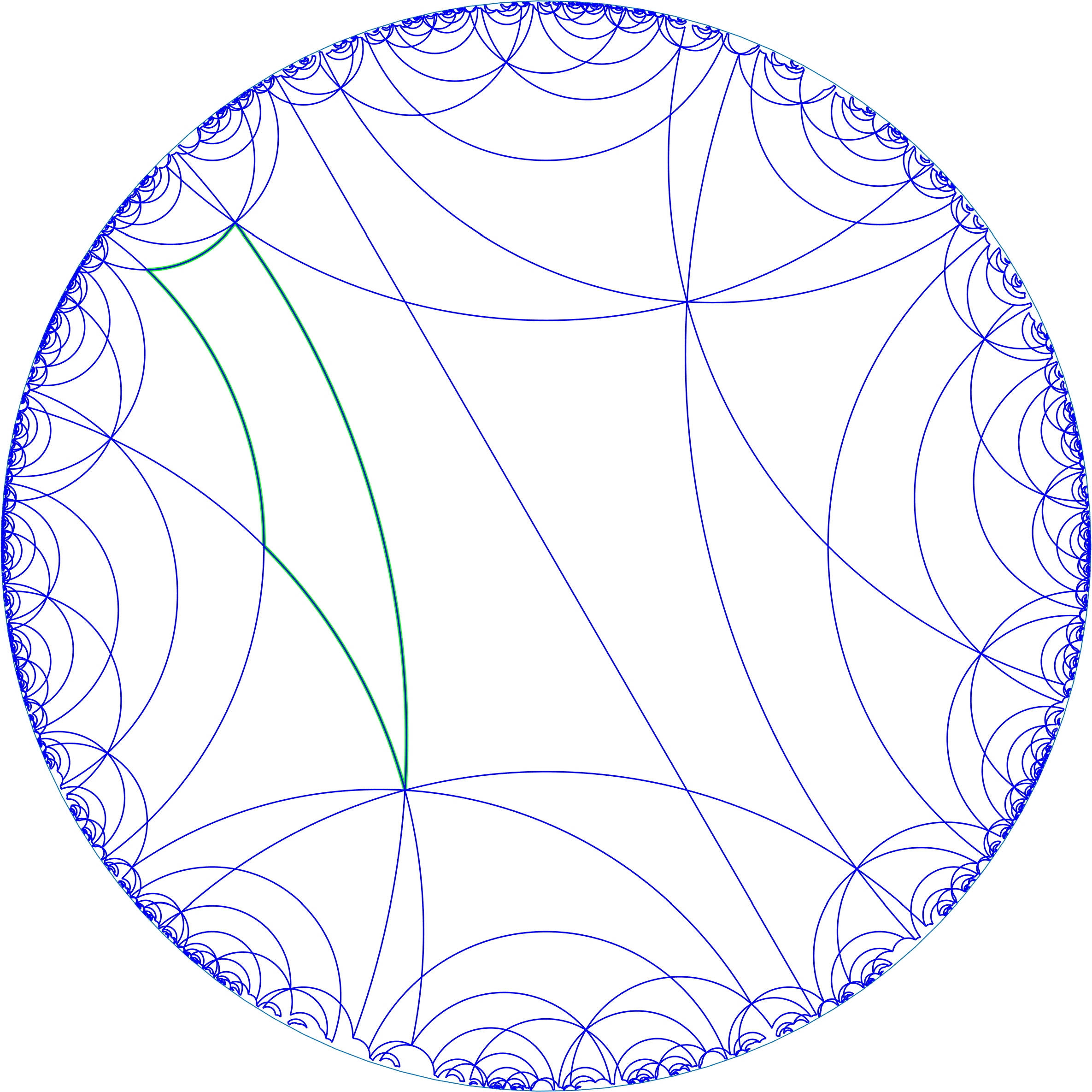}
    \caption{Tiling number $4$.}\label{fig:tiling2224c}
  \end{subfigure}
  \begin{subfigure}[t]{0.31\textwidth}
  \centering
    \includegraphics[width=\textwidth, height=\textwidth]{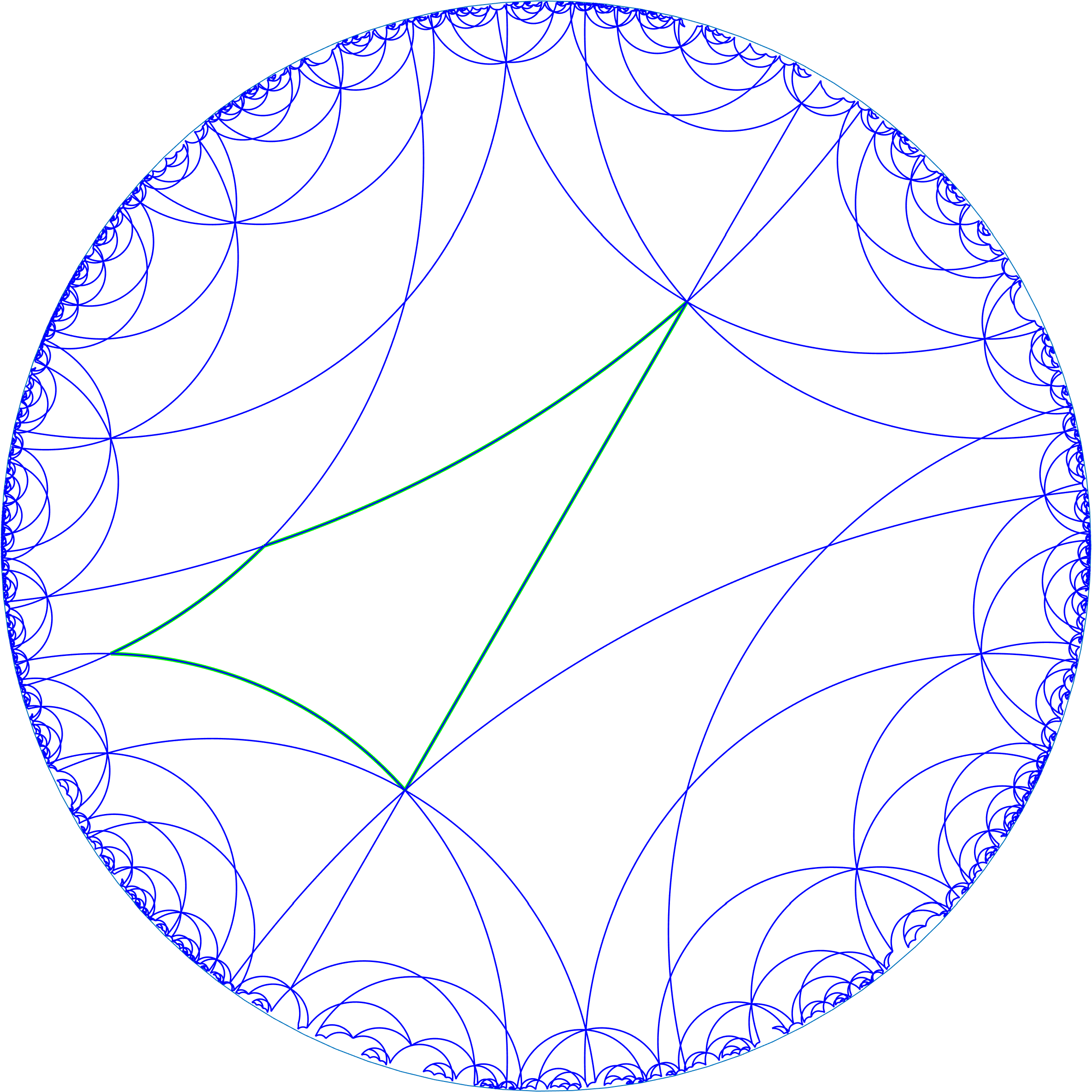}
    \caption{Tiling number $5$.}
  \end{subfigure}
  \\
  \begin{subfigure}[t]{0.31\textwidth}
  \centering
    \includegraphics[width=\textwidth, height=\textwidth]{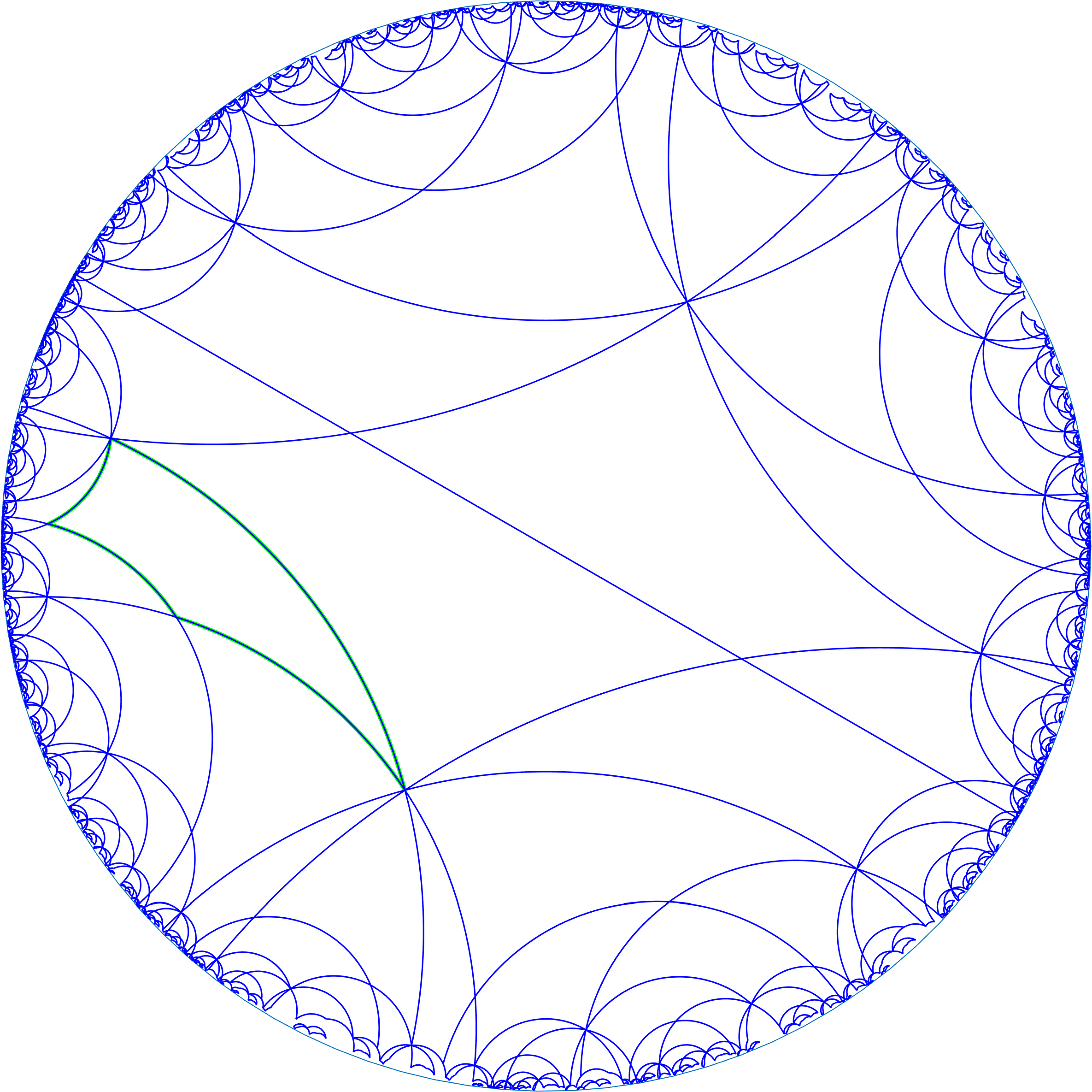}
    \caption{Tiling number $6$.}
  \end{subfigure}
  \begin{subfigure}[t]{0.31\textwidth}
  \centering
    \includegraphics[width=\textwidth, height=\textwidth]{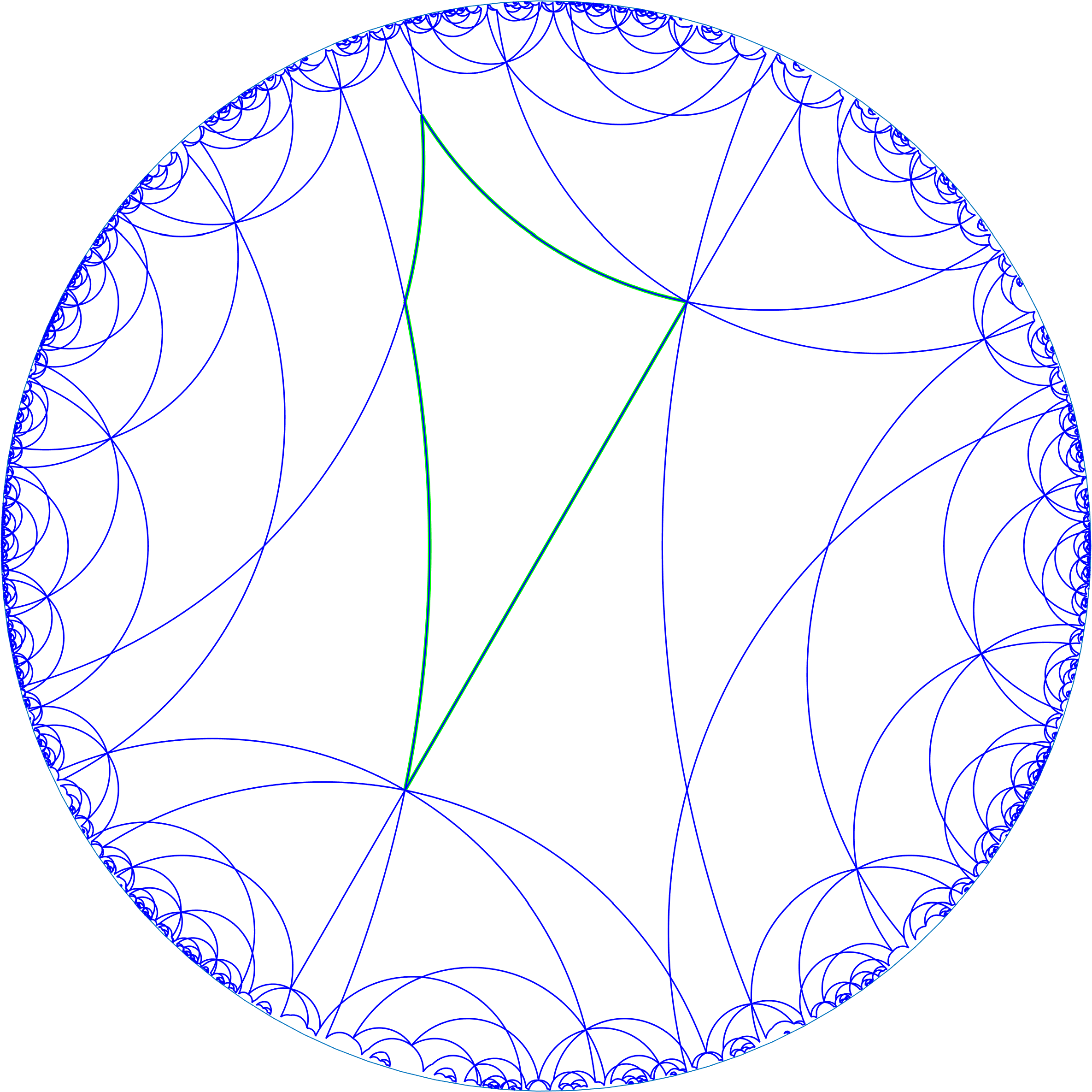}
    \caption{Tiling number $7$.}
  \end{subfigure}
  \begin{subfigure}[t]{0.31\textwidth}
  \centering
    \includegraphics[width=\textwidth, height=\textwidth]{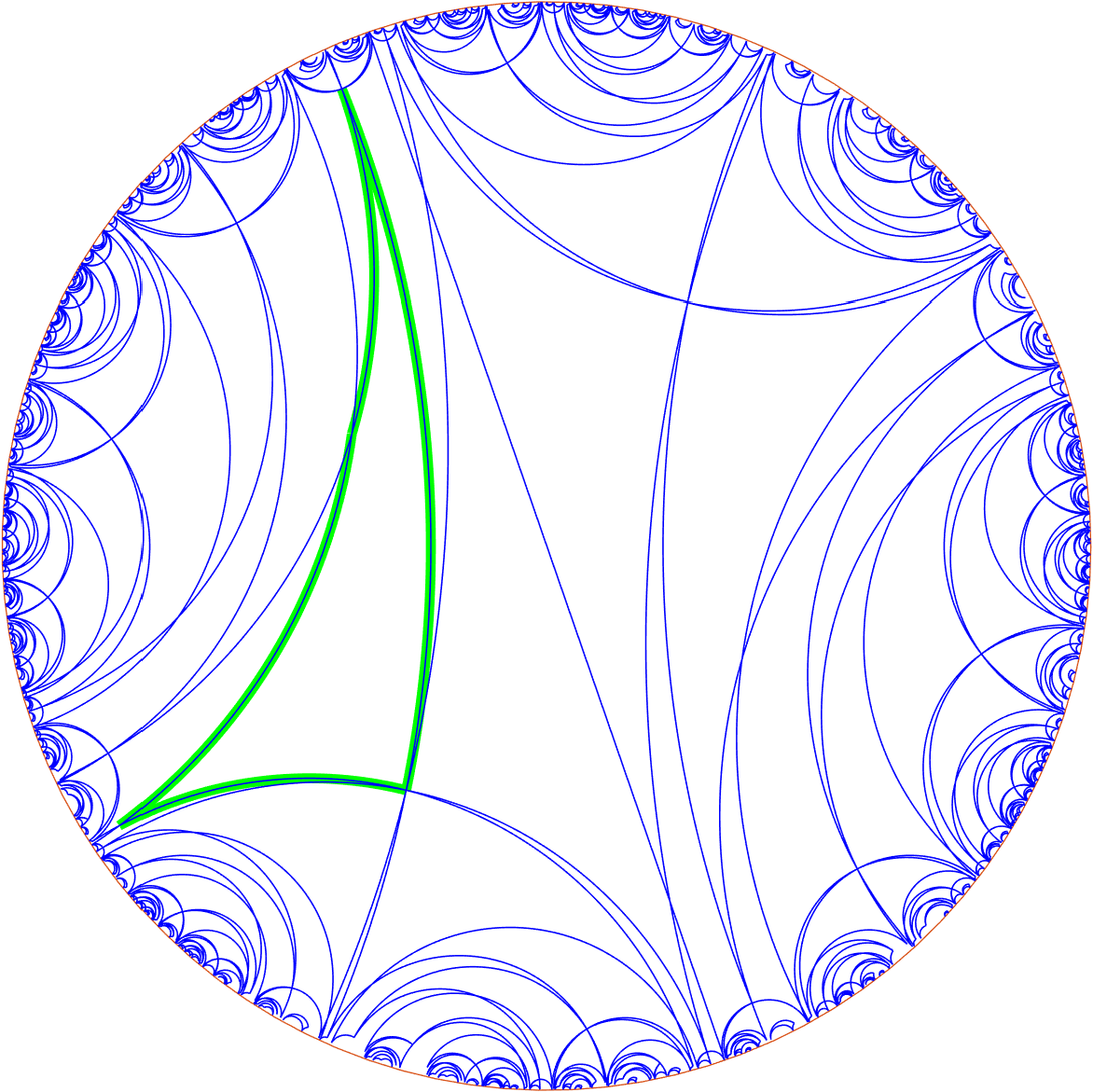}
    \caption{Tiling number $40$.}
  \end{subfigure}  
  \caption{Isotopically distinct fundamental tilings with symmetry group $2224$ and equivalent combinatorial structure. They are numbers $1$, $4$ through $7$, and $40$ in our enumeration, respectively. All tilings are commensurate with the genus-$3$ dodecagon from figure~\ref{fig:tiling3232a}.}\label{fig:tiling2224}
\end{figure}

\begin{figure}[!htbp]
\imagewidth=0.31\textwidth
\captionsetup[subfigure]{width=0.9\imagewidth,justification=raggedright}
  \begin{subfigure}[t]{0.31\textwidth}
  \centering
    \includegraphics[width=\textwidth, height=\textwidth]{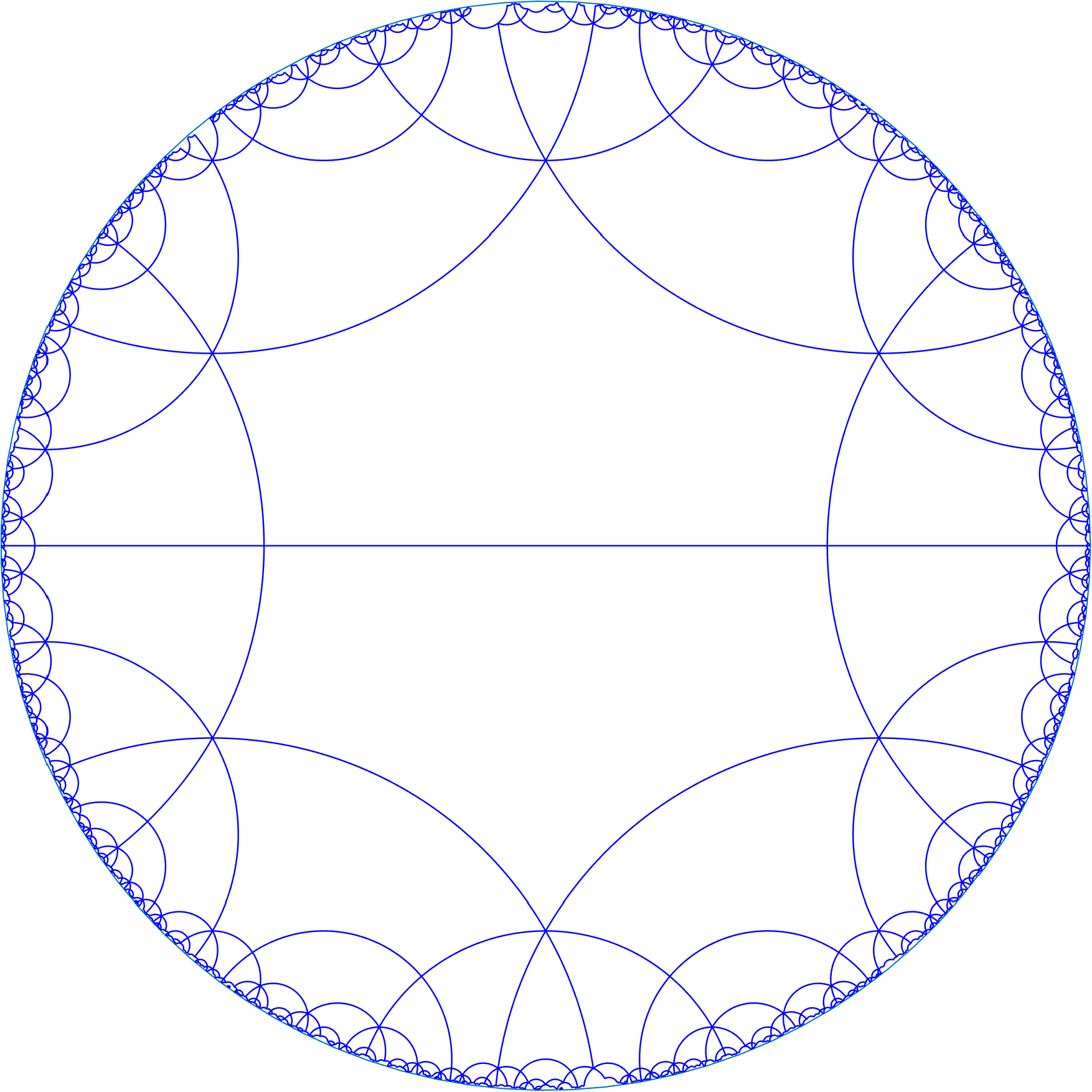}
    \caption{Tiling number $1$.}\label{fig:tiling2222254a}
  \end{subfigure}
  \begin{subfigure}[t]{0.31\textwidth}
  \centering
    \includegraphics[width=\textwidth, height=\textwidth]{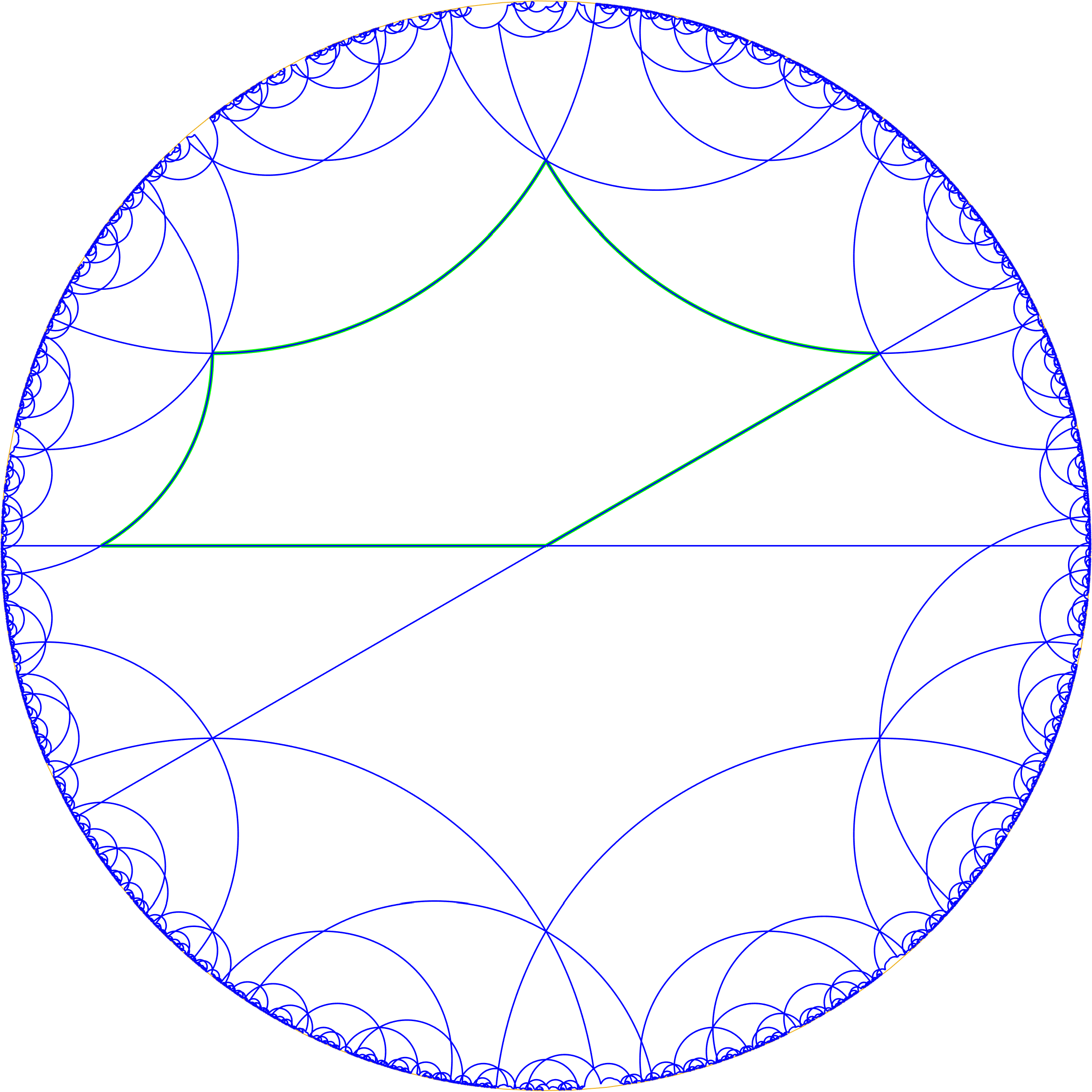}
    \caption{Tiling number $2$.}\label{fig:tiling2222254b}
  \end{subfigure}
  \begin{subfigure}[t]{0.31\textwidth}
  \centering
    \includegraphics[width=\textwidth, height=\textwidth]{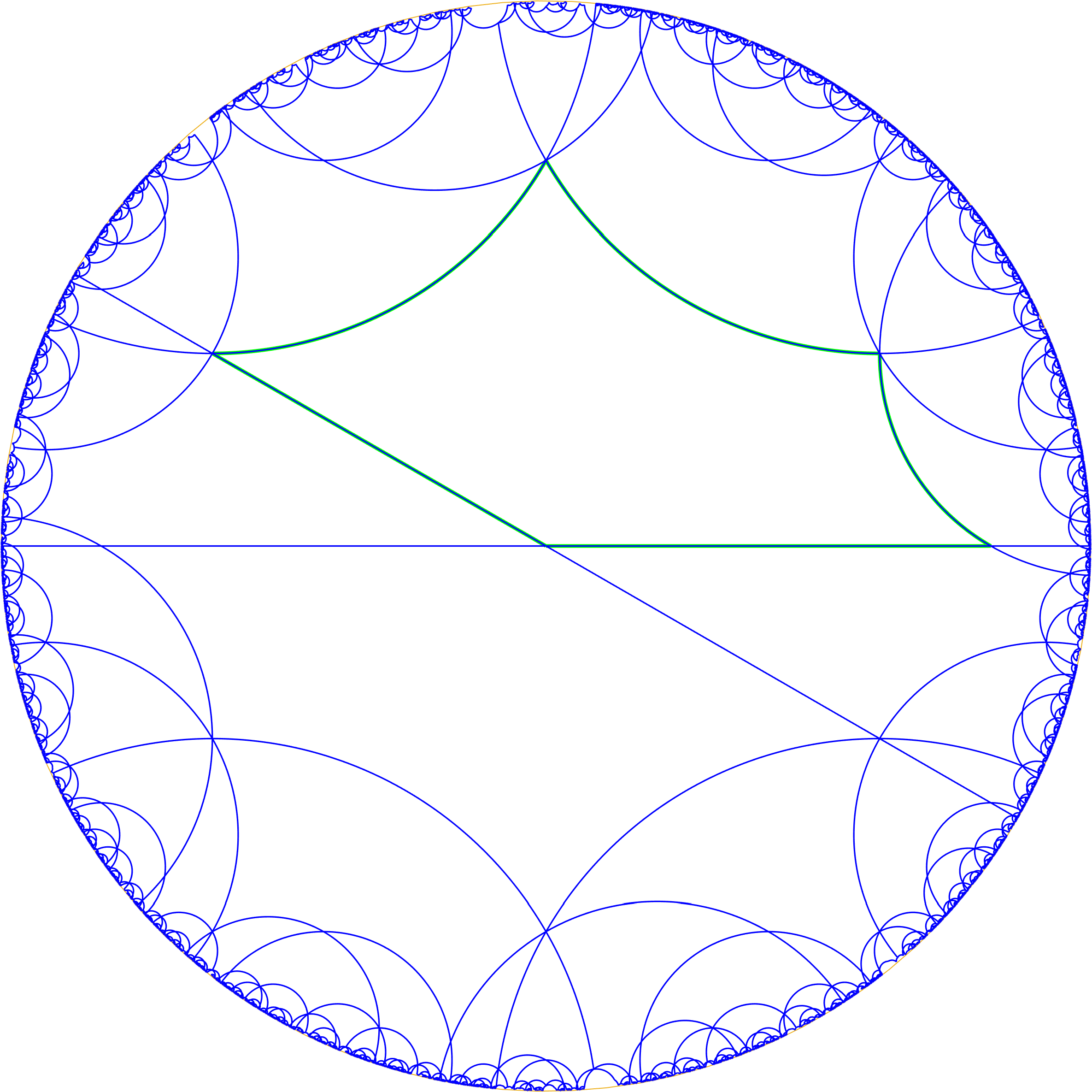}
    \caption{Tiling number $3$.}
  \end{subfigure}
  \\
    \begin{subfigure}[t]{0.31\textwidth}
  \centering
    \includegraphics[width=\textwidth, height=\textwidth]{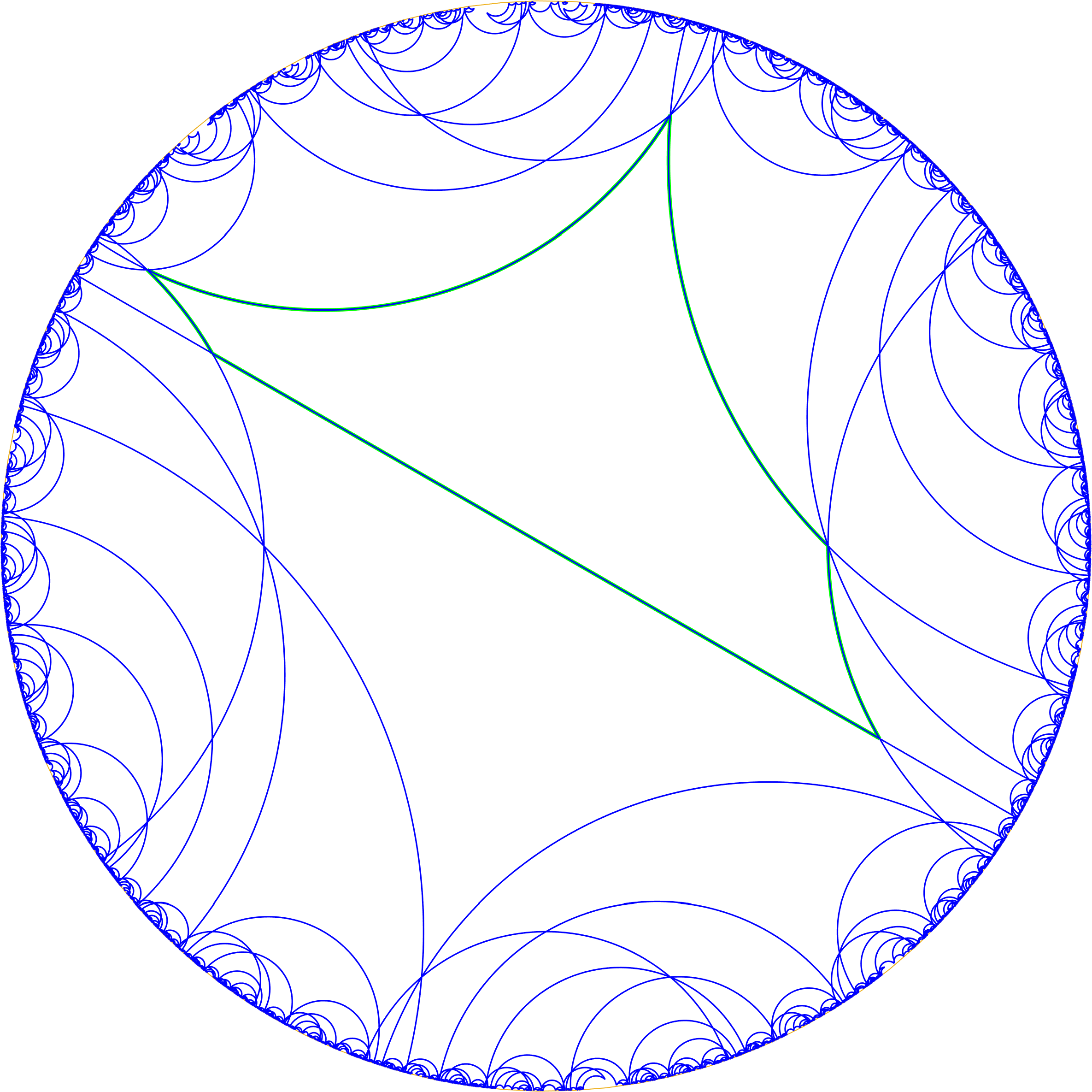}
    \caption{Tiling number $4$.}
  \end{subfigure}
  \begin{subfigure}[t]{0.31\textwidth}
  \centering
    \includegraphics[width=\textwidth, height=\textwidth]{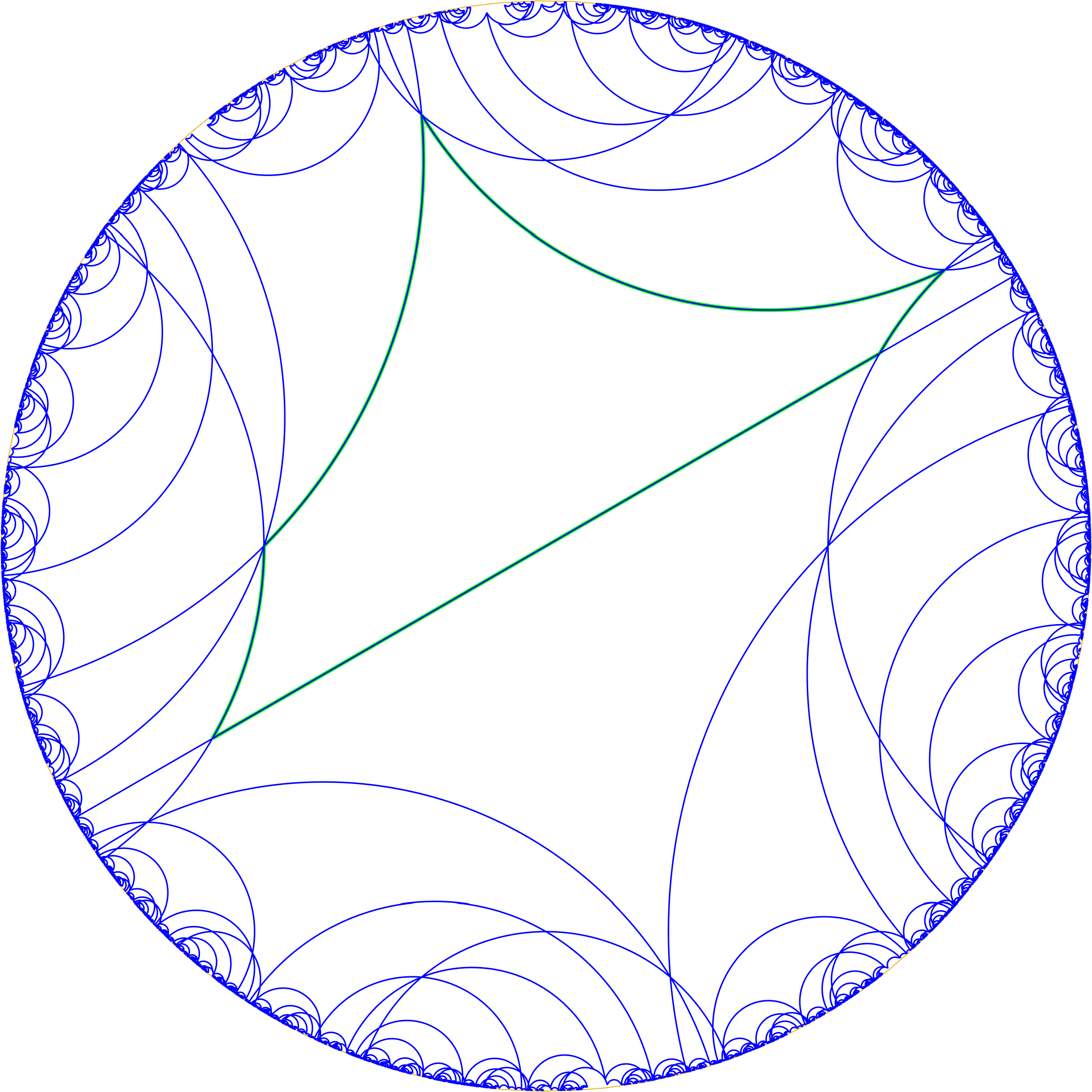}
    \caption{Tiling number $5$.}
  \end{subfigure}
  \begin{subfigure}[t]{0.31\textwidth}
  \centering
    \includegraphics[width=\textwidth, height=\textwidth]{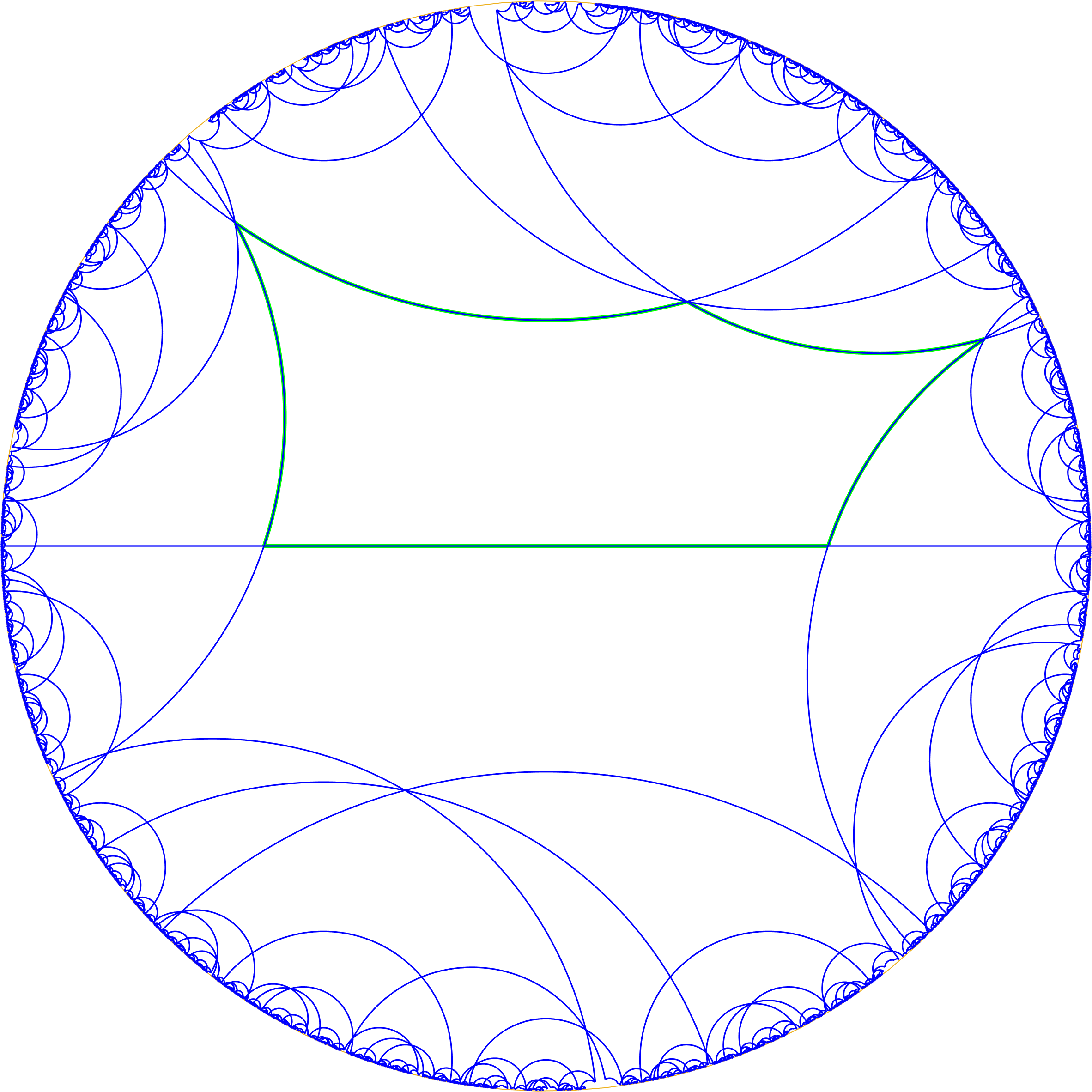}
    \caption{Tiling number $6$.}
  \end{subfigure}
   \\
    \begin{subfigure}[t]{0.31\textwidth}
  \centering
    \includegraphics[width=\textwidth, height=\textwidth]{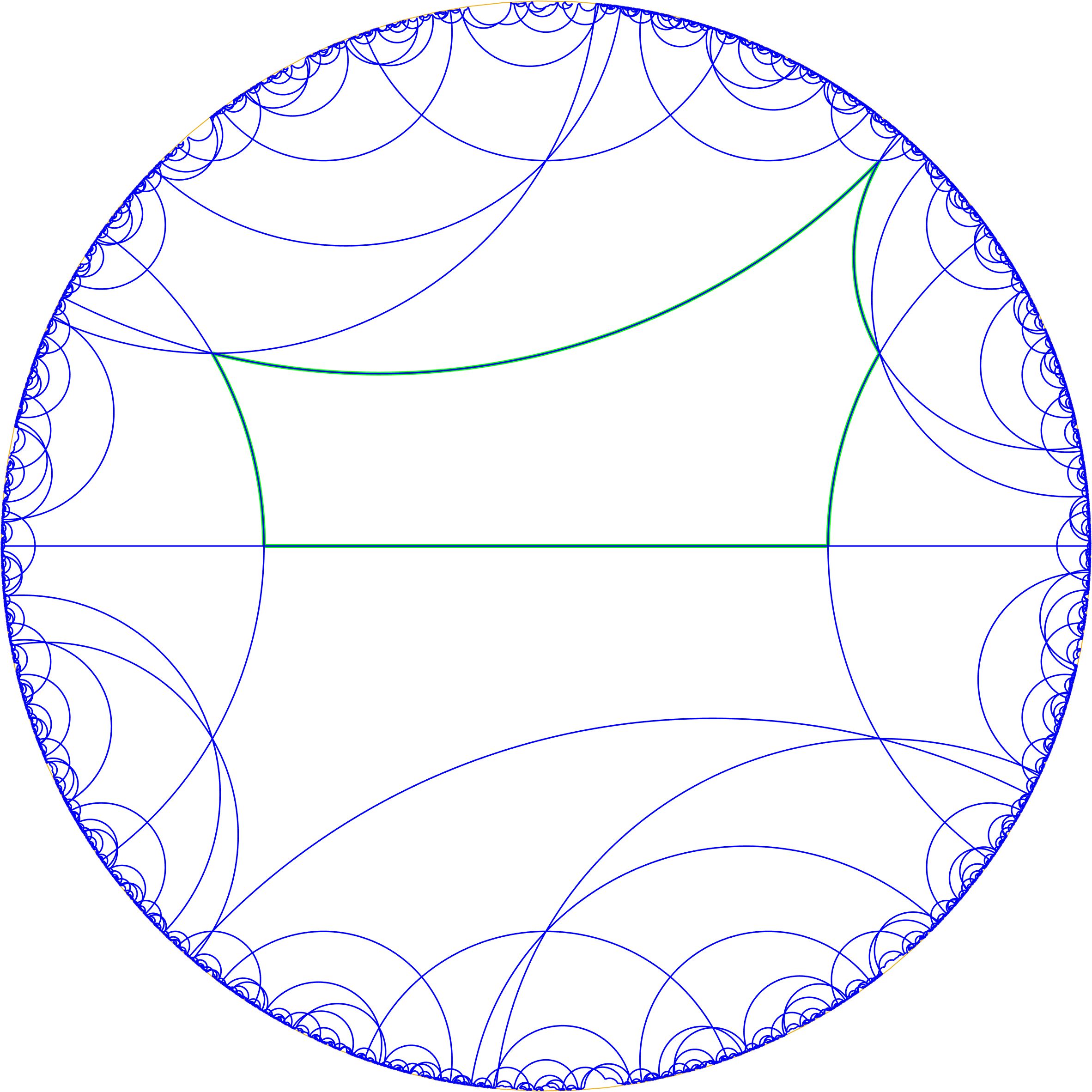}
    \caption{Tiling number $7$.}
  \end{subfigure}
  \begin{subfigure}[t]{0.31\textwidth}
  \centering
    \includegraphics[width=\textwidth, height=\textwidth]{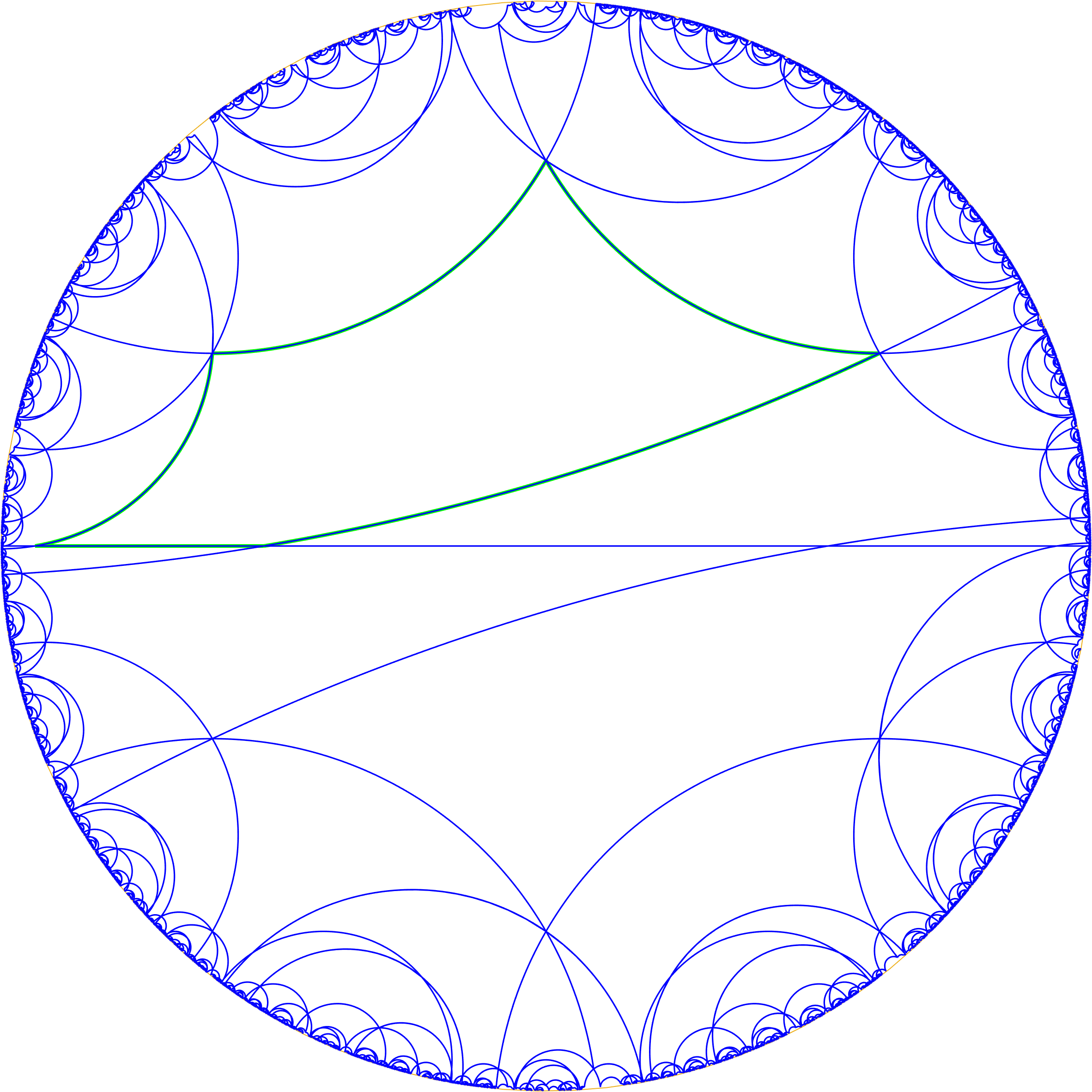}
    \caption{Tiling number $8$.}
  \end{subfigure}
  \begin{subfigure}[t]{0.31\textwidth}
  \centering
    \includegraphics[width=\textwidth, height=\textwidth]{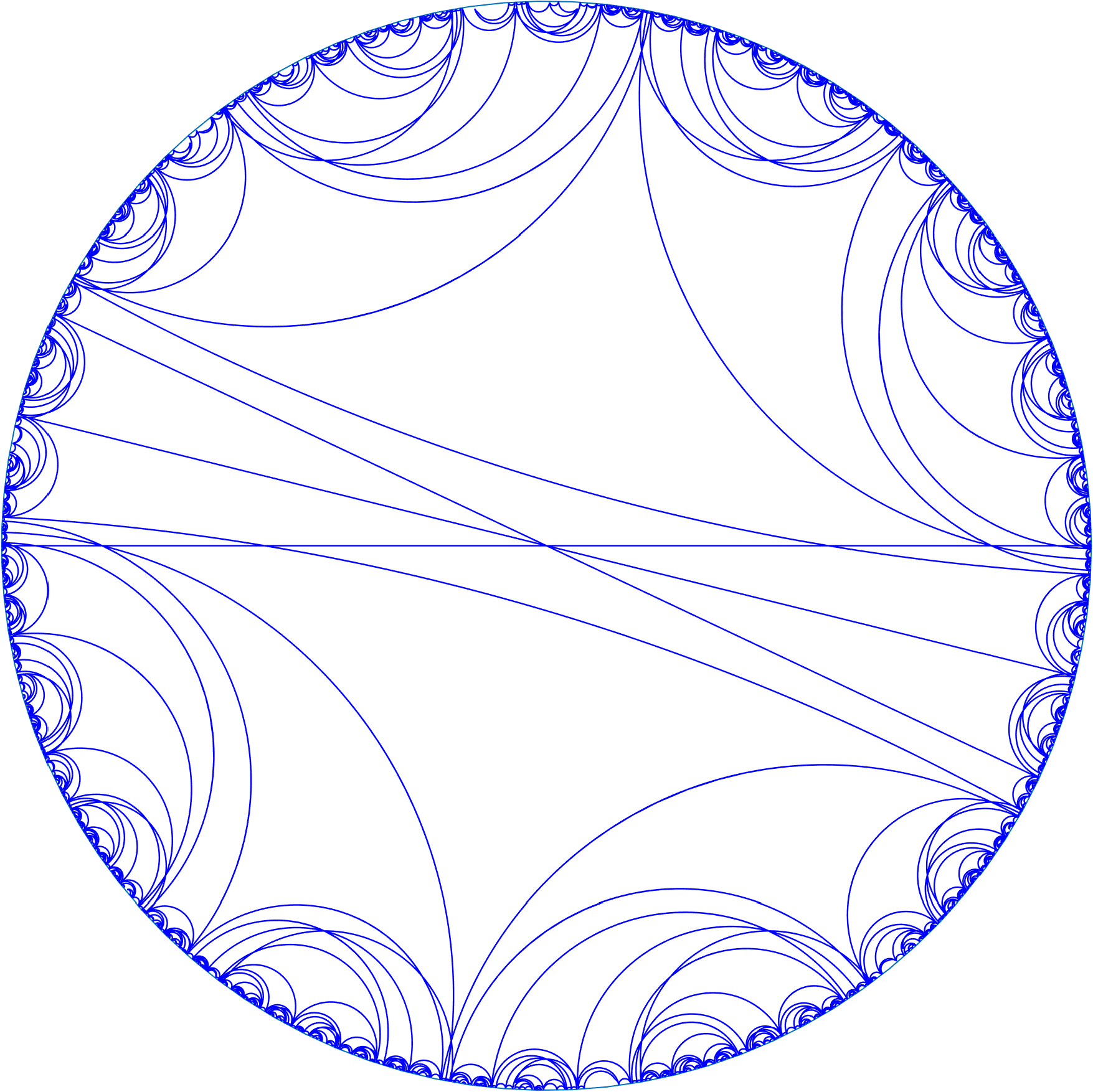}
    \caption{Tiling number $40$.}
  \end{subfigure}
  \caption{Isotopically distinct fundamental tilings of the symmetry group $22222$ (no. $54$ in \cite{Robins2004}) in one realization, with the same D-symbol. (a) shows the placements of the starting generators, with increasingly complicated shearing (nos $1$ though $8$ and $40$ at the bottom right in our enumeration). All tilings are commensurate with the genus-$3$ dodecagon from figure~\ref{fig:tiling3232a}.}\label{fig:tiling2222254}
\end{figure}

\begin{figure}[!htbp]
\imagewidth=0.31\textwidth
\captionsetup[subfigure]{width=0.9\imagewidth,justification=raggedright}
  \begin{subfigure}[t]{0.31\textwidth}
  \centering
    \includegraphics[width=\textwidth, height=\textwidth]{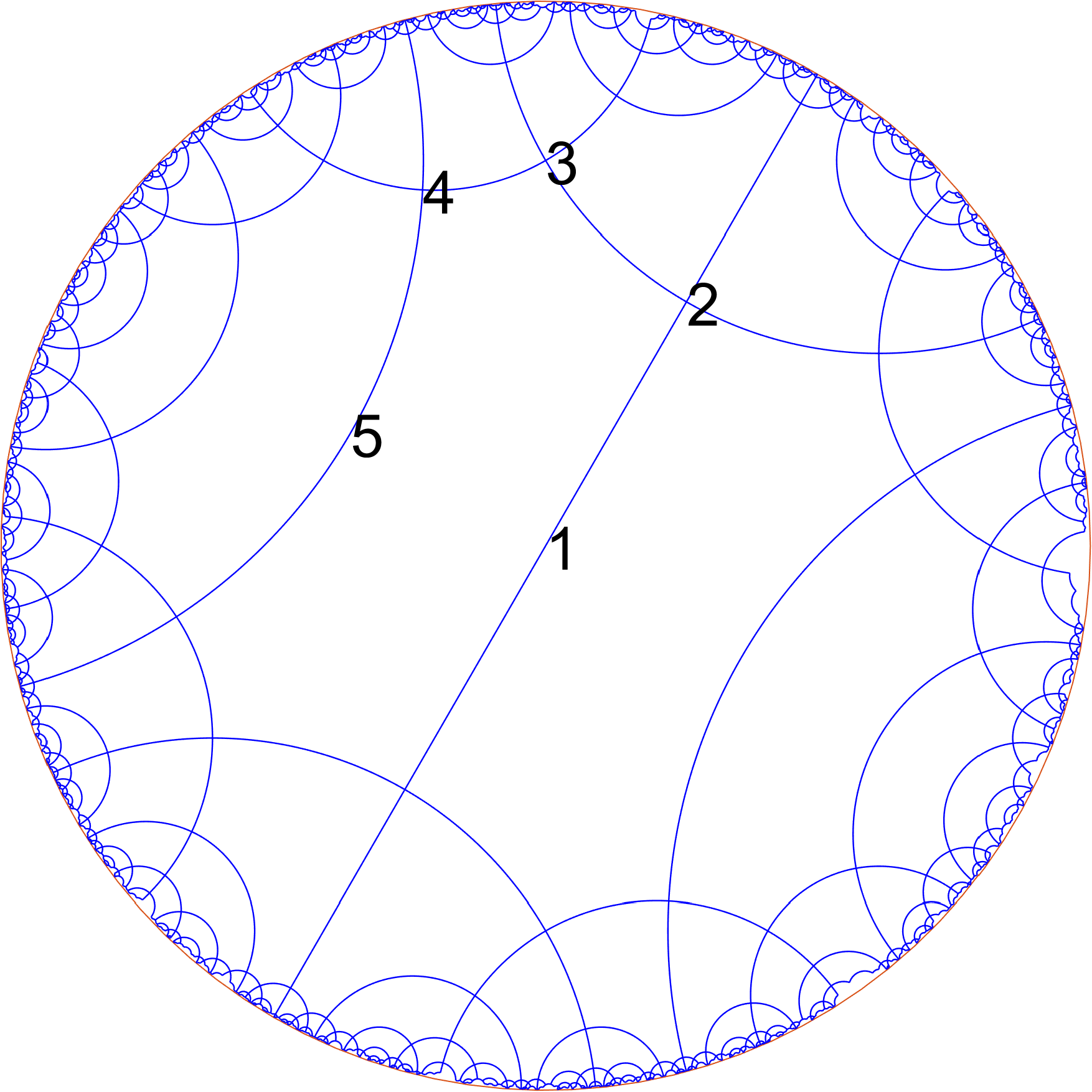}
    \caption{Tiling number $1$.}\label{fig:tiling2222276a}
  \end{subfigure}
  \begin{subfigure}[t]{0.31\textwidth}
  \centering
    \includegraphics[width=\textwidth, height=\textwidth]{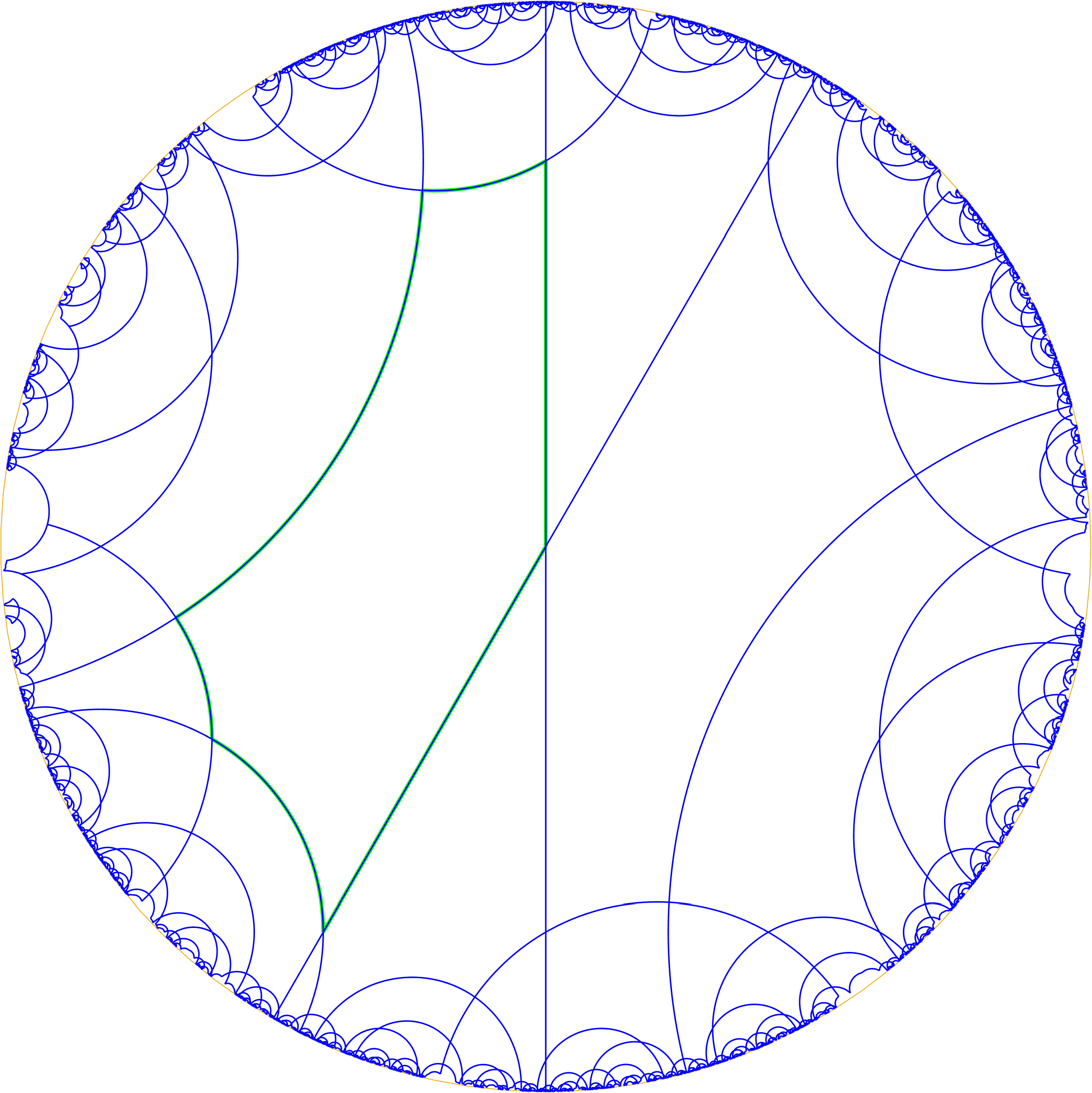}
    \caption{Tiling number $2$.}\label{fig:tiling2222276b}
  \end{subfigure}
  \begin{subfigure}[t]{0.31\textwidth}
  \centering
    \includegraphics[width=\textwidth, height=\textwidth]{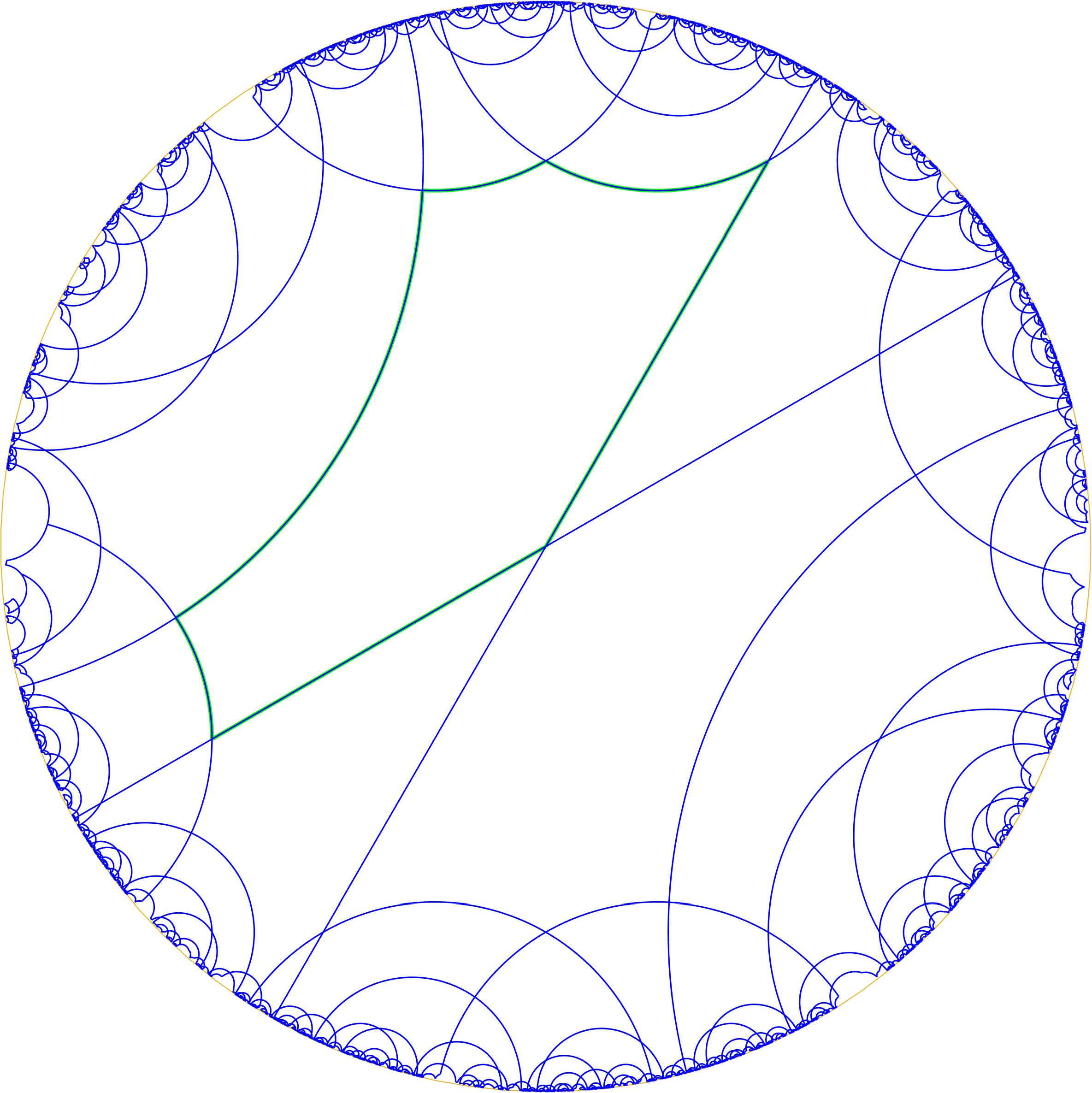}
    \caption{Tiling number $3$.}
  \end{subfigure}
  \\
    \begin{subfigure}[t]{0.31\textwidth}
  \centering
    \includegraphics[width=\textwidth, height=\textwidth]{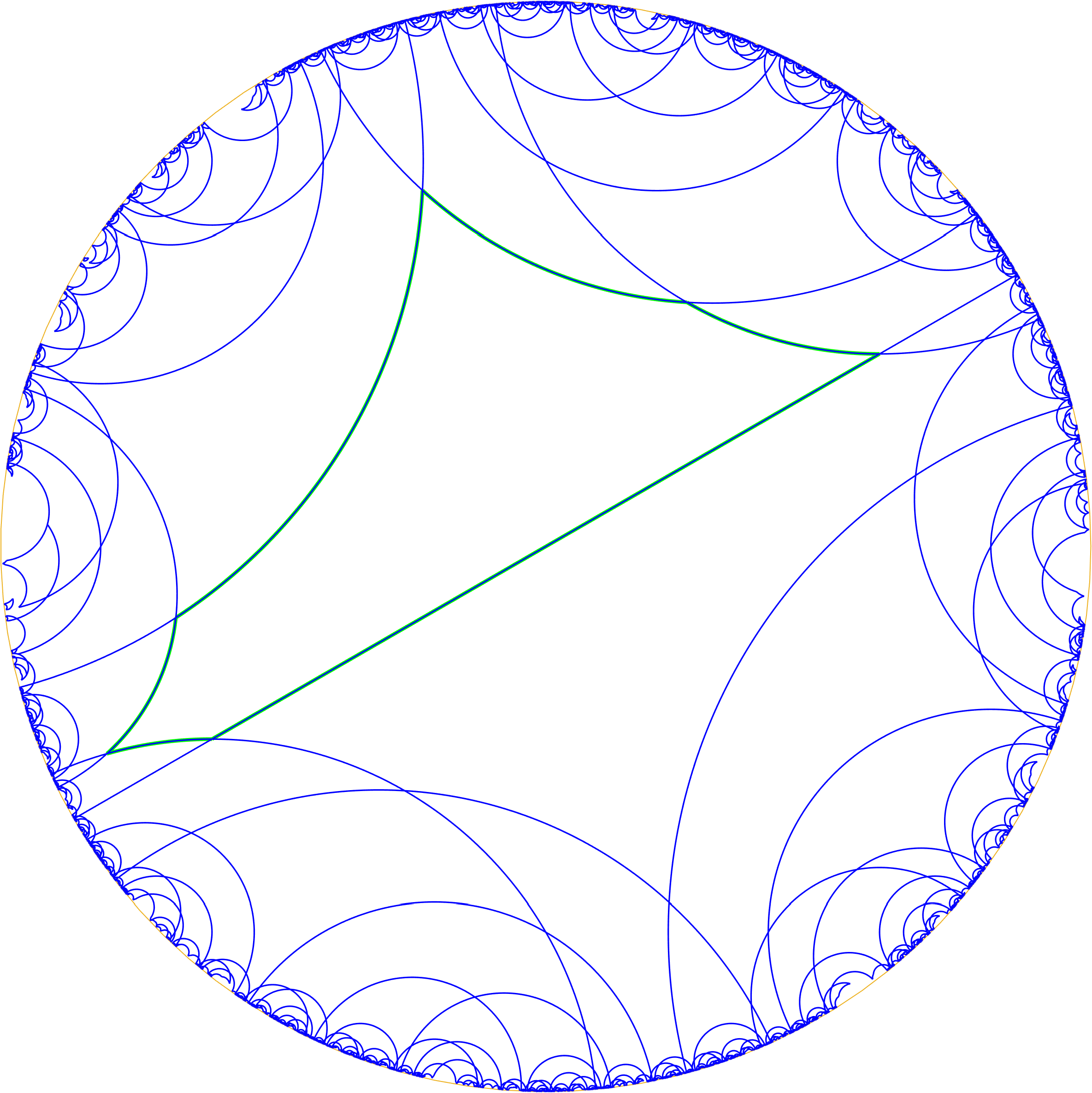}
    \caption{Tiling number $4$.}
  \end{subfigure}
  \begin{subfigure}[t]{0.31\textwidth}
  \centering
    \includegraphics[width=\textwidth, height=\textwidth]{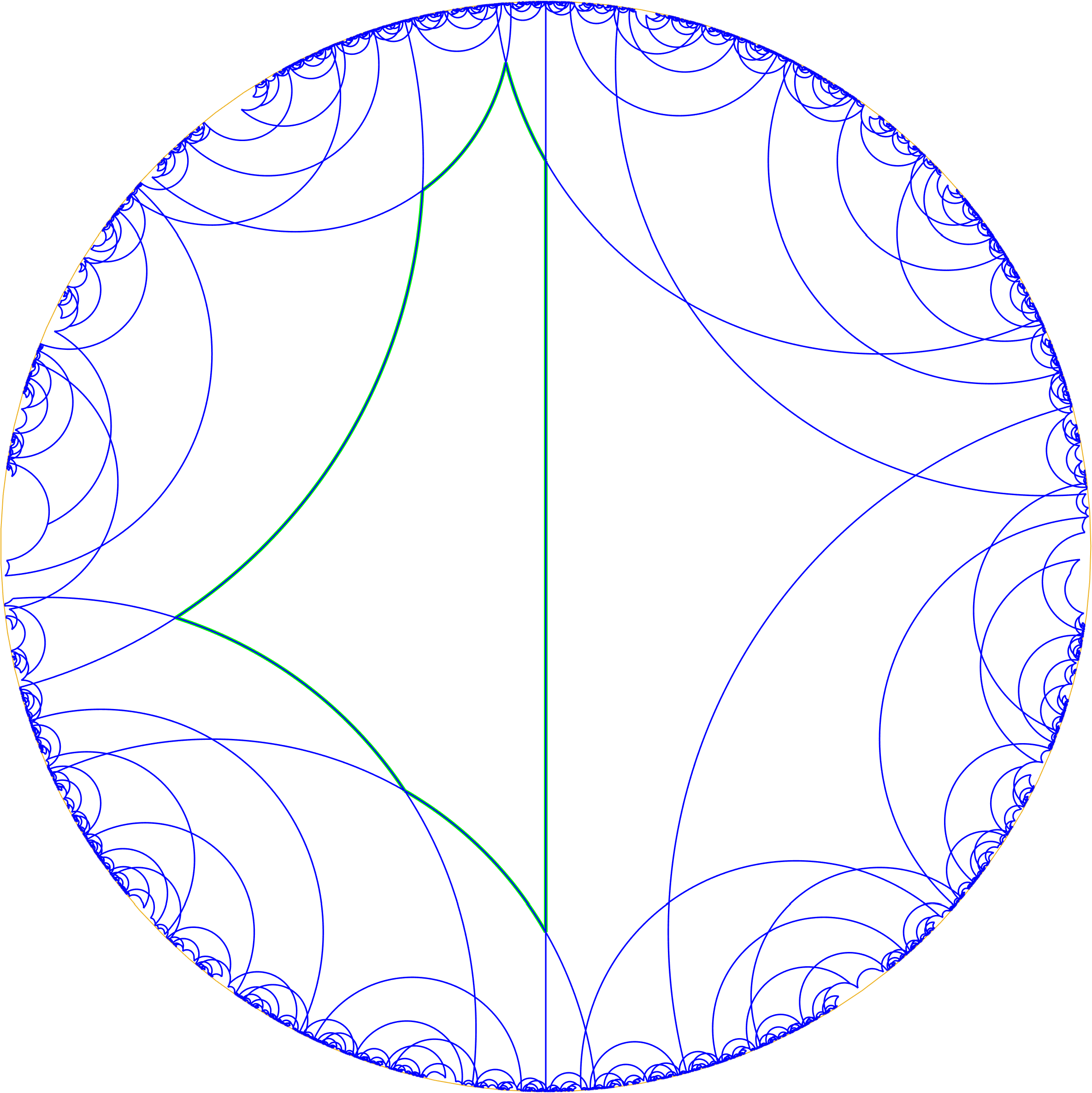}
    \caption{Tiling number $5$.}
  \end{subfigure}
  \begin{subfigure}[t]{0.31\textwidth}
  \centering
    \includegraphics[width=\textwidth, height=\textwidth]{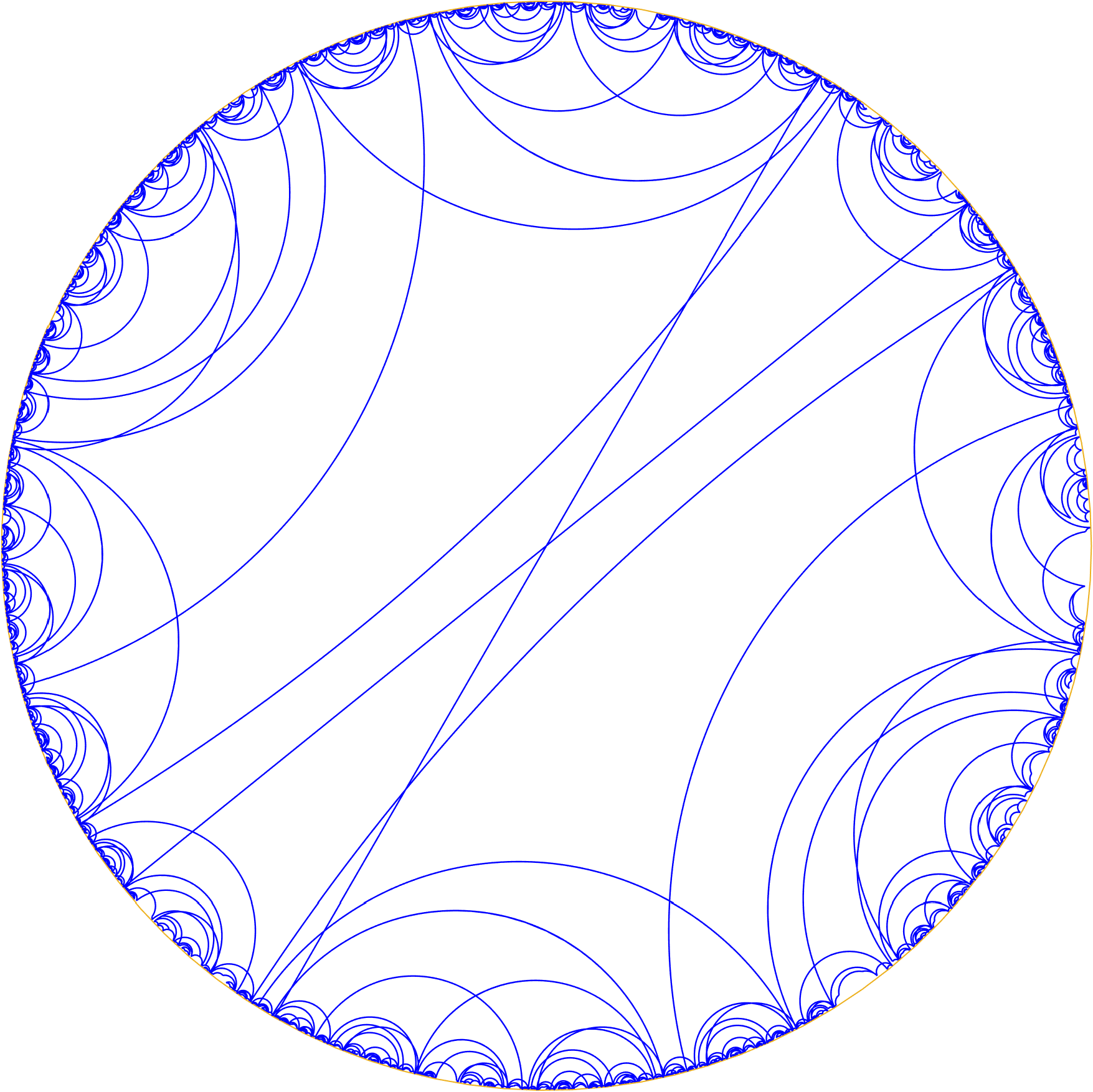}
    \caption{Tiling number $40$.}
  \end{subfigure}
  \caption{Isotopically distinct fundamental tilings of the symmetry group $22222$ (no. $76$ in \cite{Robins2004}) in one realization, with the same D-symbol. (a) shows the placements of the starting generators, with increasingly complicated shearing (nos $1$ though $5$ and $40$ at the bottom right in our enumeration). All tilings are commensurate with the genus-$3$ dodecagon from figure~\ref{fig:tiling3232a}.}\label{fig:tiling2222276}
\end{figure}

\begin{figure}[!htbp]
\imagewidth=0.31\textwidth
\captionsetup[subfigure]{width=0.9\imagewidth,justification=raggedright}
  \begin{subfigure}[t]{0.31\textwidth}
  \centering
    \includegraphics[width=\textwidth, height=\textwidth]{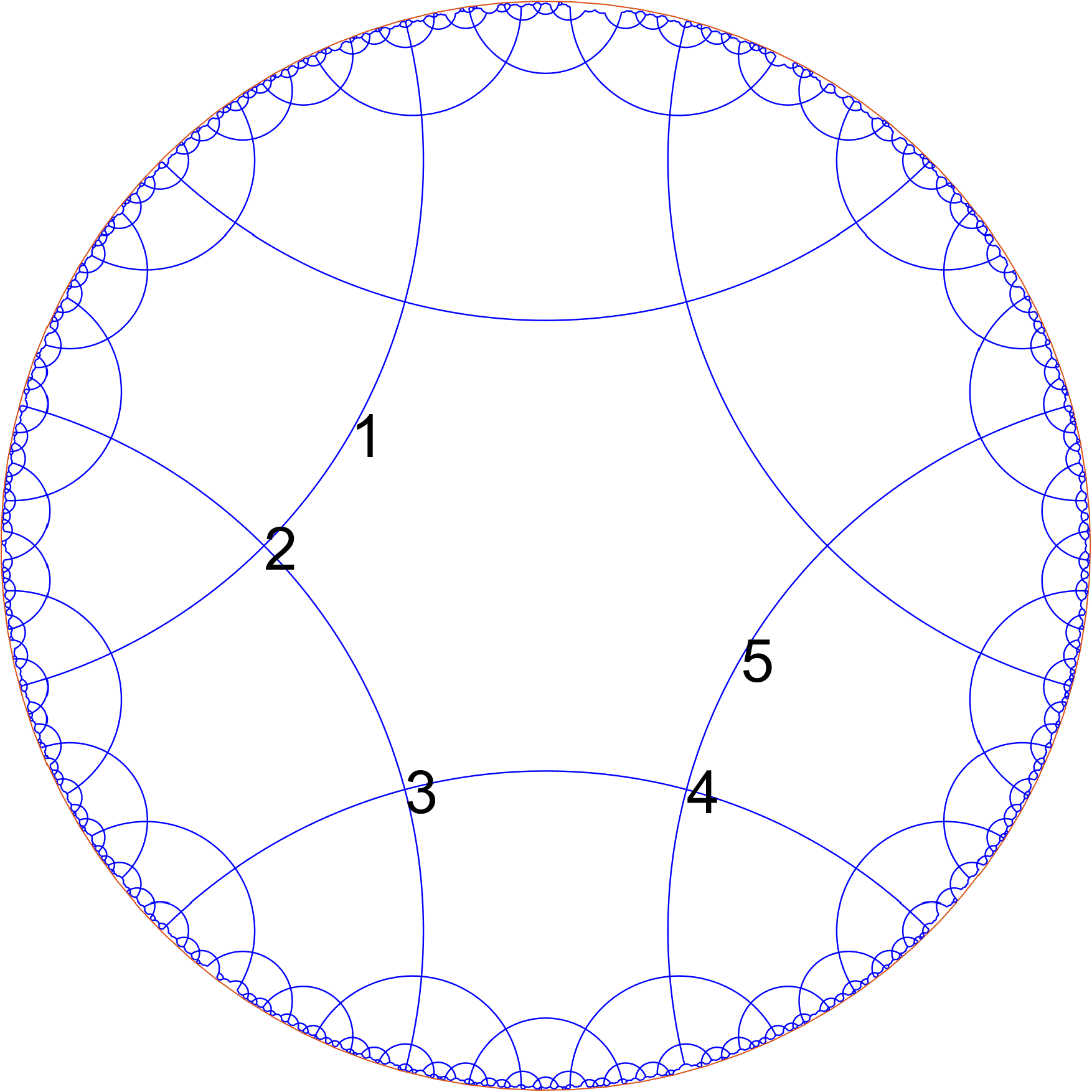}
    \caption{Tiling number $1$.}\label{fig:tiling2222277a}
  \end{subfigure}
  \begin{subfigure}[t]{0.31\textwidth}
  \centering
    \includegraphics[width=\textwidth, height=\textwidth]{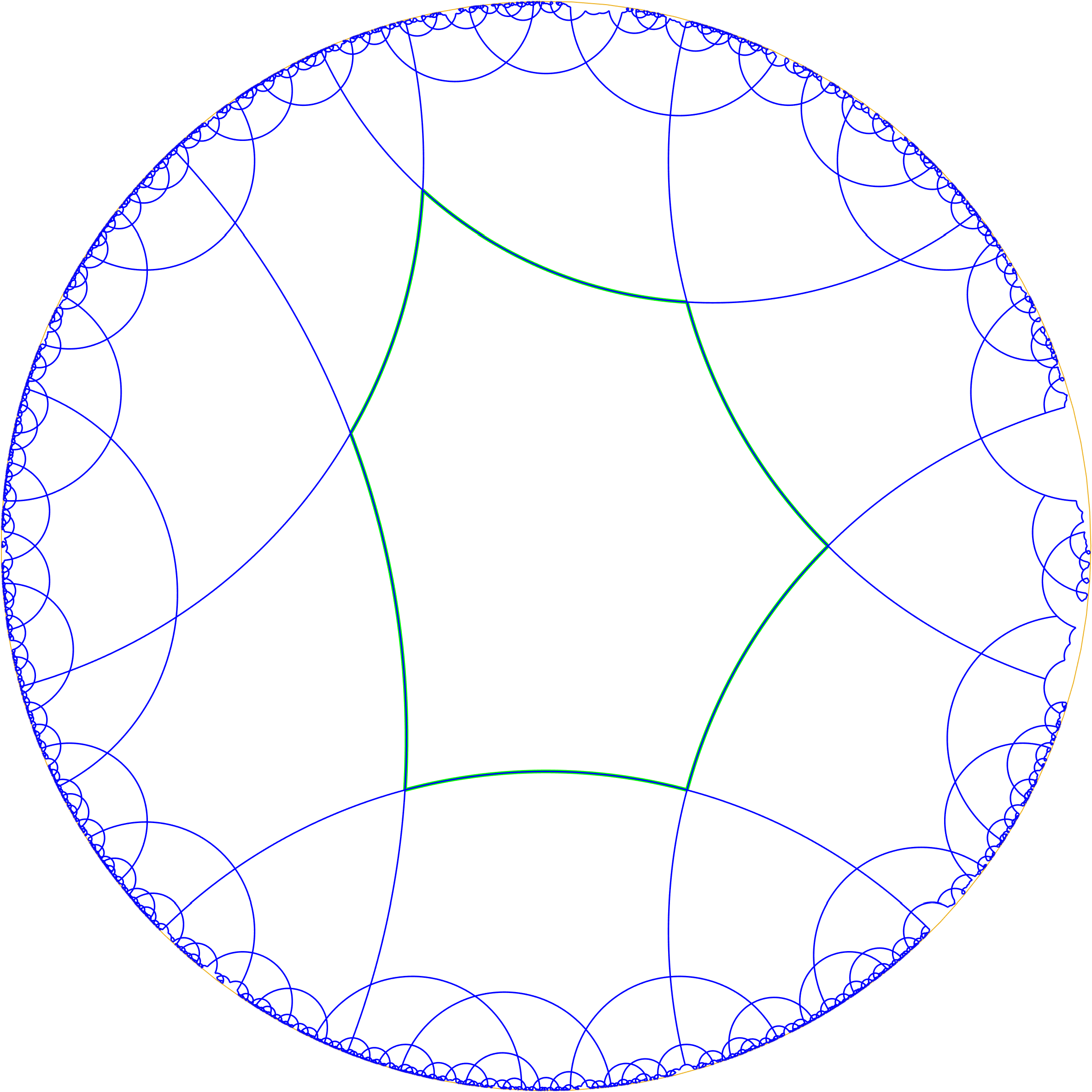}
    \caption{Tiling number $2$.}\label{fig:tiling2222277b}
  \end{subfigure}
  \begin{subfigure}[t]{0.31\textwidth}
  \centering
    \includegraphics[width=\textwidth, height=\textwidth]{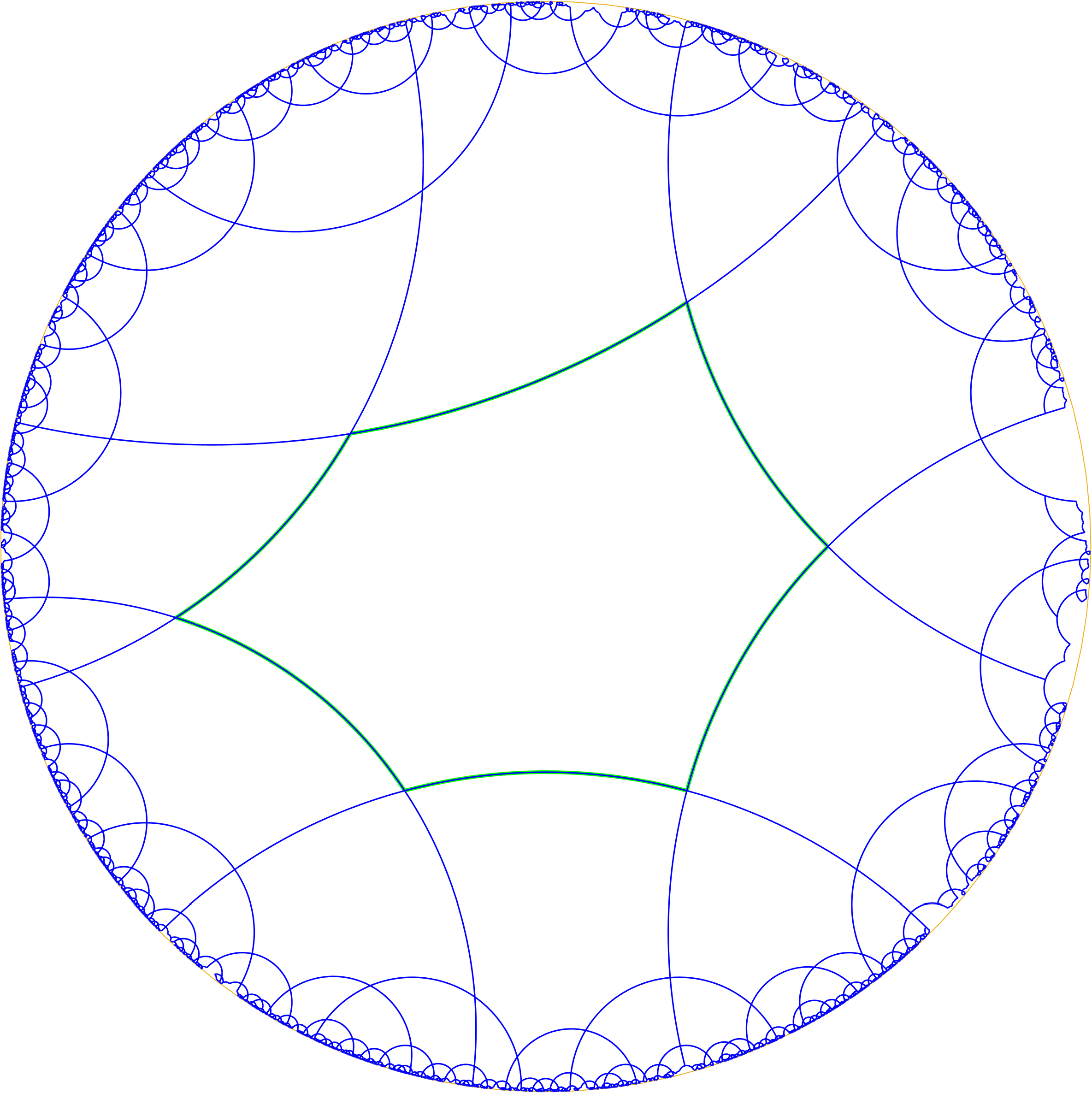}
    \caption{Tiling number $3$.}
  \end{subfigure}
  \\
    \begin{subfigure}[t]{0.31\textwidth}
  \centering
    \includegraphics[width=\textwidth, height=\textwidth]{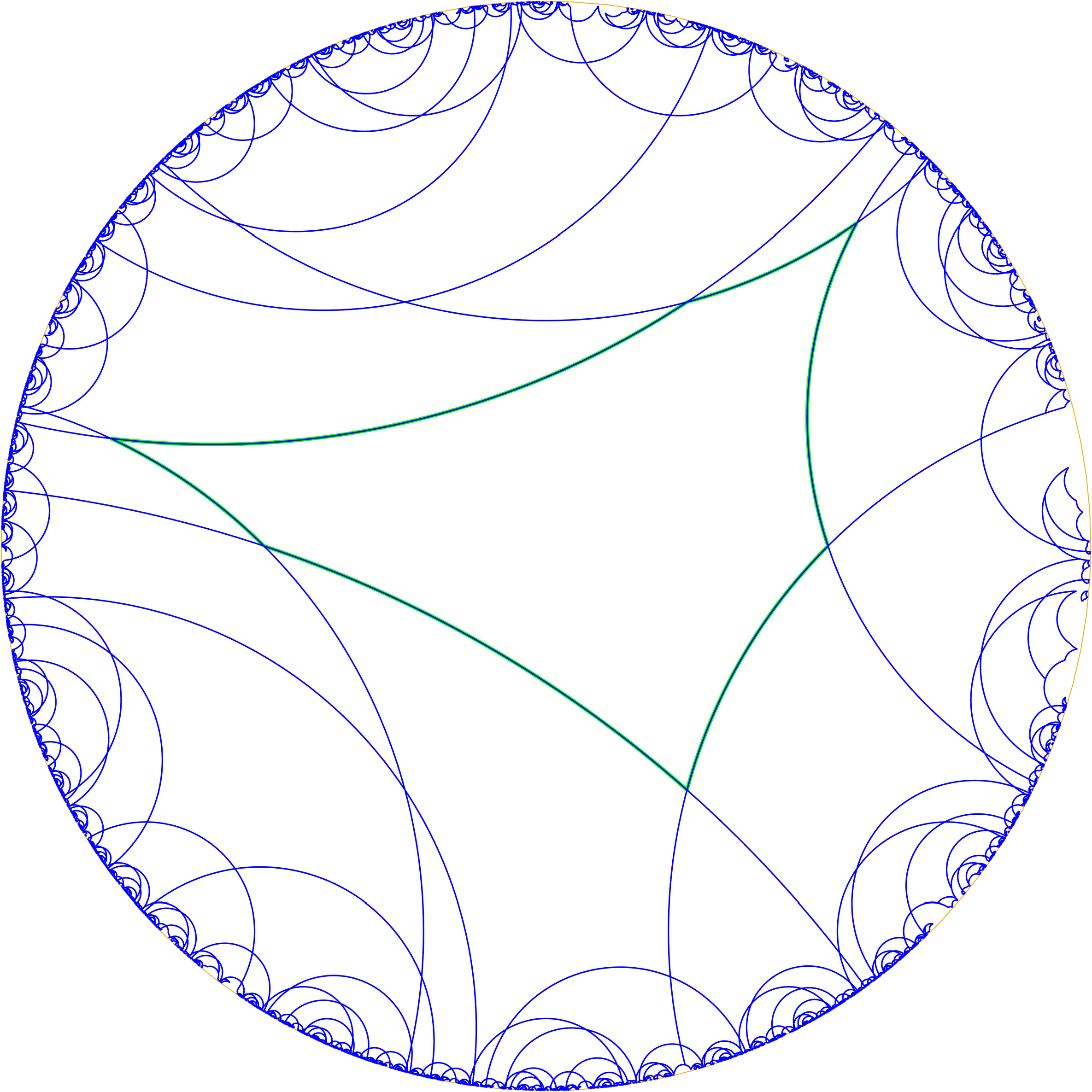}
    \caption{Tiling number $4$.}
  \end{subfigure}
  \begin{subfigure}[t]{0.31\textwidth}
  \centering
    \includegraphics[width=\textwidth, height=\textwidth]{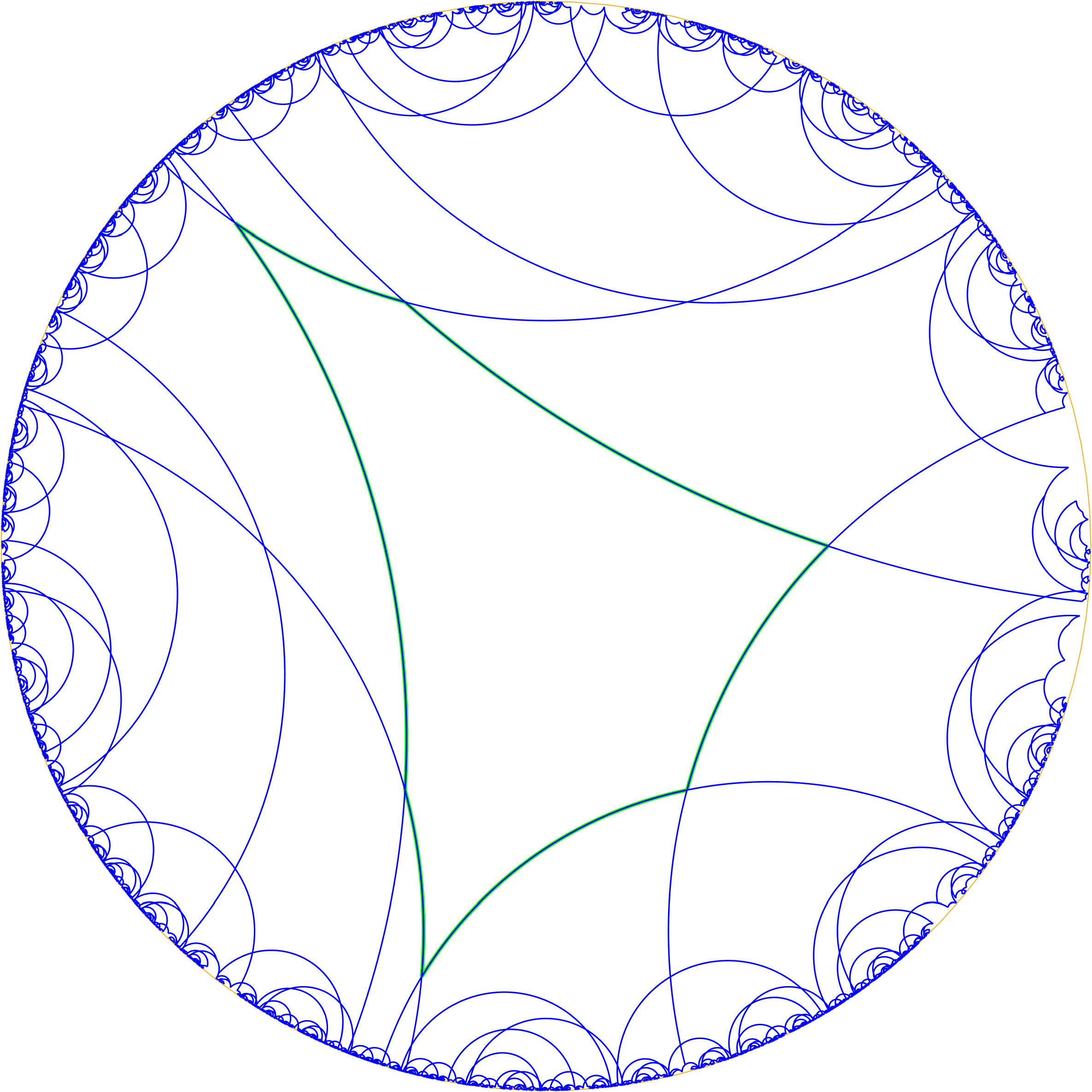}
    \caption{Tiling number $5$.}
  \end{subfigure}
  \begin{subfigure}[t]{0.31\textwidth}
  \centering
    \includegraphics[width=\textwidth, height=\textwidth]{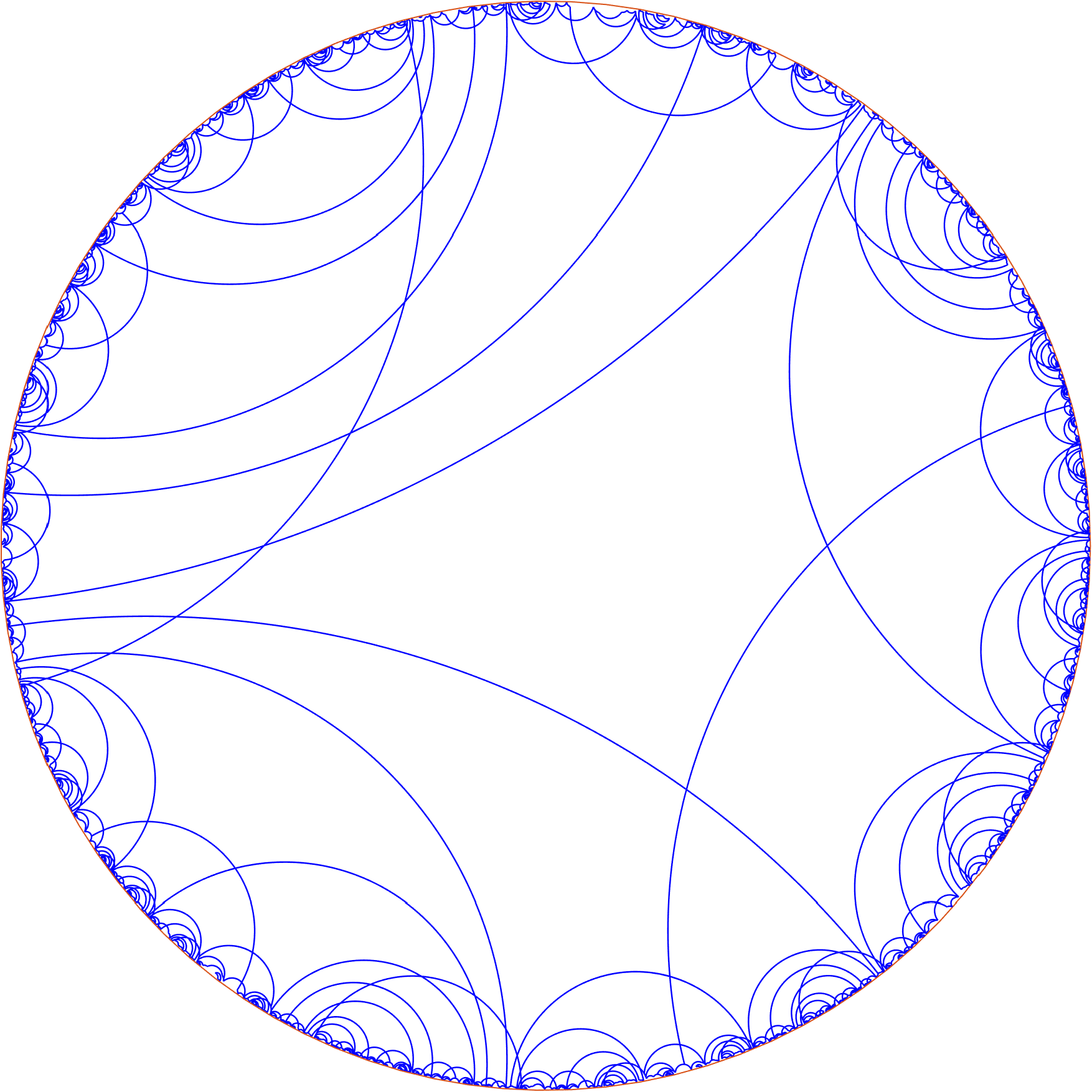}
    \caption{Tiling number $40$.}
  \end{subfigure}
  \caption{Isotopically distinct fundamental tilings of the symmetry group $22222$ (no. $77$ in \cite{Robins2004}) in one realization, with the same D-symbol. (a) shows the placements of the starting generators, with increasingly complicated shearing (nos $1$ though $5$ and $40$ at the bottom right in our enumeration). All tilings are commensurate with the genus-$3$ dodecagon from figure~\ref{fig:tiling3232a}.}\label{fig:tiling2222277}
\end{figure}

Figure~\ref{fig:tiling3232} shows isotopically distinct fundamental tilings with the same realization of the symmetry group $2233$ and the same D-symbol. Of the MCGs in table~\ref{table:stellatemcgs}, $2233$ and $2266$ use most of the theoretical results discussed in the previous sections. On the other hand, groups like $22222$ turn out to be particularly simple, using only \eqref{eq:stelmcgpresentation}, \eqref{eq:half-twists} and the arrangements of the group elements in a fundamental tiling. Figure~\ref{fig:tiling2224} shows isotopically distinct fundamental tilings with symmetry group $2224$ and the same combinatorial structure, continuing the enumeration started in figure~\ref{fig:tiling2224ht}. Figures~\ref{fig:tiling2222276},~\ref{fig:tiling2222277} and \ref{fig:tiling2222254} show isotopically distinct tilings invariant under isomorphic symmetry groups (in this case $22222$). Figures~\ref{fig:tiling2222276} and \ref{fig:tiling2222277} show isotopically distinct sets of tilings with the same D-symbol, whereas figure~\ref{fig:tiling2222254} shows a further set, with a different D-symbol, of mutually combinatorially equivalent tilings. We concentrated here on tilings with vertices only at points of increased symmetry, to better illustrate the action of the MCG on the tilings. Figure \ref{fig:tiling2222277} shows how the degree of distortion of the initial hexagonal fundamental domain increases as one steps through the enumeration. Notice that the starting placement of every set of geometric generators is induced by doubling a fundamental domain for the Coxeter supergroup of the stellate group, except for that of Figure~\ref{fig:tiling2222254a}. The fundamental domain and the starting set there yield a minimally sheared domain as a starting point for the enumeration, but the combinatorial class of the tiling w.r.t. which we found a suitable starting set of generators is different than that of Figures~\ref{fig:tiling2222276a} and \ref{fig:tiling2222277a}, illustrating that while we can narrow down the choice of starting set to a finite number, there is no unique canonical choice. On the other hand, note that all three starting fundamental domains with their placements for generators from Figures~\ref{fig:tiling2222254a}, \ref{fig:tiling2222276a} and \ref{fig:tiling2222277a} are a result of gluing $\star 246$ triangles as in Figure~\ref{fig:symmetries} together.

In the enumeration, each tiling has the same number of vertices with the same valency in the fundamental dodecagon illustrated in Figure~\ref{fig:tiling3232a} and therefore on genus $3$ hyperbolic surface $S$ that gives rise to the TPMS. On the other hand, there are many different tilings of $S$ that are not a result of applying a homeomorphism of $S$ the to the graph of the tile edges. In particular, some of the graphs of tile edges projected to $S$ will differ by entanglement arising from their embedding into $\mathbb{R}^3$, while others will be different as graphs. We leave the projection to the TPMS and analysis of the resulting structures for future endeavours.


\begin{figure}[!htbp]
\imagewidth=0.31\textwidth
\captionsetup[subfigure]{width=0.9\imagewidth,justification=raggedright}
  \begin{subfigure}[t]{0.31\textwidth}
  \centering
    \includegraphics[width=\textwidth, height=\textwidth]{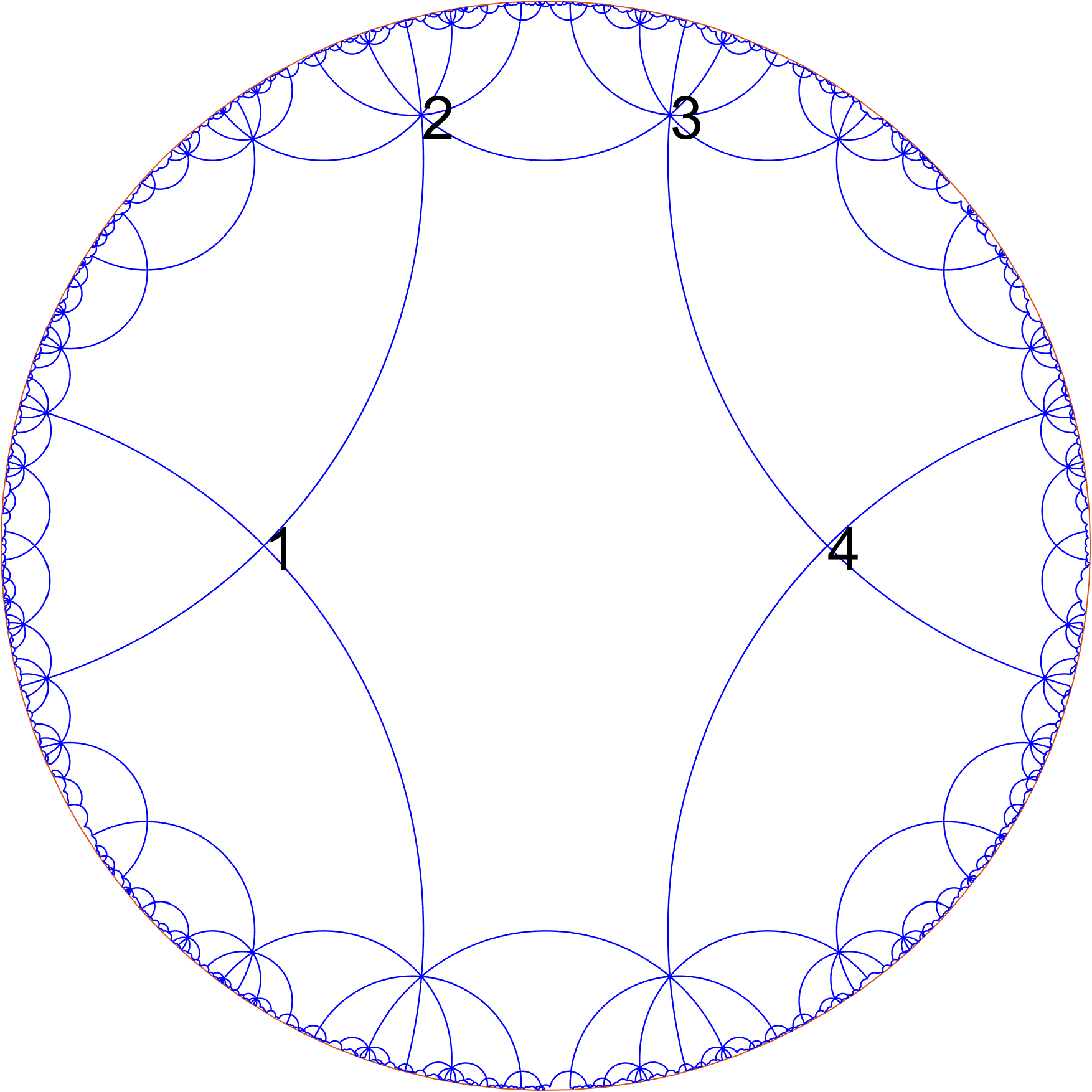}
    \caption{Tiling number $1$.}\label{fig:tiling4444a}
  \end{subfigure}
  \begin{subfigure}[t]{0.31\textwidth}
  \centering
    \includegraphics[width=\textwidth, height=\textwidth]{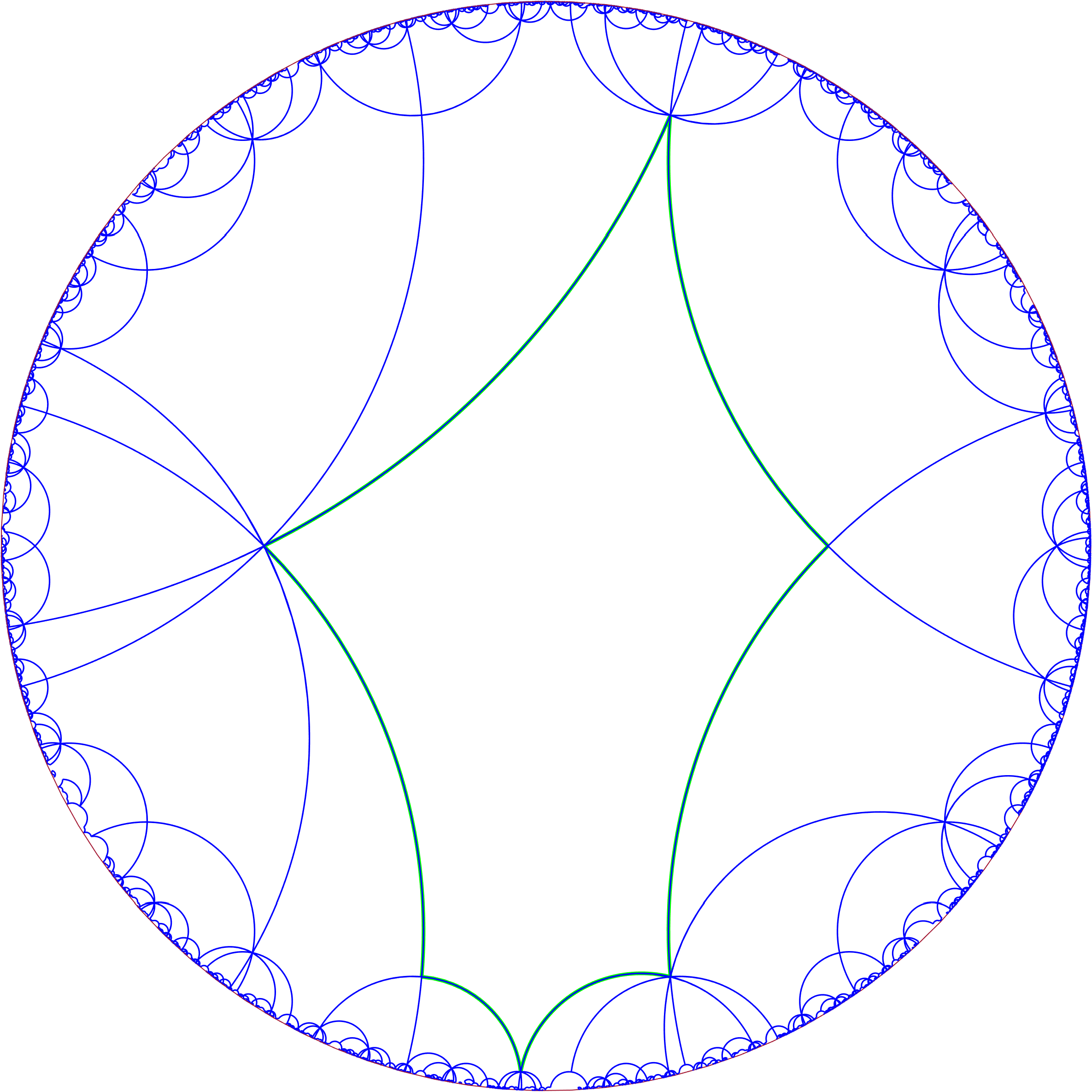}
    \caption{Tiling number $2$.}
  \end{subfigure}
  \begin{subfigure}[t]{0.31\textwidth}
  \centering
    \includegraphics[width=\textwidth, height=\textwidth]{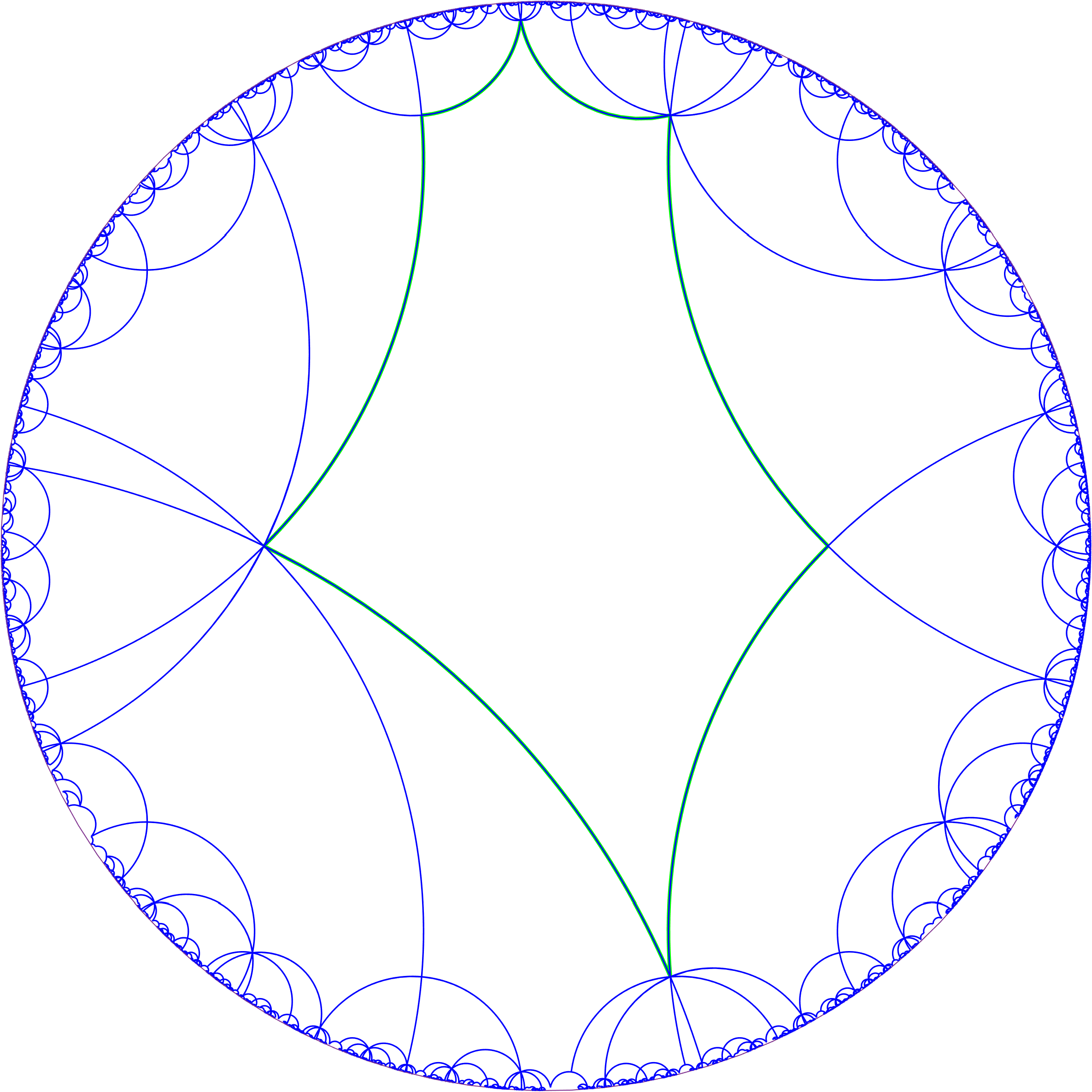}
    \caption{Tiling number $3$.}
  \end{subfigure}
  \\
    \begin{subfigure}[t]{0.31\textwidth}
  \centering
    \includegraphics[width=\textwidth, height=\textwidth]{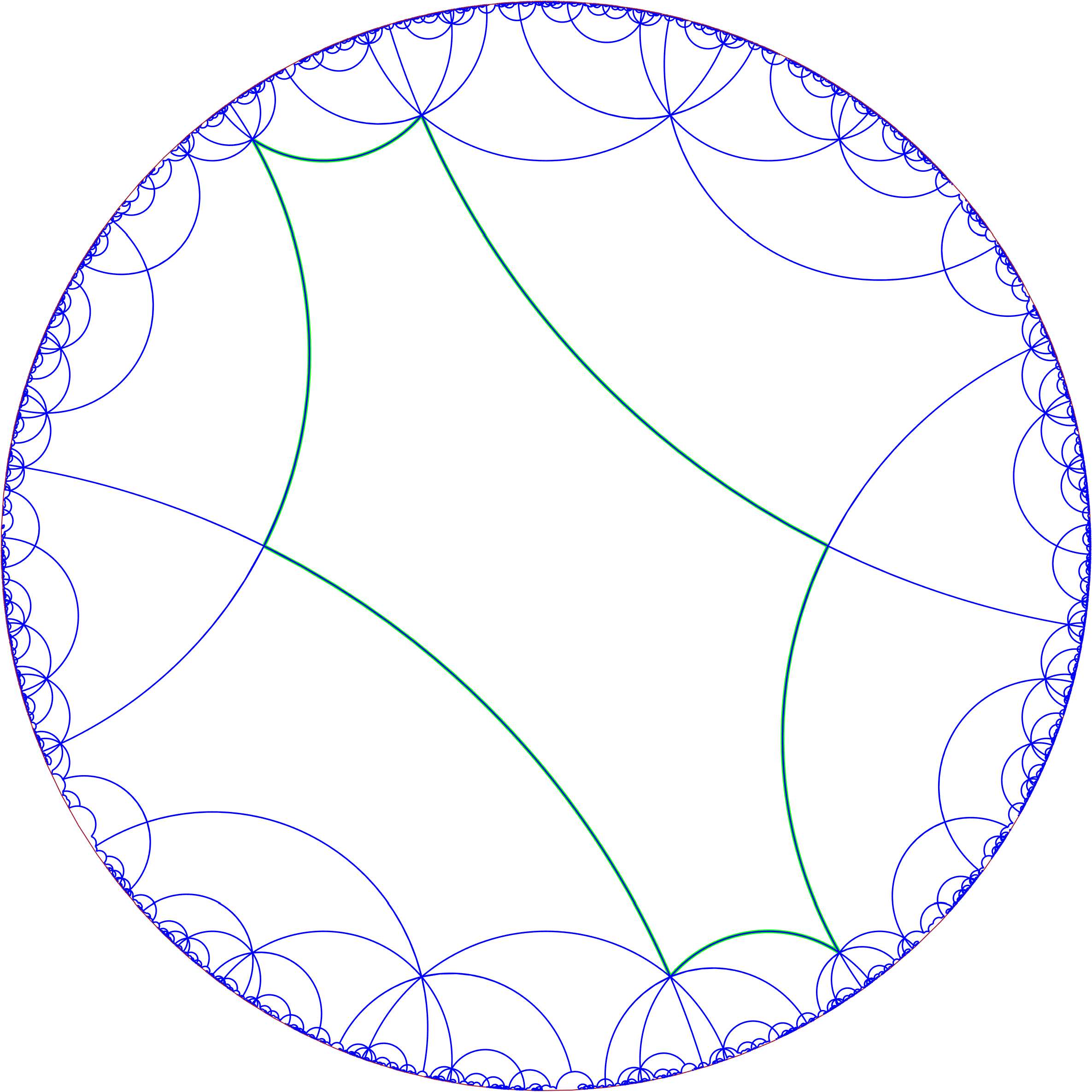}
    \caption{Tiling number $4$.}
  \end{subfigure}
  \begin{subfigure}[t]{0.31\textwidth}
  \centering
    \includegraphics[width=\textwidth, height=\textwidth]{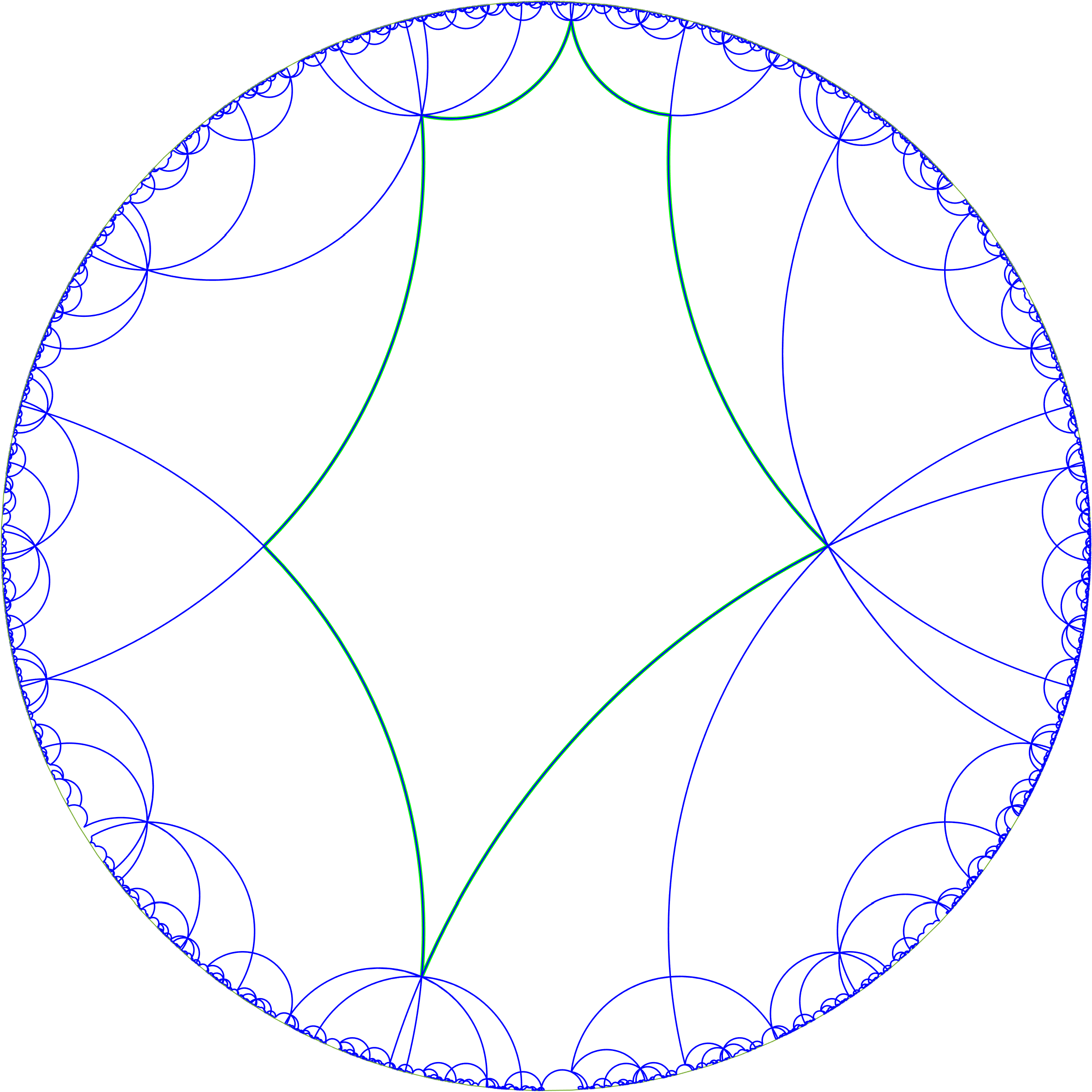}
    \caption{Tiling number $5$.}
  \end{subfigure}
  \begin{subfigure}[t]{0.31\textwidth}
  \centering
    \includegraphics[width=\textwidth, height=\textwidth]{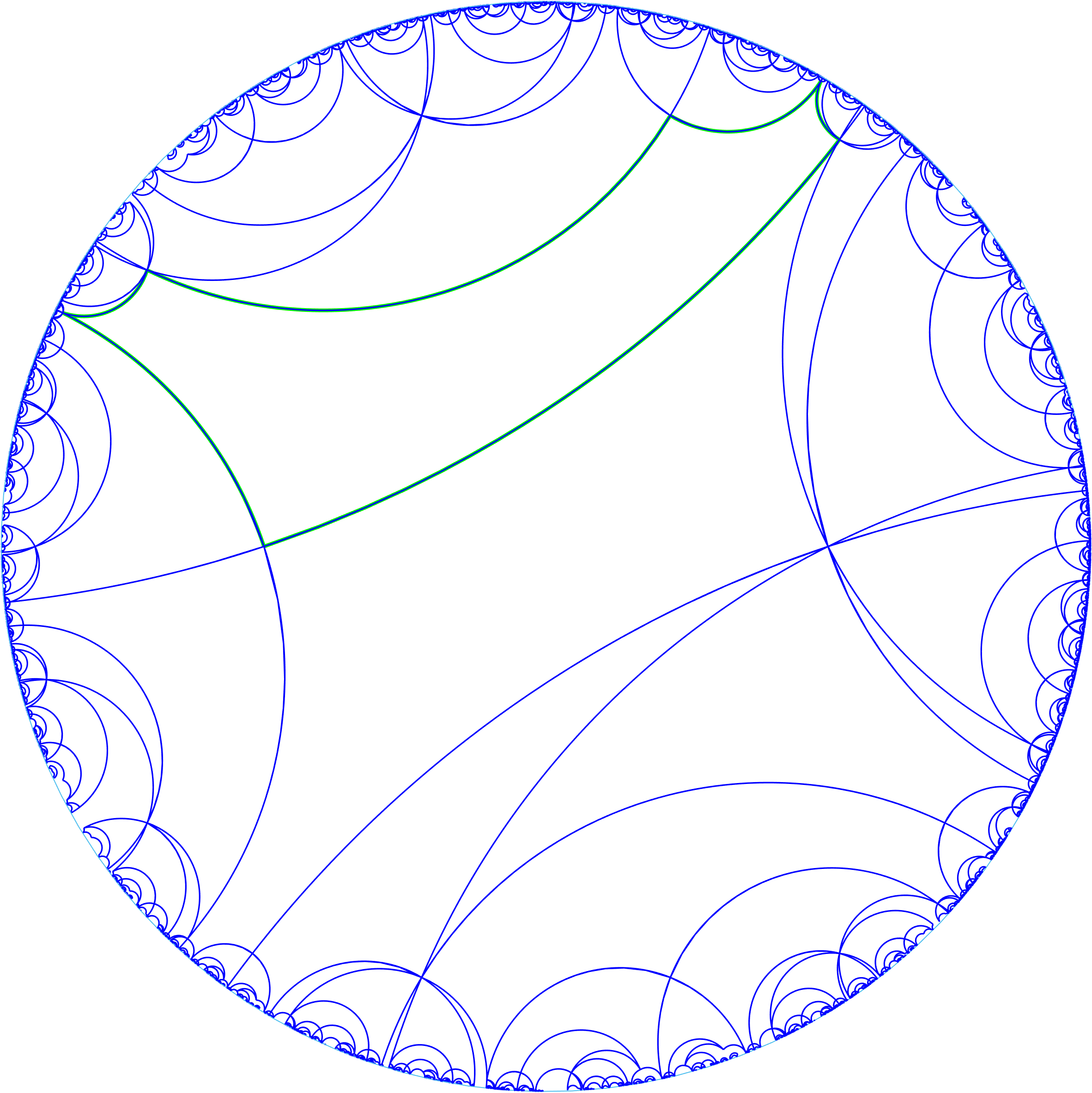}
    \caption{Tiling number $41$.}
  \end{subfigure}
  \caption{Isotopically distinct fundamental tilings of the symmetry group $4444$ with the same D-symbol. (a) shows the placements of the starting generators, with increasingly complicated (nos $1$ though $5$ and $41$ at the bottom right in our enumeration). All tilings are commensurate with the genus-$3$ dodecagon from figure~\ref{fig:tiling3232a}.}\label{fig:tiling4444}
\end{figure}

\begin{figure}[!htbp]
  \begin{subfigure}[t]{0.3\textwidth}
  \centering
    \includegraphics[width=\textwidth]{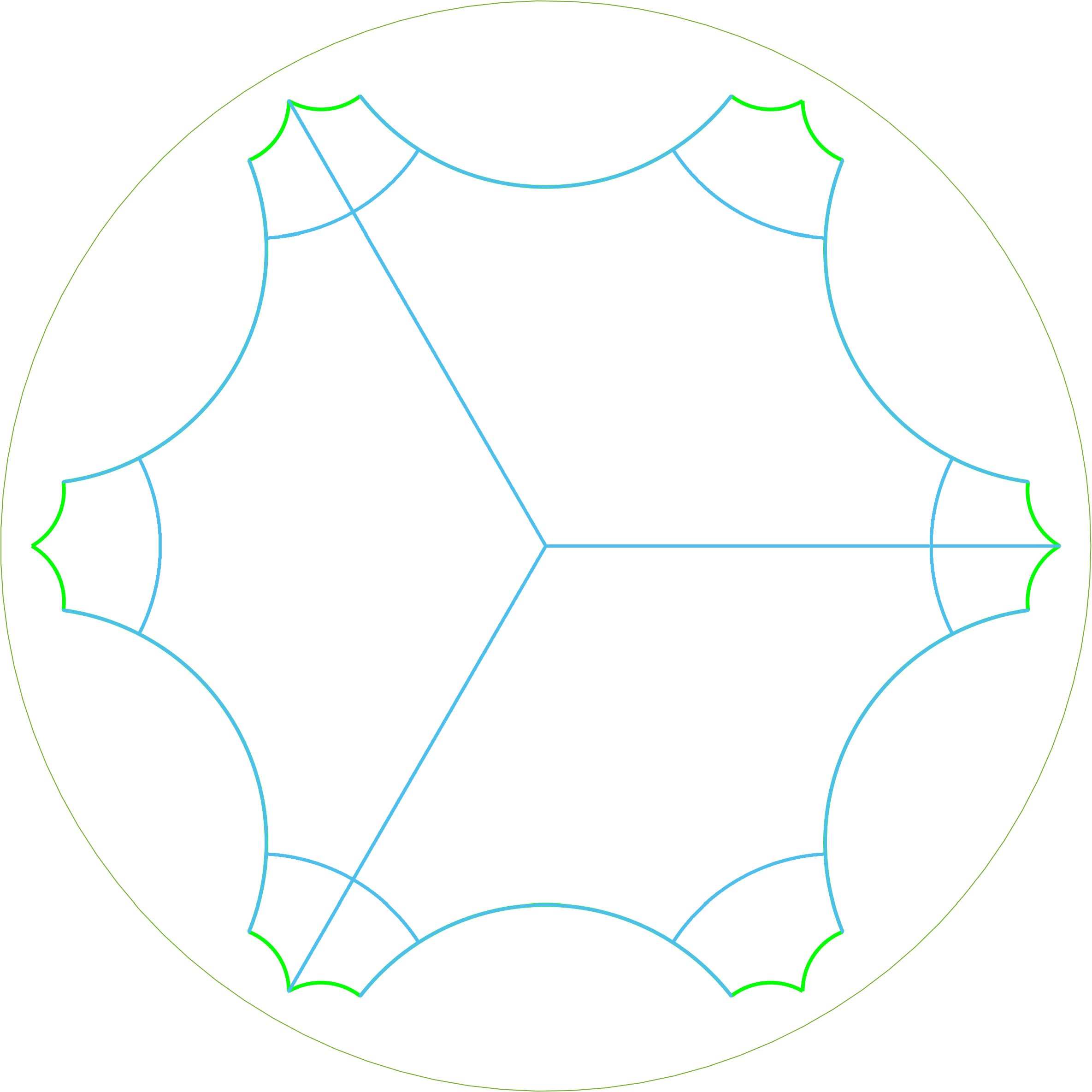}
    \caption{Tiling number $1$.}
  \end{subfigure}
  \hfill
  \begin{subfigure}[t]{0.3\textwidth}
  \centering
    \includegraphics[width=\textwidth]{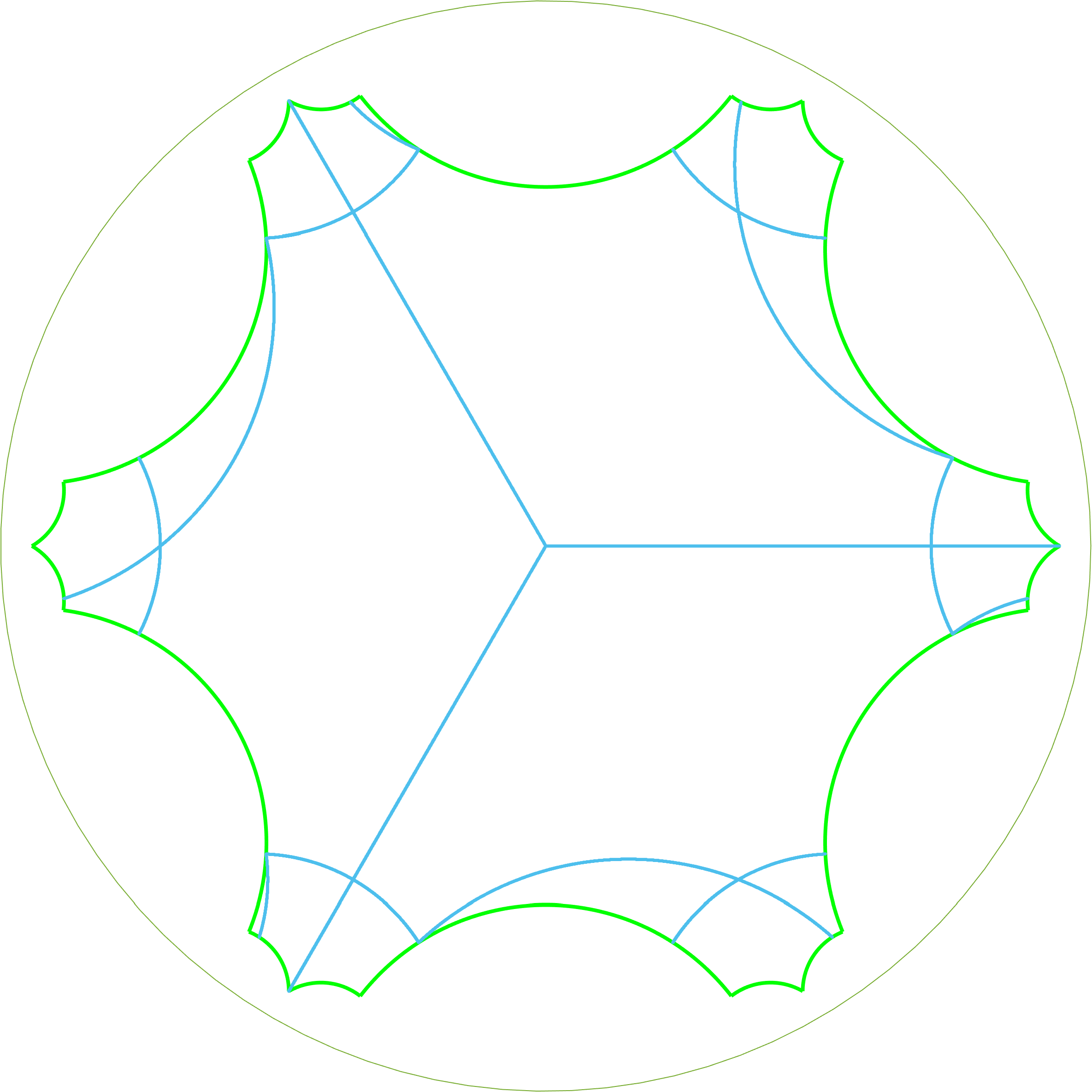}
    \caption{Tiling number $2$.}
  \end{subfigure}
  \hfill
    \begin{subfigure}[t]{0.3\textwidth}
  \centering
    \includegraphics[width=\textwidth]{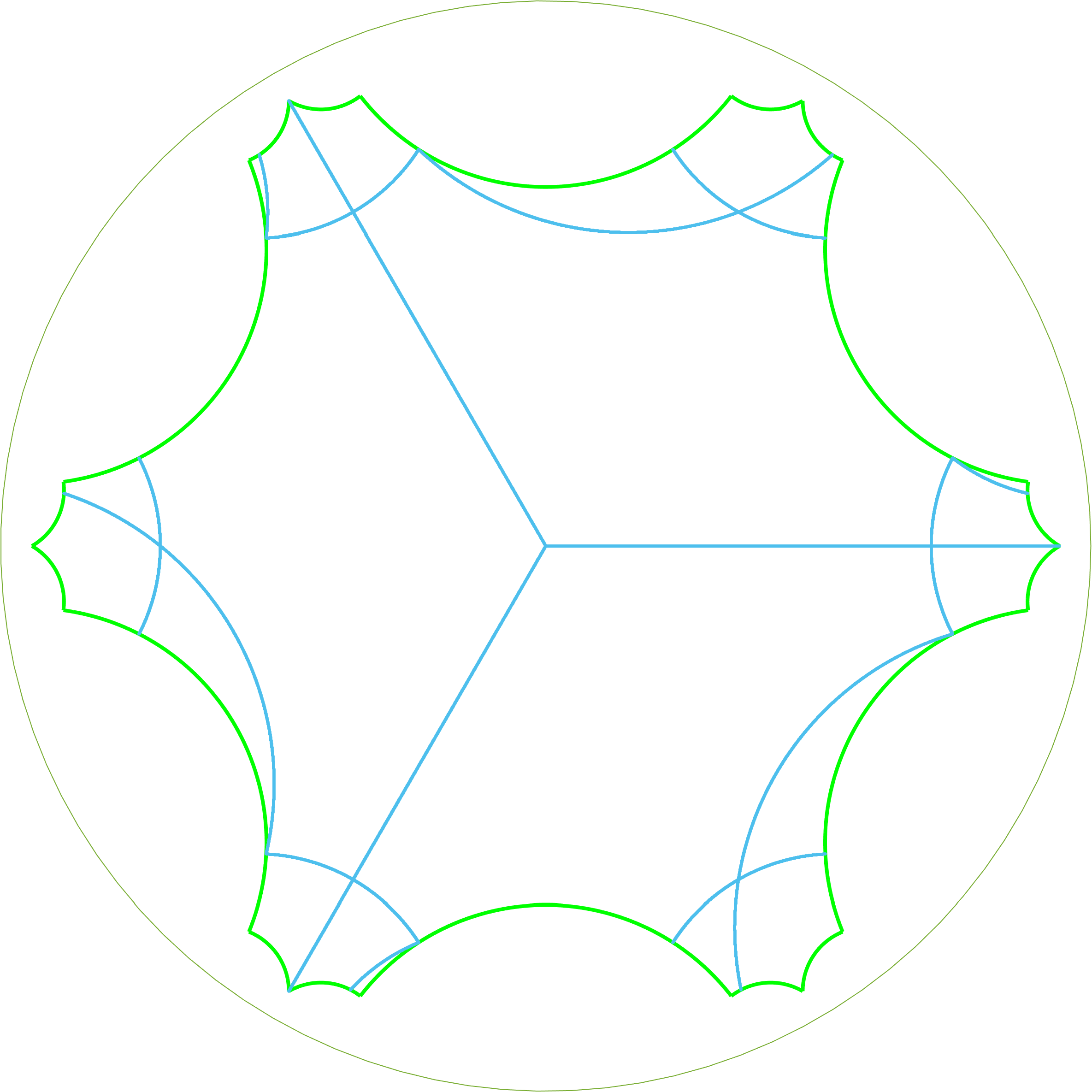}
    \caption{Tiling number $3$.}
  \end{subfigure}
  \\
  \begin{subfigure}[t]{0.3\textwidth}
  \centering
    \includegraphics[width=\textwidth]{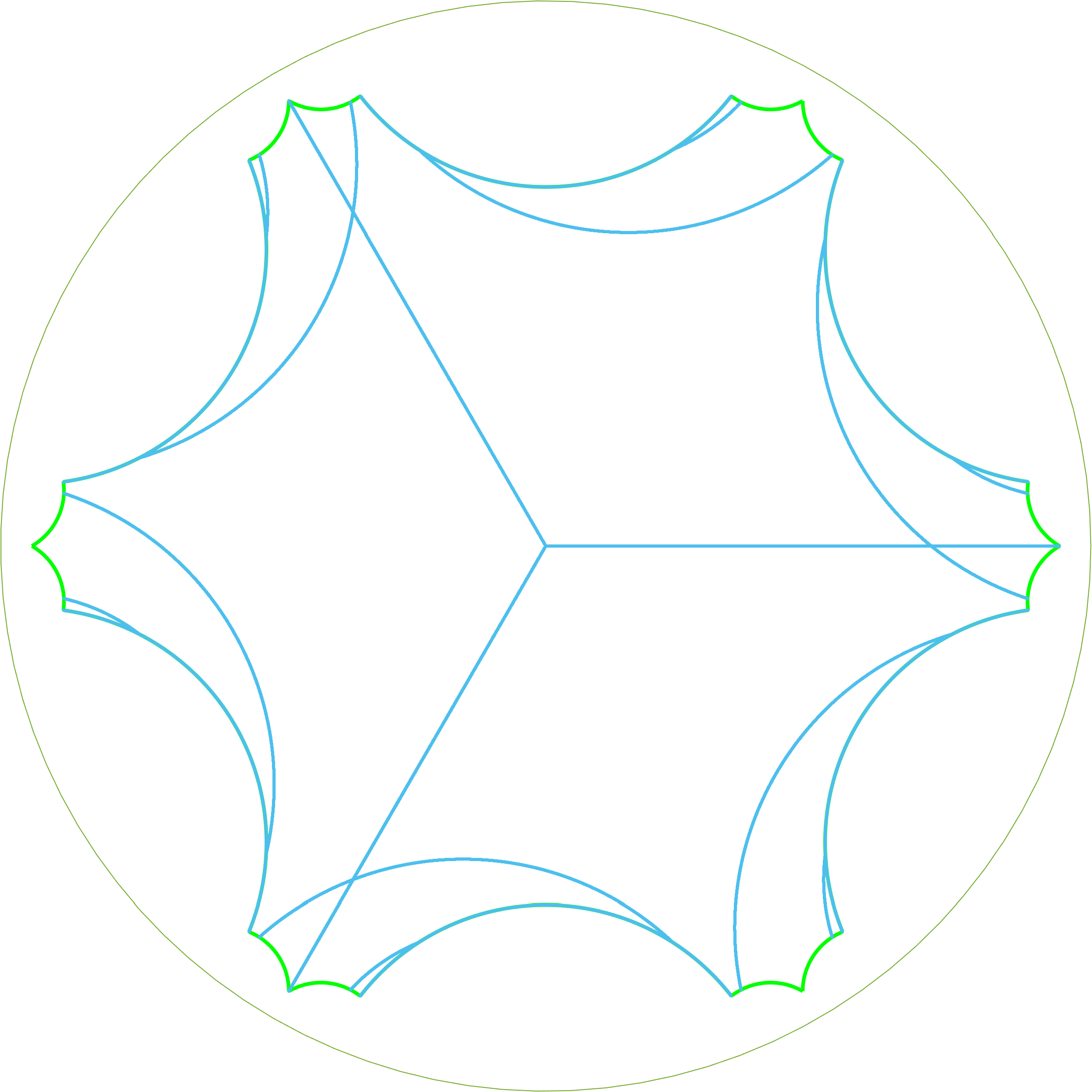}
    \caption{Tiling number $4$.}
  \end{subfigure}
  \hfill
  \begin{subfigure}[t]{0.3\textwidth}
  \centering
    \includegraphics[width=\textwidth]{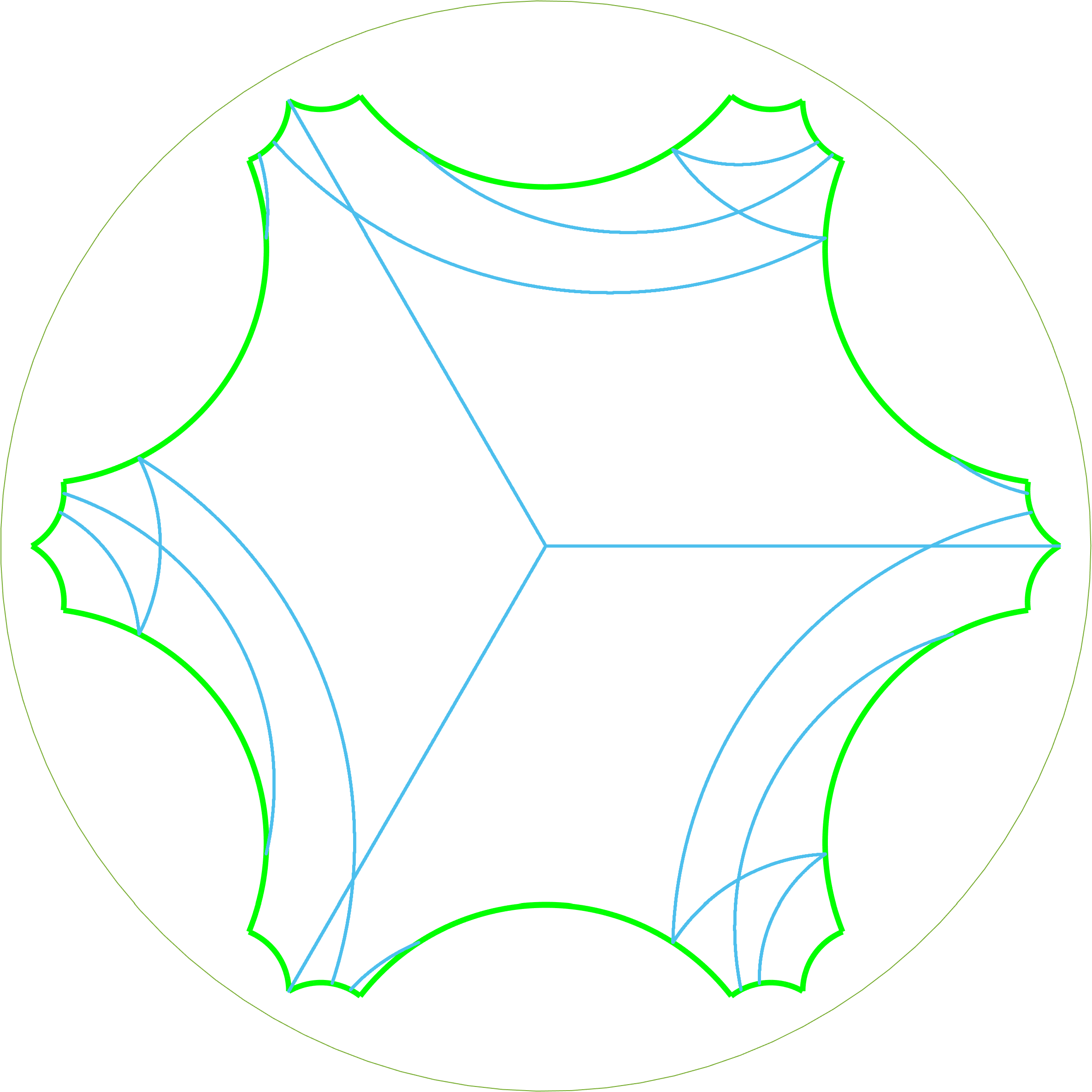}
    \caption{Tiling number $100$.}
  \end{subfigure}
  \hfill
  \begin{subfigure}[t]{0.3\textwidth}
  \centering
    \includegraphics[width=\textwidth]{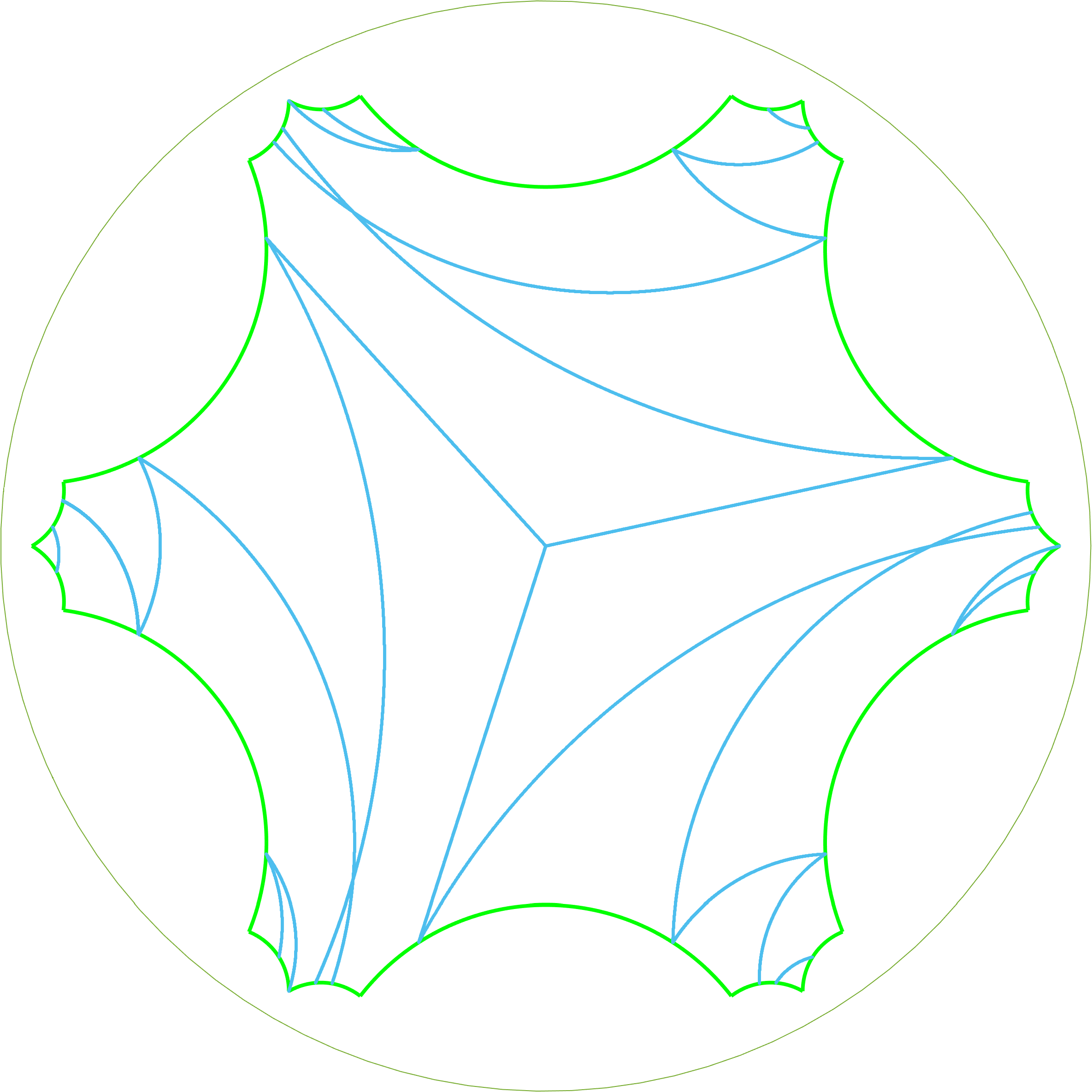}
    \caption{Tiling number $500$.}
  \end{subfigure}
    \caption{Tilings commensurate with the H-surface, with symmetry group \orb{22223}. Each is in the same combinatorial class, but isotopically distinct.}\label{fig:fds22223H_99}
\end{figure}

We conclude our presentation of enumerations of isotopy classes of tilings by presenting an example illustrating a general feature of our enumerative process. The symmetry group $22223$ in table~\ref{table:stellatemcgs} is also commensurate with another class of TPMS, namely, the associated family of the hexagonal H-surface~\cite{RobinsHsurface}. Figure~\ref{fig:fds22223H_99} shows a sample of tilings of the fundamental $18$-gon that gives rise to the H-surface with appropriate edge identifications~\cite{RobinsHsurface}, each with symmetry group $22223$ and identical D-symbol. In our enumerative process, the only difference between treating the group $22223$ of table~\ref{table:stellatemcgs} and that of the H-surface is the set of starting generators (and the shape of the fundamental polygon for the hyperbolic covering surface).
\section{Summary}
In this paper, we described a method to unambiguously produce all isotopy classes of coloured tilings of a hyperbolic surface starting from a set of geometric generators for the symmetry group of a tiling and given a construction recipe for the tiling in terms of the generators. Fundamental tile transitive tilings provide a natural starting point for an enumeration since all other classes of tilings are obtained using operations on these~\cite{BenMyf1}. The enumeration of these is achieved by first enumerating the D-symbols that encode the combinatorial types of such tilings. These in turn yield decorations of the orbifold associated to the hyperbolic symmetry group in terms of representations of the symmetry group $\Gamma$ in $\Iso(\mathbb{H}^2)$, which are induced by a set of generators for $\Gamma$ in $\Iso(\mathbb{H}^2).$ The action of the MCG on sets of generators of $\Gamma$ yields all possible sets of generators satisfying the same relation. Using the construction of the tiling in terms of generators, all isotopy classes of tilings are produced from the action of the MCG on the sets of generators. The action of prominent generators of all possible MCGs of hyperbolic symmetry groups was derived under the correspondence of (type-preserving) outer automorphisms of the symmetry group and the MCG. 
 
Subsequently, the results were used in section~\ref{sec:steltileex} to illustrate the general approach to a systematic enumeration of isotopy classes of tilings by producing an enumeration of fundamental tile-$1$-transitive tilings related to symmetry groups generated by rotational symmetries. The importance of the class of examples comes from the fact that the associated tilings all fit onto the genus-$3$ hyperbolic surface that constitutes the (uniformized) parts of the diamond, primitive, and gyroid minimal surfaces in a unit cell, with hyperbolic structure induced by the hyperbolic orbifold group $\star 246$. 




\section*{Acknowledgments}
We  thank  Vanessa  Robins,  Stephen  Hyde,  and  Stuart  Ramsden  of the  Australian  National  University  for  detailed  discussion  and  guidance.  This  research  was funded  by  the  Emmy  Noether  Programme  of  the  Deutsche  Forschungsgemeinschaft.  B.K. was  supported  by  the  Deutscher  Akademischer  Austauschdienst  for  a  research  stay  at  the Australian  National  University. 

\bibliographystyle{siamplain}
\bibliography{bibliography}
\end{document}